%% file: Ward2026_combined.tex
\RequirePackage{xr-hyper}
\documentclass[pdflatex,sn-mathphys-num]{sn-jnl}

\input{Shared/shared__preamble.tex}
\input{Article/article__preamble}
\newif\ifsuplementaryincluded
\suplementaryincludedtrue

\begin{document}

\maketitle

\addtocontents{toc}{\protect\setcounter{tocdepth}{-5}}

\input{Article/article__sections}

\newpage
\null\vskip21pt
{\parindent0pt\Titlefont Supplementary Material:\hb\papertitle\par}
\vskip24pt

\tableofcontents
\addtocontents{toc}{\protect\setcounter{tocdepth}{2}}

\newpage
\activatesmnumbering

\input{Supplementary/supplementary__sections}


\input{Article/article_acknowledgements}
\bibliography{references}

\end{document}

%% file: Shared/shared__preamble.tex
\input{Shared/shared_packages}
\input{Shared/shared_setup}

\input{Shared/shared_commands}
\input{Shared/shared_title}

\input{Shared/shared_authors}
\input{Shared/shared_abstract}
\input{Shared/shared_keywords}

%% file: Shared/shared_packages.tex
\usepackage{amsmath}  
\usepackage{amsfonts}  
\usepackage{amssymb}  
\usepackage{amsopn}   
\usepackage{accents}  
\usepackage{mathtools}  

\usepackage{algpseudocode}  
\usepackage{algorithm}  

\usepackage{array}  
\usepackage{multirow}  
\usepackage{booktabs}  
\let\cline\cmidrule

\usepackage{enumitem}  
\usepackage{cancel}  
\usepackage{url}  

\usepackage{graphicx} 
\usepackage{comment}  

\let\orcidlogo\undefined

\usepackage{orcidlink}

\usepackage{xcolor}
\usepackage{ifpdf}
\usepackage[depth=2]{bookmark}

\usepackage{cleveref}

%% file: Shared/shared_setup.tex
\definecolor{dark_red}{rgb}{0.6,0,0}
\definecolor{dark_green}{rgb}{0,0.4,0}
\definecolor{dark_blue}{rgb}{0,0,0.5}
\definecolor{colour_royal_purple}{rgb}{0.47, 0.32, 0.66}

\definecolor{c1}{HTML}{954641}
\definecolor{c2}{HTML}{40518C}
\definecolor{c3}{HTML}{40844D}
\definecolor{c4}{HTML}{8A43A4}
\definecolor{c5}{HTML}{949494}
\definecolor{c6}{HTML}{B5984E}

\hypersetup{hidelinks}

\hypersetup{bookmarksnumbered=true}

\theoremstyle{thmstyleone}
\newtheorem{theorem}{Theorem}[section]
\newtheorem{proposition}[theorem]{Proposition}
\newtheorem{lemma}[theorem]{Lemma}

\theoremstyle{thmstylethree}
\newtheorem{definition}[theorem]{Definition}
\theoremstyle{thmstyletwo}
\newtheorem{remark}[theorem]{Remark}
\newtheorem{example}[theorem]{Example}

\theoremstyle{thmstyleone}
\numberwithin{equation}{section}

\crefname{hypothesis}{Hypothesis}{Hypotheses}
\crefname{example}{Example}{Examples}

\makeatletter
\renewcommand*\l@section[2]{%
  \ifnum \c@tocdepth >\z@
    \addpenalty\@secpenalty
    \addvspace{1.0em \@plus\p@}%
    \setlength\@tempdima{3.2em}%
    \begingroup
      \parindent \z@ \rightskip \@pnumwidth
      \parfillskip -\@pnumwidth
      \leavevmode \bfseries
      \advance\leftskip\@tempdima
      \hskip -\leftskip
      #1\nobreak\hfil
      \nobreak\hb@xt@\@pnumwidth{\hss #2\kern-\p@\kern\p@}\par
    \endgroup
  \fi}
\renewcommand*\l@subsection{\@dottedtocline{2}{3.2em}{3.4em}}
\renewcommand*\l@subsubsection{\@dottedtocline{3}{6.6em}{4.2em}}
\makeatother

\newcommand{\activatesmnumbering}{%
  \setcounter{section}{0}%
  \setcounter{subsection}{0}%
  \setcounter{subsubsection}{0}%
  \setcounter{equation}{0}%
  \setcounter{figure}{0}%
  \setcounter{table}{0}%
  \setcounter{theorem}{0}%
  \setcounter{algorithm}{0}%
  \setcounter{footnote}{0}%
  \renewcommand{\thesection}{SM\arabic{section}}%
  \renewcommand{\thefigure}{SM\arabic{figure}}%
  \renewcommand{\thetable}{SM\arabic{table}}%
  \renewcommand{\thealgorithm}{SM\arabic{algorithm}}%
  \renewcommand{\theHsection}{SM.\arabic{section}}%
  \renewcommand{\theHsubsection}{\theHsection.\arabic{subsection}}%
  \renewcommand{\theHsubsubsection}{\theHsubsection.\arabic{subsubsection}}%
  \renewcommand{\theHequation}{SM.\arabic{section}.\arabic{equation}}%
  \renewcommand{\theHfigure}{SM.\arabic{figure}}%
  \renewcommand{\theHtable}{SM.\arabic{table}}%
  \renewcommand{\theHalgorithm}{SM.\arabic{algorithm}}%
  \renewcommand{\theHtheorem}{SM.\arabic{theorem}}%
  \renewcommand{\theHfootnote}{SM.\arabic{footnote}}%
}

\DeclareMathOperator{\sign}{sign}
\DeclareMathOperator*{\argmin}{argmin}

\DeclareMathOperator*{\hingeloss}{HingeLoss}
\DeclareMathOperator*{\minimise}{minimise}
\DeclareMathOperator*{\maximise}{maximise}
\DeclareMathOperator*{\FischerBurmeister}{fb}
\DeclareMathOperator*{\MStationaryFunc}{msf}
\DeclareMathOperator*{\subjectto}{subject~to}

%% file: Shared/shared_commands.tex
\newcommand*{\R}{\mathbb{R}}
\newcommand*{\N}{\mathbb{N}}

\newcommand*{\Q}{\mathcal{Q}}
\newcommand*{\Y}{\mathcal{Y}}
\newcommand*{\id}{\mathcal{I}}

\newcommand*{\aaa}{{\boldsymbol{a}}}
\newcommand*{\ddd}{{\boldsymbol{d}}}
\newcommand*{\zzz}{{\boldsymbol{z}}}
\newcommand*{\www}{{\boldsymbol{w}}}
\newcommand*{\xxx}{{\boldsymbol{x}}}
\newcommand*{\xibold}{{\boldsymbol{\xi}}}

\newcommand*{\zone}{{{\boldsymbol{z}}_{\boldsymbol{1}}}}{}
\newcommand*{\ztwo}{{{\boldsymbol{z}}_{\boldsymbol{2}}}}{}

\newcommand{\penaltyparam}{\pi}
\newcommand{\relaxparam}{\tau}

\newcommand*{\vvhead}[2]{\rotatebox[origin=c]{90}{\hspace{.5em}\shortstack[c]{#1\\#2}\hspace{.5em}}}
\newcommand*{\vhead}[1]{\rotatebox[origin=c]{90}{ \hspace{.5em}\shortstack[c]{#1}\hspace{.5em}}}

\makeatletter
\def\mlt@hassub#1_#2\mlt@end{\def\mlt@rest{#2}}
\newcommand{\multiplier}[4][]{%
  \mlt@hassub#2_\mlt@end
  \ifx\mlt@rest\empty
    \ifx\relax#1\relax
      {\color{#4}\boldsymbol{#3}#2}%
    \else
      {\color{#4}\hat{\boldsymbol{#3}}#2}%
    \fi
  \else
    \ifx\relax#1\relax
      {\color{#4}#3#2}%
    \else
      {\color{#4}\hat{#3}#2}%
    \fi
  \fi
}
\makeatother
\newcommand{\multg}[2][]{\multiplier[#1]{#2}{\lambda}{dark_red}}
\newcommand{\multh}[2][]{\multiplier[#1]{#2}{\mu}{dark_red}}
\newcommand{\multG}[2][]{\multiplier[#1]{#2}{\nu}{dark_red}}
\newcommand{\multH}[2][]{\multiplier[#1]{#2}{\xi}{dark_red}}
\newcommand{\multGH}[2][]{\multiplier[#1]{#2}{\omega}{dark_red}}
\newcommand*{\ablue}[2][]{\multiplier[#1]{#2}{\alpha}{dark_green}}
\newcommand{\ubar}[1]{\underaccent{\bar}{#1}}
\newcommand*{\mupper}[2][]{\multiplier[#1]{#2}{\bar{v}}{dark_green}}
\newcommand*{\mlower}[2][]{\multiplier[#1]{#2}{\ubar{v}}{dark_green}}
\newcommand*{\uu}[1]{\textcolor{dark_green}{u{#1}}}
\newcommand*{\etablue}[2][]{\multiplier[#1]{#2}{\eta}{dark_green}}
\newcommand*{\multip}[1]{\multiplier{#1}{\Lambda}{dark_red}}

\newcommand*{\fupper}{f_{\text{up}}}
\newcommand*{\gupper}{g_{\text{up}}}
\newcommand*{\hupper}{h_{\text{up}}}
\newcommand*{\flower}{f_{\text{low}}}
\newcommand*{\glower}{g_{\text{low}}}
\newcommand*{\hlower}{h_{\text{low}}}

\newcommand{\link}[2]{\hyperref[#1]{#2}}
\makeatletter
 \newcommand{\target}[1]{\Hy@raisedlink{\hypertarget{#1}{}}}
\makeatother

\newcommand{\MODEL}{\eqref{eq:final-model}}
\newcommand{\NLP}{\hyperlink{def:NLP}{NLP}}
\newcommand{\MPCC}{\hyperlink{def:MPCC}{MPCC}}
\newcommand{\LICQ}{\hyperlink{def:LICQ}{LICQ}}
\newcommand{\MFCQ}{\hyperlink{def:MFCQ}{MFCQ}}
\newcommand{\MPCCLICQ}{\link{def:MPCC-LICQ}{MPCC-LICQ}}
\newcommand{\MPCCMFCQT}{\link{def:MPCC-MFCQ-T}{MPCC-MFCQ-T}}
\newcommand{\MPCCMFCQR}{\link{def:MPCC-MFCQ-R}{MPCC-MFCQ-R}}
\newcommand{\MPCCMFCQC}{\link{def:MPCC-MFCQ-C}{MPCC-MFCQ-C}}
\newcommand{\penalisation}[1]{(\link{eq:penalisation}{$\mathcal{P}_{#1}$})}
\newcommand{\relaxation}[1]{(\link{eq:relaxation}{$\mathcal{R}_{#1}$})}
\newcommand{\RNLP}[1]{(\link{eq:RNLP}{RNLP$_{#1}(\bar{\zzz})$})}

\newcommand*{\veczero}{\mathbf{0}}
\newcommand*{\vecone}{\mathbf{1}}

\newcommand{\norm}[1]{\left\lVert#1\right\rVert}

\newcommand{\derivative}[1]{\nabla_{\hspace{-.2em}#1}{}}

%% file: Shared/shared_title.tex
\newcommand{\papertitle}{Mathematical programs with complementarity constraints and application to hyperparameter tuning}
\newcommand{\papershorttitle}{MPCCs and application to hyperparameter tuning}
\title[\papershorttitle]{\papertitle}

%% file: Shared/shared_authors.tex
\author*[1]{%
    \fnm{Samuel J.}
    \sur{Ward}
    \orcidlink{0009-0001-3084-7099}
}
\email{s.ward@soton.ac.uk}

\author[1]{%
    \fnm{Alain} 
    \sur{Zemkoho}
    \orcidlink{0000-0003-1265-4178}
}
\email{a.b.zemkoho@soton.ac.uk}

\author[1]{%
    \fnm{Selin Damla} 
    \sur{Ahipa\c{s}ao\u{g}lu}
    \orcidlink{0000-0003-1371-315X}
}
\email{s.d.ahipasaoglu@soton.ac.uk}

\affil[1]{%
    \orgdiv{School of Mathematical Sciences},
    \orgname{University of Southampton},
    \orgaddress{%
        \city{Southampton},
        \postcode{SO17~1BJ},
        \country{United Kingdom}%
    }%
}

%% file: Shared/shared_abstract.tex
\newcommand{\paperabstract}{%
We consider the Mathematical Program with Complementarity Constraints (MPCC). One of the main challenges in solving this problem is the systematic failure of standard constraint qualifications. Carefully accounting for the combinatorial nature of the complementarity constraints, tractable versions of the Mangasarian--Fromovitz Constraint Qualification (MFCQ) have been designed and widely studied in the literature.  This paper looks closely at two such MPCC-MFCQs and their influence on MPCC algorithms. As a key contribution, we
introduce a third MPCC-MFCQ, lying between these two, under which we prove the convergence of the inexact sequential penalisation and Scholtes relaxation algorithms.
This significantly weakens the assumptions currently used in the literature.
We then form the problem of tuning hyperparameters of a nonlinear Support Vector Machine (SVM), a fundamental machine learning problem for classification, as an MPCC. For this application, we establish that the aforementioned
new qualification holds under
very mild restrictions. Moreover, we program robust implementations and carry out comprehensive numerical experiments on real-world data sets, where we show that the sequential penalisation method applied to the MPCC formulation for tuning SVM hyperparameters can outperform both the Scholtes relaxation technique and the derivative-free methods from the machine learning literature.%
}

%% file: Shared/shared_keywords.tex
\newcommand{\paperkeywords}{%
Mathematical program with complementarity constraints,
Hyperparameter optimisation,
Bilevel optimisation,
Penalisation,
Relaxation,
Support vector machines%
}

\newcommand{\papermsc}{%
    49M37,  
    65K05,  
    65K15,  
    90C30,  
    90C33
}

%% file: Article/article__preamble.tex
\ifpdf
\hypersetup{
  pdftitle={Mathematical programs with complementarity constraints and application to hyperparameter tuning},
  pdfauthor={Samuel J. Ward, Alain Zemkoho, Selin Damla Ahipasaoglu},
  pdfsubject={Mathematics}
}
\fi

\abstract{\paperabstract}

\keywords{\paperkeywords}

\pacs[MSC Classification]{\papermsc}

%% file: Article/article__sections.tex
\input{Article/article_10_introduction}

\input{Article/article_20_preliminaries}

\input{Article/article_21_constraint_quali}

\input{Article/article_22_stationarity}
\input{Article/article_30_solution_methods}
\input{Article/article_31_penalisation}

\input{Article/article_32_relaxation}
\input{Article/article_40_application}
\input{Article/article_41_model}
\input{Article/article_42_properties}

\input{Article/article_50_experiments}
\input{Article/article_51_implementation}

\input{Article/article_52_behaviour}
\input{Article/article_53_performance}

\input{Article/article_60_conclusion}

\input{Article/article_declarations}

%% file: Article/article_10_introduction.tex
\section{Introduction}
\label{section:introduction}
% ...
The Mathematical Program with Complementarity Constraints~(\MPCC{}) is now a well-established optimisation problem, with a wide range of applications including in economics~\cite{Kenneth1954}, chemical engineering~\cite{Baumrucker2008}, energy distribution~\cite{Gabriel2010}, transportation~\cite{Lawphongpanich2004}, and machine learning (see Section~\ref{section:Application} for an important problem class, as well as some relevant references).

Most standard optimisation algorithms are ineffective at finding solutions of MPCCs, in part due to the systematic failure of classical Constraint Qualifications~(CQs) and the combinatorial nature of the problem's feasible set~\cite[Chapter 3]{Luo1996}.
To address the issue with standard CQs, many alternative tractable ones tailored to the complementarity constraints have been developed in the literature (see, e.g.,~\cite{Flegel2005,Pang1999,Scheel2000,schwartz2011mathematical,Ye2005}).
Considering the importance of the Mangasarian--Fromovitz Constraint Qualification~(\MFCQ{}) in optimisation, it has been particularly at the centre of attention. Two categories of MFCQ-type CQs for the MPCC have emerged in the literature; one is based on the decomposition of the complementarity constraints, especially the constraints associated with the biactive index sets, to build suitable versions of MFCQ that can hold for MPCCs~\cite{Flegel2005,Scheel2000,Scholtes2001}.
A second category, constructed from the lens of variational analysis, consists of reformulating the complementarity constraints either with non-smooth operator functions or sets that lead to normal cone operations that restore the fulfilment of such CQs; see, e.g.,~\cite{DempeZemkoho2012,gfrerer2017new,Ye2005}.

In the literature, the \MPCCMFCQT{}, where we use the label ``T'' to represent the fact that it is built from a tightened decomposition of the complementarity constraints, has emerged as one of the most used CQs for MPCCs; see, e.g.,~\cite{Flegel2005,Kim2020,Scheel2000,schwartz2011mathematical}.  However, looking at the \MPCCMFCQT{} very closely, one can observe that it boils down to the \MPCCLICQ{}, the well-known tailored linear independence constraint qualification for MPCCs, when the feasible set does not involve any inequality constraint apart from the ones in the complementarity constraints; cf. Remark~\ref{remark:MPCCMFCQT_boils_down_to_MPCCLICQ} in Section~\ref{section:preliminaries}. For this reason, and considering the fact that the MPCC-LICQ is very strong, the \MPCCMFCQT{} fails for many MPCC problems (see Example~\ref{eq:Example01} and Proposition~\ref{proposition:fail_mpcc_licq}).

Alternatively, there is another version of MPCC-tailored MFCQ, studied in~\cite{Qiu2015,Ralph2004}, based on a relaxation of the complementarity constraints associated with the biactive index set. Hence, we denote it as \MPCCMFCQR{} with the label ``R'' representing the corresponding relaxation framework. 
We show in this paper that the \MPCCMFCQR{} is strictly weaker than the \MPCCMFCQT{}.
We then introduce a third qualification, the \MPCCMFCQC{}, which we show lies strictly between these two and which we denote with a ``C'' because it is associated with the Lagrange multiplier pattern of C-stationarity.
It is related to conditions in~\cite{Scheel2000,Ye2005} but, to the best of our knowledge, has not been stated explicitly before.

Based on this observation, we revisit two classical algorithms for MPCCs that rely on MPCC-type CQs for their convergence results, namely, the partial penalisation and Scholtes relaxation methods, which have been widely studied and used in the literature (see, e.g.,~\cite{Anitescu2000,Leyffer2006interior,Ralph2004} and~\cite{Demiguel2005,Hoheisel2013,Lin2005,Scholtes2001,Steffensen2010}, respectively). 

The partial penalisation method consists of moving the product of the complementarity constraint functions from the feasible set to the objective.  Two options have been pursued for this approach: (i) the \textit{exact} penalisation, where the penalty parameter is fixed~\cite{Ralph2004}, and (ii) the \textit{sequential} version, where the penalty parameter increases sequentially throughout the algorithm~\cite{Leyffer2006interior}.  As is well known in constrained optimisation, a big challenge with exact penalisation is identifying a suitable penalty parameter which is highly impactful on the numerical performance. Crucially, for sequential penalisation algorithms, the selection of the initial penalisation parameter is not too critical, as theoretically, under suitable assumptions, the algorithm gets to a point of stability where a solution is found.  

Considering this practical advantage of the sequential approach, we focus on it in the analysis in this paper. The article~\cite{Leyffer2006interior} shows that the sequential partial penalisation algorithm, in conjunction with an interior-point method for the subproblems, converges to an S-stationary point, under the \MPCCLICQ{}, as a key assumption.  In this paper,  we establish the convergence of the sequential partial penalisation method to C-stationary points under the
\MPCCMFCQC{}, and to S-stationary points under the \MPCCMFCQR{}
together with an additional assumption at the biactive indices; cf. Theorem~\ref{theorem:penalisation_convergence}. Not only
are the \MPCCMFCQC{} and \MPCCMFCQR{} far weaker than the \MPCCLICQ{}, but the subproblems can be solved with any method, and not necessarily with an interior-point method, as required in~\cite{Leyffer2006interior}.

As for the Scholtes relaxation algorithm, its convergence has been established under the \MPCCMFCQT{}; see, e.g.,~\cite{schwartz2011mathematical}. Here, we show that the
weaker \MPCCMFCQC{} is sufficient to establish this convergence%
, and that S-stationarity follows under the \MPCCMFCQR{} together with an additional assumption at the biactive indices; cf. Theorem~\ref{theorem:relaxation_convergence}.   Moreover, to make the method more practically relevant, in our proof, we only require the computation of approximate Karush--Kuhn--Tucker (KKT) points for the relaxed problem at each iteration. This is a significant departure from the existing convergence results, where it is usually assumed that the KKT points of the relaxed problem are computed exactly, something that is typically not possible in a practical implementation.
A notable exception is~\cite{Kanzow2015}, which allows errors in all the KKT conditions of the relaxed problem but requires the \MPCCMFCQT{} and a tolerance that tends to zero faster than the relaxation parameter; see Table~\ref{table:comparison} for a comparison of the assumptions.

Subsequently, we explore the application of our aforementioned results for general MPCCs to the problem of computing optimal hyperparameters for a Support Vector Machine (SVM).  
Note that the SVM represents one of the most potent tools for classification in machine learning. 
Since its introduction in 1995 by Cortes and Vapnik~\cite{Cortes1995,Vapnik1999}, it continues to be prolifically used to solve a multitude of practical problems, including facial authentication~\cite{Tefas2001}, protein sequencing~\cite{Rahman2018} and speech recognition~\cite{Smith2001}; it also sets a competitive benchmark in all classification problems.

One of the biggest challenges with SVMs is choosing their hyperparameters. 
Unlike model parameters, which are learnt from patterns in the training data, hyperparameters define the learning process itself and so must be considered separately. 
The two most common hyperparameters are the \textit{regularisation parameter}, which balances the empirical loss against the complexity of the model, and the \textit{choice of kernel}.
There is no one model that will be best suited to all scenarios~\cite{Wolpert1996}. 

Traditionally, in the machine learning literature, derivative-free search algorithms such as grid search and Bayesian optimisation have been used to select the hyperparameters that provide the best validation accuracy~\cite{Wainer2021}.
However, these methods all assume that training is a black box -- that a new model must, each time, be trained in its entirety before being evaluated. 
They are general methods that neglect the specific structure of the SVM problem and scale very poorly with multiple hyperparameters.

More recently, a new approach to hyperparameter selection using bilevel optimisation has been proposed~\cite{Franceschi2018,Li2022unified,Okuno2021}. 
It enables the simultaneous training of the model while hyperparameter tuning is being conducted.
The bilevel optimisation technique has been applied thoroughly and with great success to linear models. For example, the series of papers~\cite{Moore2011,Moore2010,Moore2009} proposes some non-smooth optimisation methods for selecting linear support vector regression and classification hyperparameters. 
Other authors have created a similar framework but propose different approaches to solving the bilevel program; e.g., a value function-based descent, global relaxation, and stochastic gradient descent methods have been proposed in~\cite{Gao2022},~\cite{Li2022linear} and~\cite{Jiang2020}, respectively. 
The aforementioned works are on linear models; thus, enabling some mathematical simplicity, but such a framework is less useful in practice, as real-world data sets are rarely linearly separable.

Of particular interest to us is the sequence of papers~\cite{Bennett2006,Kunapuli2008b,Kunapuli2008a}, as they suggest an extension of the bilevel optimisation framework to nonlinear kernel models, while proposing an MPCC reformulation for the corresponding problem.  
However, they do not provide much theoretical insight into the problem.
The paper~\cite{Coniglio2022} begins addressing this gap.
Our paper expands on the latter work by first rigorously constructing the MPCC model of tuning the hyperparameters of nonlinear kernel SVMs.
We then prove for the first time that the
\MPCCMFCQC{}, and hence the \MPCCMFCQR{}%
, automatically holds at any feasible point, provided that the corresponding pair of hyperparameters is positive%
, at least one training point of each fold is a free support vector, and the training inputs of each fold are distinct; cf. Proposition~\ref{proposition:satisfy_mpccmfcqr}.
Note that this positiveness is not a strong requirement, as for SVMs, those values cannot be zero in practice%
; the two other restrictions are also mild~\cite{KeerthiLin2003asymptotic}.
As a consequence of this result, the sequential penalisation and Scholtes relaxation algorithms, discussed above, converge automatically, almost for free,
to C-stationary points for the MPCC model.
It is important to note that we also show that for the MPCC hyperparameter optimisation model, the \MPCCMFCQT{} fails at any of its feasible points with more than two biactive indices; see Proposition~\ref{proposition:fail_mpcc_licq}.

As another key contribution of this paper, we conduct a robust implementation of the partial penalisation and Scholtes relaxation methods and compare them to each other, and also to classical hyperparameter tuning algorithms in the machine learning literature. 
First, we observe that the sequential penalisation method, on average, not only runs faster but also converges to a stationary point with a lower objective function value than the Scholtes relaxation algorithm.
Moreover, we show that the sequential partial penalisation technique can outperform the traditional derivative-free methods such as grid, random, Bayesian and pattern search. 
Note that for the experiments, we use 18 well-known classification datasets from the literature, where the median problem size is 1773 variables, 1015 constraints (excluding complementarity constraints) and 1010 pairs of complementarity constraints. 
To the best of our knowledge, these problems are far larger than any others that have been used for experiments in the context of non-convex MPCCs. 
For illustration, observe that for the MacMPEC collection~\cite{MacMPEC}, which is classical for MPCC numerical experimentation, the median problem size is 49 variables, 20 constraints (excluding complementarity constraints) and 12 pairs of complementarity constraints.

For the remainder of the paper, note that Section~\ref{section:preliminaries} covers preliminaries with particular attention on the difference between the two MPCC-tailored CQs that are a key focus in this paper; i.e., the \MPCCMFCQT{} and \MPCCMFCQR{}%
, as well as the new \MPCCMFCQC{} lying between them.
In Section~\ref{section:solution_methods}, we present the sequential inexact penalisation and relaxation methods and prove new convergence results which allow these methods to be applicable to a broader range of problems.
In Section~\ref{section:Application}, we construct the MPCC optimisation-based hyperparameter tuning model for SVMs and prove that it satisfies
\MPCCMFCQC{}, and hence \MPCCMFCQR{}%
, but not \MPCCMFCQT{}.
Finally, Section~\ref{section:experiments} presents our numerical experiments and a demonstration that the MPCC model solved with sequential penalisation can outperform traditional derivative-free methods from the machine learning literature. 
Many technical details, which could further facilitate the comprehension and implementation of some aspects of the paper, are provided in 
\ifsuplementaryincluded
the supplementary materials at the end of this document\nocite{Ward2025}.
\else
online supplementary materials~\cite{Ward2025}.
\fi

%% file: Article/article_20_preliminaries.tex
\section{Preliminaries on MPCCs}
\label{section:preliminaries}
% ...
\MakeLinkTarget*{def:NLP}{}%
For a number of concepts used throughout the paper, we refer to a general Nonlinear Program~(\NLP{}) defined, for given dimensions $n,p,q\in \N$, as
\begin{equation*}
\begin{aligned}
\minimise_{\zzz\in \R^n} \quad
&f(\zzz)&
\\
\subjectto \quad
& g_i(\zzz) \geq 0&
\quad\text{for } i=1,...,p,
\\
&h_i(\zzz) = 0&
\quad\text{for } i=1,...,q,
\end{aligned}
\end{equation*}
where the functions $f:\R^n\rightarrow\R, \ g:\R^n\rightarrow\R^p, \ h:\R^n\rightarrow\R^q$ representing the objective, inequality, and equality constraint functions, respectively, are continuously differentiable. 
\MakeLinkTarget*{def:MPCC}{}%
However, the main focus of our analysis is the Mathematical Program with
Complementarity Constraints~(\MPCC{}),
defined for given dimensions $n, p, q, r\in \N$, as 
\begin{equation*}
\begin{aligned}
\minimise_{\zzz\in \R^n} \quad
&f(\zzz)&
\\
\subjectto \quad
& g_i(\zzz) \geq 0& 
\quad\text{for } i=1,...,p,
\\
&h_i(\zzz) = 0&
\quad\text{for } i=1,...,q,
\\
&G_i(\zzz) \geq 0,\ H_i(\zzz) \geq 0,\ G_i(\zzz)H_i(\zzz)=0 &
\quad\text{for } i=1,...,r,
\end{aligned}
\end{equation*}
where the functions
$f:\R^n \rightarrow \R$,
$g:\R^n \rightarrow \R^p$,
$h:\R^n \rightarrow \R^q$ and
$G, H:\R^n \rightarrow \R^r$ are continuously differentiable. 
$G$ and $H$ are referred to as complementarity  constraint functions and $G_i(\zzz)H_i(\zzz)$ as the product term.
The last line of problem~(\MPCC{}) can be written compactly as $0\leq G(\zzz) \bot H(\zzz) \geq 0$.

%% file: Article/article_21_constraint_quali.tex
\subsection{Constraint qualifications}
\label{section:constraint_qualification}
% ...
In continuous constrained optimisation, constraint qualifications (CQs) are usually required to construct necessary optimality conditions and prove the convergence of numerical algorithms. 
To introduce a CQ for problem (\NLP{}), we define the \emph{index sets}
\begin{equation*}
    I^g(\bar{\zzz}) \coloneqq \left\{i \in \{1,...,p\}:\; g_i(\bar{\zzz}) = 0 \right\} \, \mbox{ and }\; I^h \coloneqq \left\{1,...,q\right\}.
\end{equation*}
\MakeLinkTarget*{def:LICQ}{}%
\noindent Based on this notation, the \emph{Linear Independence Constraint Qualification} (\LICQ{}) is satisfied at a feasible point $\bar{\zzz}\in\R^n$ of the~(\NLP{}) iff the family of gradients
\begin{align*}
    &\left\{
        \nabla h_i(\bar{\zzz}): i\in I^h
    \right\} \cup 
    \left\{
        \nabla g_i(\bar{\zzz}): i\in I^g(\bar{\zzz})
    \right\}
    \text{ is linearly independent,}
\end{align*}
\MakeLinkTarget*{def:MFCQ}{}%
while, the \emph{Mangasarian-Fromovitz Constraint Qualification} (\MFCQ{}) will be said to hold at a feasible point $\bar{\zzz}\in\R^n$ of the~(\NLP{}) iff
\begin{equation*}
    \begin{array}{l}
    \left\{\nabla h_i(\bar{\zzz}): i\in I^h \right\} \text{ is linearly independent and}\\
    \exists \ \ddd\in \R^n: \quad
    \begin{cases}
        \nabla h_i(\bar{\zzz}) ^{\top} \ddd = 0&  \forall  i\in I^h,\\
        \nabla g_i(\bar{\zzz}) ^{\top} \ddd > 0&  \forall  i\in I^g(\bar{\zzz}).
    \end{cases}
    \end{array}
\end{equation*}
Obviously, the \LICQ{} is stronger than the \MFCQ{}.  These CQs play a fundamental role in nonlinear constrained optimisation thanks to their versatility in deriving fundamental properties for problem \NLP{}; see the excellent tutorial in \cite{Solodov2010} for an overview of these properties.

Unfortunately, it is well known that \LICQ{} and \MFCQ{} are systematically violated at every feasible point of problem  {\MPCC{}} when viewed as a constrained problem with equality and inequality constraints in the format described in {\NLP{}} \cite{Chen1995}. To observe this, consider
an index $i=1, \dots, r$, a feasible point $\bar{\zzz}$ and the three constraints $G_i(\bar{\zzz})\geq 0$, $H_i(\bar{\zzz})\geq 0$ and $G_i(\bar{\zzz})H_i(\bar{\zzz})=0$. Assuming, without loss of generality, that $H_i(\bar{\zzz})$ is zero, the gradient of the second constraint $\nabla H_i(\bar{\zzz})$ and the gradient of the third constraint $G_i(\bar{\zzz})\nabla H_i(\bar{\zzz})$ (by the product rule) are not linearly independent.

The main challenge in addressing problem \MPCC{} comes from the combinatorial nature of the complementarity constraints. For each index $i=1,\dots,r$, either $H_i$ or $G_i$ may be positive, but not both. To further study this, we partition these cases by defining the following index sets:
\begin{equation}
\label{eq:index_sets}
\begin{aligned}
    I^{G}(\bar{\zzz})  \coloneqq & \left\{i \in \{1,...,r\}  : \quad G_i(\bar{\zzz}) = 0,  H_i(\bar{\zzz}) > 0 \right\},\\
    I^{H}(\bar{\zzz})  \coloneqq & \left\{i \in \{1,...,r\}  : \quad G_i(\bar{\zzz}) > 0,  H_i(\bar{\zzz}) = 0 \right\},\\
    I^{GH}(\bar{\zzz}) \coloneqq & \left\{i \in \{1,...,r\}  : \quad G_i(\bar{\zzz}) = 0,  H_i(\bar{\zzz}) = 0 \right\}.
\end{aligned}
\end{equation}
Based on these index sets, at a given feasible point $\bar{\zzz}$ of problem \MPCC{}, we first introduce the corresponding \emph{Tightened Nonlinear Program}~\eqref{eq:TNLP}:
\begin{equation}
\tag{TNLP\mbox{$(\bar{\zzz})$}} 
\label{eq:TNLP}
\begin{aligned}
\minimise_{\zzz\in\R^n} \quad
&f(\zzz)&
\\
\subjectto \quad
& g_i(\zzz) \geq 0& 
&\text{ for } i=1,\dots,p,&
\\
&h_i(\zzz) = 0&
&\text{ for } i=1,\dots,q,&
\\
&G_{i}(\zzz)=0,    \quad H_{i}(\zzz)\geq0 \quad &
& \text{ for } i\in I^{G}(\bar{\zzz}),\\
&G_{i}(\zzz)\geq0, \quad H_{i}(\zzz)=0    \quad &
& \text{ for } i\in I^{H}(\bar{\zzz}),\\
&G_{i}(\zzz)=0,    \quad H_{i}(\zzz)=0    \quad &
& \text{ for } i\in I^{GH}(\bar{\zzz}).
\end{aligned}
\end{equation}
It is very important to remember that these index sets depend on the point $\bar{\zzz}$, so the~\eqref{eq:TNLP} may be a different program when constructed around a different point.
The decomposition of the complementarity constraints in \MPCC{} makes the fulfilment of standard CQs more likely. Hence, problem~\eqref{eq:TNLP} is commonly used to define CQs tailored to problem \MPCC{} as follows.

\begin{definition}
\label{def:MPCC-MFCQ-T}
The \MPCCMFCQT{} holds at a feasible point $\bar{\zzz}$ of the \MPCC{} if the corresponding tightened program~\eqref{eq:TNLP} satisfies the \MFCQ{} at the same point $\bar{\zzz}$. Equivalently, the \MPCCMFCQT{} holds at a feasible point $\bar{\zzz}$ iff
\begin{equation*}
    \Big\{\nabla h_i(\bar{\zzz}):i \in I^h\Big\} \cup
    \Big\{\nabla G_i(\bar{\zzz}):i\in {I^{G}(\bar{\zzz})\cup I^{GH}(\bar{\zzz})}\Big\} \cup
    \Big\{\nabla H_i(\bar{\zzz}):i\in{I^{H}(\bar{\zzz}) \cup I^{GH}(\bar{\zzz})}\Big\}\notag
\end{equation*}
is linearly independent and
\begin{equation*}
    \exists \ \ddd\in \R^n: \quad
    \begin{cases}
        \nabla h_i(\bar{\zzz}) ^{\top} \ddd = 0 &  \forall i \in I^h,\\
        \nabla G_i(\bar{\zzz}) ^{\top} \ddd = 0 &  \forall i \in I^{G}(\bar{\zzz}) \cup I^{GH}(\bar{\zzz}),\\
        \nabla H_i(\bar{\zzz}) ^{\top} \ddd = 0 &  \forall i \in I^{H}(\bar{\zzz}) \cup I^{GH}(\bar{\zzz}),\\
        \nabla g_i(\bar{\zzz}) ^{\top} \ddd > 0 &  \forall i \in I^{g}(\bar{\zzz}).
    \end{cases}
\end{equation*}
\end{definition}
% ...
Definition~\ref{def:MPCC-MFCQ-T} has been widely used in the literature (see, e.g., \cite[Definition 3.1]{Flegel2005}, \cite[Definition 4]{Kim2020}, \cite[p. 4]{Scheel2000} and \cite[Definition 5.4]{schwartz2011mathematical}) and is usually referred to as the MPCC-MFCQ. However, we denote this as \MPCCMFCQT{} to emphasise the fact that it is based on the tightened program. In addition, this avoids conflict with the alternative definition that we will present next.

An alternative MPCC-tailored MFCQ relies on the Relaxed Nonlinear Program~\eqref{eq:RNLP} defined, for a given feasible point $\bar{\zzz}$ of problem \MPCC{}
and a subset $J\subseteq I^{GH}(\bar{\zzz})$ of the biactive indices, as
\begin{equation}
\tag{RNLP\mbox{$_J(\bar{\zzz})$}}
\label{eq:RNLP}
\begin{aligned}
\minimise_{\zzz\in\R^n} \quad
&f(\zzz)&
\\
\subjectto \quad
& g_i(\zzz) \geq 0&
&\text{ for } i=1,\dots,p,&
\\
&h_i(\zzz) = 0&
&\text{ for } i=1,\dots,q,&
\\
&G_{i}(\zzz)=0,    \quad H_{i}(\zzz)\geq0 \quad &
& \text{ for } i\in I^{G}(\bar{\zzz}),\\
&G_{i}(\zzz)\geq0, \quad H_{i}(\zzz)=0    \quad &
& \text{ for } i\in I^{H}(\bar{\zzz}),\\
&G_{i}(\zzz)\geq0,    \quad H_{i}(\zzz)\geq0    \quad &
& \text{ for } i\in I^{GH}(\bar{\zzz})\setminus{}J,\\
&%
G_{i}(\zzz)\leq0,    \quad H_{i}(\zzz)\leq0    \quad&
&
\text{ for } i\in    J.
\end{aligned}
\end{equation}
This problem is similar to the tightened program but for the biactive indices $i\in I^{GH}(\bar{\zzz})$, where inequality constraints rather than equality constraints are added for $G_i$ and $H_i$.
For $J=\emptyset$, as in the following Definition~\ref{def:MPCC-MFCQ-R}, these are the constraints $G_i(\zzz)\geq0$, $H_i(\zzz)\geq0$ for all $i\in I^{GH}(\bar{\zzz})$.
On the other hand, for a nonempty $J$, as in the following Definition~\ref{def:MPCC-MFCQ-C}, the signs of the constraints for $i\in J$ are reversed.
In either case, $\bar{\zzz}$ is feasible for~\eqref{eq:RNLP}, since $G_i(\bar{\zzz})=H_i(\bar{\zzz})=0$ for all $i\in I^{GH}(\bar{\zzz})$. Analogously, we use the notation~\eqref{eq:RNLP} to make clear that this program is defined at a feasible point $\bar{\zzz}$ and for a subset $J$. This leads to two interesting variants of the MPCC-tailored constraint qualification.

\begin{definition}The \MPCCMFCQR{} holds at a feasible point $\bar{\zzz}$ of the \MPCC{} iff the corresponding relaxed program~\RNLP{\emptyset} satisfies the \MFCQ{} at the same point $\bar{\zzz}$. Equivalently, the \MPCCMFCQR{} holds at a feasible point $\bar{\zzz}$ of the \MPCC{} iff
\label{def:MPCC-MFCQ-R}
\begin{align}
    &\Big\{\nabla h_i(\bar{\zzz}):i\in I^h\Big\} \cup
    \Big\{\nabla  G_i(\bar{\zzz}):i\in{I^{G}(\bar{\zzz})}\Big\} \cup
    \Big\{\nabla  H_i(\bar{\zzz}):i\in{I^{H}(\bar{\zzz})}\Big\}\notag\\
    &\text{is linearly independent and}
    \label{eq:MPCCMFCQR_linearly_independent_vectors}\\
    &\exists \ \ddd\in \R^n: \quad
    \begin{cases}
        \nabla h_i(\bar{\zzz}) ^{\top} \ddd = 0 &\forall i\in I^h,\\
        \nabla G_i(\bar{\zzz}) ^{\top} \ddd = 0 &\forall i\in I^{G}(\bar{\zzz}),\\
        \nabla H_i(\bar{\zzz}) ^{\top} \ddd = 0 &\forall i\in I^{H}(\bar{\zzz}),\\
        \nabla g_i(\bar{\zzz}) ^{\top} \ddd > 0 &\forall i\in I^{g}(\bar{\zzz}),\\
        \nabla G_i(\bar{\zzz}) ^{\top} \ddd > 0,\ \nabla H_i(\bar{\zzz}) ^{\top} \ddd > 0 &\forall i\in I^{GH}(\bar{\zzz}).\\
    \end{cases}
    \label{eq:MPCCMFCQR_direction_d}
\end{align}
\end{definition}

Definition~\ref{def:MPCC-MFCQ-R} is less common. It can be found in \cite[Definition 2.5]{Qiu2015} and \cite[Definition 2.5]{Ralph2004}. These papers call it \emph{MPCC-MFCQ} as well, but to avoid confusion, we will call it \MPCCMFCQR{}, denoting that it is based on the relaxed program.
We now introduce a third MPCC-tailored MFCQ, which lies strictly between the \MPCCMFCQT{} and the \MPCCMFCQR{} (see Theorem~\ref{theorem:cq_implications}) and is the one needed for the relaxation method of Section~\ref{section:Relaxation}.

\begin{definition}
\leavevmode\label{def:MPCC-MFCQ-C}%
The \MPCCMFCQC{} holds at a feasible point $\bar{\zzz}$ of the \MPCC{} iff the relaxed program~\eqref{eq:RNLP} satisfies the \MFCQ{} at $\bar{\zzz}$ for every subset $J\subseteq I^{GH}(\bar{\zzz})$. Equivalently, the \MPCCMFCQC{} holds at $\bar{\zzz}$ iff the family in~\eqref{eq:MPCCMFCQR_linearly_independent_vectors} is linearly independent and, for every $J\subseteq I^{GH}(\bar{\zzz})$,
\begin{equation}
\label{eq:MPCCMFCQC_direction_d}
\exists \ \ddd\in \R^n: \quad
\begin{cases}
    \nabla h_i(\bar{\zzz}) ^{\top} \ddd = 0 &\forall i\in I^h,\\
    \nabla G_i(\bar{\zzz}) ^{\top} \ddd = 0 &\forall i\in I^{G}(\bar{\zzz}),\\
    \nabla H_i(\bar{\zzz}) ^{\top} \ddd = 0 &\forall i\in I^{H}(\bar{\zzz}),\\
    \nabla g_i(\bar{\zzz}) ^{\top} \ddd > 0 &\forall i\in I^{g}(\bar{\zzz}),\\
    \nabla G_i(\bar{\zzz}) ^{\top} \ddd > 0,\ \nabla H_i(\bar{\zzz}) ^{\top} \ddd > 0 &\forall i\in I^{GH}(\bar{\zzz})\setminus J,\\
    \nabla G_i(\bar{\zzz}) ^{\top} \ddd < 0,\ \nabla H_i(\bar{\zzz}) ^{\top} \ddd < 0 &\forall i\in J.
\end{cases}
\end{equation}
\end{definition}

Why is such a definition worth introducing?
% ...
Firstly, it is precisely the CQ needed in our proof that the sequential Scholtes relaxation of Section~\ref{section:Relaxation} converges.
At a KKT point of the relaxed problem, the multipliers $\multG{_i}$ and $\multH{_i}$ of $\nabla G_i$ and $\nabla H_i$ for an index $i$ that is biactive in the limit come from one of two places:
either from the constraints $G_i(\zzz)\geq0$, $H_i(\zzz)\geq0$, in which case $\multG{_i}, \multH{_i}$ are both non-negative,
or from the constraint $G_i(\zzz)H_i(\zzz)\leq\relaxparam{}$, in which case $\multG{_i}, \multH{_i}$ are both non-positive.
These two cases are mutually exclusive, since $G_i(\zzz)=0$ or $H_i(\zzz)=0$ gives $G_i(\zzz)H_i(\zzz)=0<\relaxparam{}$. 
In the limit, the two multipliers therefore share a sign and $\multG{_i}\multH{_i}\geq0$, which is exactly the sign condition~\eqref{eq:C-stationary} of C-stationarity.
The \MPCCMFCQR{} only excludes non-negative multipliers, which is why it does not suffice.
The \MPCCMFCQT{} excludes all multipliers, which is overly strong.
The \MPCCMFCQC{}, in contrast, excludes exactly the sign pattern~$\multG{_i}\multH{_i}\geq0$, which is precisely what this argument requires.

Secondly, the \MPCCMFCQC{} expresses the same structure for other reformulations.
Applying Clarke's multiplier rule to the reformulation $\min(G_i(\zzz),H_i(\zzz))=0$ of the complementarity constraints gives Fritz John-type conditions whose multipliers have exactly the sign pattern~$\multG{_i}\multH{_i}\geq0$ excluded by the \MPCCMFCQC{}; the \MPCCMFCQC{} thus relates closely to the work of Scheel and Scholtes~\cite[Lemma~1]{Scheel2000}.
We note that a similar constraint qualification, called MPEC-NNAMCQ, was introduced by Ye~\cite[Definition~2.10]{Ye2005}, which excludes the multipliers with the sign pattern of M-stationarity.
Since the sign pattern of M-stationarity is a special case of that of C-stationarity, the \MPCCMFCQC{} implies the MPEC-NNAMCQ (Theorem~\ref{theorem:cq_implications}), so that any local minimiser of \MPCC{} at which the \MPCCMFCQC{} holds is M-stationary~\cite[Theorem~2.2]{Ye2005}.
To the best of our knowledge, the condition associated with the sign pattern of C-stationarity has not been stated explicitly before.

\begin{remark}
\label{def:MPCC-LICQ}
The \MPCCLICQ{} holds at point $\bar{\zzz}$ of problem \MPCC{} iff any of the following three equivalent statements holds:
\begin{enumerate}[leftmargin=2em,label=(\roman*)]
    \item The corresponding tightened program~\eqref{eq:TNLP} satisfies the \LICQ{} at the same point $\bar{\zzz}$.
    \item The corresponding relaxed program~\RNLP{\emptyset} satisfies the \LICQ{} at the same point $\bar{\zzz}$.
    \item The following family of vectors is linearly independent:
\begin{equation}
\label{eq:MPCCLICQ_linearly_independent_vectors}
\begin{split}
    \Big\{\nabla     h_i(\bar{\zzz}):   i\in I^h\Big\}
    \cup\Big\{\nabla g_i(\bar{\zzz}):   i\in I^{g}(\bar{\zzz})\Big\}\cup
    \Big\{\nabla     G_i(\bar{\zzz}):   i\in{I^{G}(\bar{\zzz})\cup I^{GH}(\bar{\zzz})}\Big\}\ \\
    \cup\Big\{\nabla H_i(\bar{\zzz}):   i\in{I^{H}(\bar{\zzz})\cup I^{GH}(\bar{\zzz})}\Big\}.
\end{split}
\end{equation}
\end{enumerate}
This means that the concept of the \MPCCLICQ{} is the same regardless of whether it is defined from the tightened or relaxed program.
\end{remark}

Observe that the \MPCCMFCQT{}%
, \MPCCMFCQR{} and \MPCCMFCQC{} are defined above in primal space due to the involvement of feasible directions. Next, we introduce the dual forms of these CQs, which are usually more handy in practice.
To proceed, consider two sets of vectors $A \coloneqq \{a_i: i\in I_1\}$ and $B\coloneqq \{b_i: i\in I_2\}$. We say that $(A\textbf{, }B)$ is \emph{positive-linearly dependent} iff there exist scalars $\{\multg{_i}\}_{i\in I_1}$ and $\{\multh{_i}\}_{i\in I_2}$, not all of them zero, such that
\[
\sum_{i\in I_1}\multg{_i} a_i + \sum_{i\in I_2}\multh{_i} b_i=0, \text{ where }
\multg{_i} \geq 0 \,\mbox{ for }\, i\in I_1 \,\mbox{ and } \,
\multh{_i} \text{ free} \, \mbox{ for } i \in I_2.
\]If no such scalars exist, we say $(A,B)$ is \emph{positive-linearly independent}.
% ...
The concept of positive-linear independence, which appears to have been introduced in \cite[p. 965, Definition 2.1]{Qi2000}, has since been widely used, in particular in the MPCC literature; see, e.g., \cite[p. 260, Definition 2.1]{Hoheisel2013} and \cite[p. 63, Definition 4.5]{schwartz2011mathematical}.
We now state a lemma that will be useful for many of the proofs below.
\begin{lemma}The following statements hold true: 
\label{lemma:positive-linearly-independent}
\begin{enumerate}[leftmargin=*,label=(\roman*)]
    \item The \MFCQ{} holds at a feasible point $\bar{\zzz}$ of problem~(\NLP{}) if and only if $(A,B)$ is positive-linearly independent where: \vspace{-.5em}
    \begin{align*}
    \label{eq:MFCQ_dual_form_part_1}
        &A=\{\nabla g_i(\bar{\zzz}):\, i\in I^g(\bar{\zzz})\},&
        &B=\{\nabla h_i(\bar{\zzz}):\, i\in I^h\}_{}.&
    \end{align*}
    \item The \MPCCMFCQT{} is satisfied at a feasible point $\bar{\zzz}$ of problem \MPCC{} if and only if $(A,B)$ is positive-linearly independent where:  \vspace{-.5em}
    \begin{equation}
    \label{eq:MFCQ_dual_form_part_2}
    \begin{aligned}
        &A=\{\nabla g_i(\bar{\zzz}) : i\in I^g(\bar{\zzz}) \}\textbf{, }\\
        & B=\left\{\nabla h_i(\bar{\zzz}):i\in{I^h}\right\} \cup 
        \left\{\nabla G_i(\bar{\zzz}):i\in{I^{G}(\bar{\zzz}) \cup I^{GH}(\bar{\zzz})}\right\}
        \cup \left\{\nabla H_i(\bar{\zzz}):i\in{I^{H}(\bar{\zzz})\cup I^{GH}(\bar{\zzz})}\right\}.
    \end{aligned}
    \end{equation}
    \item The \MPCCMFCQR{} is satisfied at a feasible point $\bar{\zzz}$ of problem \MPCC{} if and only if $(A,B)$ is positive-linearly independent where:  \vspace{-.5em}
    \begin{equation}
    \label{eq:MFCQ_dual_form_part_3}
    \begin{aligned}
        &A=\{\nabla g_i(\bar{\zzz}) : i\in I^g(\bar{\zzz}) \}\cup \left\{\nabla G_i(\bar{\zzz}), \nabla H_i(\bar{\zzz}):i\in{I^{GH}(\bar{\zzz})}\right\}\textbf{, }\\
        &B=\left\{\nabla h_i(\bar{\zzz}):i\in{I^h}\right\} \cup \left\{\nabla G_i(\bar{\zzz}):i\in{I^{G}(\bar{\zzz})}\right\} \cup \left\{\nabla H_i(\bar{\zzz}):i\in{I^{H}(\bar{\zzz})}\right\}.
    \end{aligned}
    \end{equation}
\end{enumerate}
\begin{enumerate}[leftmargin=*,label=(\roman*),resume]
    \item The \MPCCMFCQC{} is satisfied at a feasible point $\bar{\zzz}$ of problem \MPCC{} if and only if, for every subset $J\subseteq I^{GH}(\bar{\zzz})$, $(A_J,B)$ is positive-linearly independent where $B$ is given in~\eqref{eq:MFCQ_dual_form_part_3} and \vspace{-.5em}
    \begin{equation}
    \label{eq:MFCQ_dual_form_part_4}
    \begin{aligned}
        A_J=\{\nabla g_i(\bar{\zzz}) : i\in I^g(\bar{\zzz}) \}
        &\cup \left\{\nabla G_i(\bar{\zzz}), \nabla H_i(\bar{\zzz}):i\in{I^{GH}(\bar{\zzz})\setminus J}\right\}\\
        &\cup \left\{-\nabla G_i(\bar{\zzz}), -\nabla H_i(\bar{\zzz}):i\in J\right\}.
    \end{aligned}
    \end{equation}
    Equivalently, there exist no scalars $\{\multg{_i}\}_{i\in I^g(\bar{\zzz})}$, $\{\multh{_i}\}_{i\in I^h}$, $\{\multG{_i}\}_{i\in I^G(\bar{\zzz})\cup I^{GH}(\bar{\zzz})}$ and $\{\multH{_i}\}_{i\in I^H(\bar{\zzz})\cup I^{GH}(\bar{\zzz})}$, not all of them zero, such that
    \begin{equation}
    \label{eq:MPCCMFCQC_multipliers}
    \begin{gathered}
    \sum_{i\in I^g(\bar{\zzz})} \multg{_i}\nabla g_i(\bar{\zzz})
    +\sum_{i\in I^h} \multh{_i}\nabla h_i(\bar{\zzz})
    +\sum_{\mathclap{i\in I^G(\bar{\zzz})\cup I^{GH}(\bar{\zzz})}} \multG{_i}\nabla G_i(\bar{\zzz})
    +\sum_{\mathclap{i\in I^H(\bar{\zzz})\cup I^{GH}(\bar{\zzz})}} \multH{_i}\nabla H_i(\bar{\zzz})=0,\\
    \multg{_i}\geq 0 \ \text{ for } i\in I^g(\bar{\zzz}),\qquad
    \multG{_i}\multH{_i}\geq 0 \ \text{ for } i\in I^{GH}(\bar{\zzz}).
    \end{gathered}
    \end{equation}
    \item The MPEC-NNAMCQ~\cite[Definition~2.10]{Ye2005} holds at a feasible point $\bar{\zzz}$ of problem \MPCC{} iff there exist no scalars, not all of them zero, that satisfy the equality in~\eqref{eq:MPCCMFCQC_multipliers} with $\multg{_i}\geq0$ for $i\in I^g(\bar{\zzz})$ and, for every $i\in I^{GH}(\bar{\zzz})$, either $\multG{_i}>0$ and $\multH{_i}>0$, or $\multG{_i}\multH{_i}=0$.
\end{enumerate}
\end{lemma}

\begin{proof}
This follows directly from applying Motzkin's theorem of the alternative \cite[p. 27]{Mangasarian1969} to the definitions.
For the second statement in~(iv), observe that $\multG{_i}\multH{_i}\geq0$ holds iff $\multG{_i}$ and $\multH{_i}$ are both non-negative or both non-positive. Given scalars, not all zero, satisfying~\eqref{eq:MPCCMFCQC_multipliers}, let $J$ be the set of indices $i\in I^{GH}(\bar{\zzz})$ with $\multG{_i}<0$ or $\multH{_i}<0$. Then the coefficients $-\multG{_i}$, $-\multH{_i}$ of $-\nabla G_i(\bar{\zzz})$, $-\nabla H_i(\bar{\zzz})$ for $i\in J$ and $\multG{_i}$, $\multH{_i}$ of $\nabla G_i(\bar{\zzz})$, $\nabla H_i(\bar{\zzz})$ for $i\in I^{GH}(\bar{\zzz})\setminus J$ are all non-negative, so $(A_J,B)$ is positive-linearly dependent. Conversely, a positive-linear dependence of some $(A_J,B)$ yields scalars satisfying~\eqref{eq:MPCCMFCQC_multipliers} after changing the signs of the coefficients of $-\nabla G_i(\bar{\zzz})$ and $-\nabla H_i(\bar{\zzz})$ for $i\in J$.
\end{proof}

We have the following sequence of implications.
The first is well known, see for example~\mbox{\cite[Corollary~3.2]{Flegel2005}}, and the implication from \MPCCMFCQT{} to MPEC-NNAMCQ was given in~\mbox{\cite[Theorem~3.2]{Ye2005}}.
\begin{theorem}
\leavevmode\label{theorem:cq_implications}%
For any feasible point $\bar{\zzz}$ of problem \MPCC{}, the following implications hold:
\begin{equation*}
        \text{\MPCCLICQ{}}
        \implies
        \text{\MPCCMFCQT{}}
        \implies
        \text{\MPCCMFCQC{}}%
        \implies
        \text{MPEC-NNAMCQ}
        \implies
        \text{\MPCCMFCQR{}}
\end{equation*}
\end{theorem}

\begin{proof}
For the second implication, assume that the \MPCCMFCQC{} does not hold at $\bar{\zzz}$.
By Lemma~\ref{lemma:positive-linearly-independent}(iv), there exist scalars, not all of them zero, satisfying~\eqref{eq:MPCCMFCQC_multipliers}.
In particular, the coefficients $\multg{_i}$ of $\nabla g_i(\bar{\zzz})$, $i\in I^g(\bar{\zzz})$, are non-negative, while every other vector in~\eqref{eq:MPCCMFCQC_multipliers} belongs to the set $B$ in~\eqref{eq:MFCQ_dual_form_part_2}, whose coefficients are unrestricted in sign.
This contradicts the positive-linear independence of~\eqref{eq:MFCQ_dual_form_part_2}. Therefore, \MPCCMFCQT{} also fails.
For the third implication, note that $\multG{_i}>0$ and $\multH{_i}>0$, or $\multG{_i}\multH{_i}=0$, implies $\multG{_i}\multH{_i}\geq0$. Hence, any scalars, not all of them zero, that violate the MPEC-NNAMCQ also satisfy~\eqref{eq:MPCCMFCQC_multipliers}, and the \MPCCMFCQC{} fails by Lemma~\ref{lemma:positive-linearly-independent}(iv).
For the fourth implication, assume that the \MPCCMFCQR{} does not hold at $\bar{\zzz}$. By Lemma~\ref{lemma:positive-linearly-independent}(iii), the pair $(A,B)$ in~\eqref{eq:MFCQ_dual_form_part_3} is positive-linearly dependent; that is, there exist scalars, not all of them zero, that satisfy the equation in~\eqref{eq:MPCCMFCQC_multipliers} with $\multg{_i}\geq0$ for $i\in I^g(\bar{\zzz})$ and $\multG{_i},\multH{_i}\geq0$ for $i\in I^{GH}(\bar{\zzz})$. Each such pair $(\multG{_i},\multH{_i})$ either has both components positive or satisfies $\multG{_i}\multH{_i}=0$, so the MPEC-NNAMCQ fails.
\end{proof}
% ...
The implications in Theorem~\ref{theorem:cq_implications} are strict; for the first one this is well known, and for the others it is shown in the following example.
\begin{example} 
\label{eq:Example01} Consider the following~(\MPCC{}): 
\begin{equation*}
\begin{aligned}
\minimise_{\zzz\in \R^2} \quad
&z_1 + z_2&
\\
\subjectto \quad
&z_1 + z_2 \geq 0,& \\
&0 \leq z_1 \bot z_2 \geq 0& 
\end{aligned}
\end{equation*}
The inequality constraint is $g(\zzz) \coloneqq z_1 + z_2$ and complementarity is $G(\zzz) \coloneqq z_1$ and $H(\zzz) \coloneqq z_2$.
At the point $\bar{\zzz}=[0,0]^\top$, we are in the biactive case where $G(\bar{\zzz})=H(\bar{\zzz})=0$. Consider the derivatives of the constraints: $\nabla g(\bar{\zzz})=[1,1]^\top$, $\nabla G(\bar{\zzz})=[1,0]^\top$ and $\nabla H(\bar{\zzz})=[0,1]^\top$. Choosing, for example, $\ddd=[3,4]^\top$, we get $\ddd^{\top}\nabla g(\bar{\zzz})=7>0$ and $\ddd^{\top}\nabla G(\bar{\zzz})=3>0$ and $ \ddd^{\top}\nabla H(\bar{\zzz})=4>0$. Thus, the  \MPCCMFCQR{} holds at $\bar{\zzz}$. However, there does not exist any direction $\ddd\in\R^2$ such that $\ddd^{\top}\nabla g(\bar{\zzz})>0$ and $\ddd^{\top}\nabla G(\bar{\zzz})=0$ and $\ddd^{\top}\nabla H(\bar{\zzz})=0$. This implies the failure of the  \MPCCMFCQT{} at $\bar{\zzz}$.

The same example shows that \MPCCMFCQR{} does not imply \MPCCMFCQC{}.
Consider $J=I^{GH}(\bar{\zzz})$, no direction $\ddd\in\R^2$ satisfies $\ddd^{\top}\nabla g(\bar{\zzz})>0$, $\ddd^{\top}\nabla G(\bar{\zzz})<0$ and $\ddd^{\top}\nabla H(\bar{\zzz})<0$.
Equivalently, $\multg{_{}}=1$ and $\multG{_{}}=\multH{_{}}=-1$ satisfy~\eqref{eq:MPCCMFCQC_multipliers}.
In fact, every solution of the equation in~\eqref{eq:MPCCMFCQC_multipliers} for this example has $\multG{_{}}=\multH{_{}}=-\multg{_{}}$ with $\multg{_{}}\geq0$, so no nonzero solution has the sign pattern of the MPEC-NNAMCQ. The example therefore satisfies the MPEC-NNAMCQ, which consequently does not imply the \MPCCMFCQC{}.
Adding a second copy of the constraint $0\leq z_1\bot z_2\geq 0$ to the example preserves the \MPCCMFCQR{}, since $\ddd=[3,4]^\top$ still satisfies~\eqref{eq:MPCCMFCQR_direction_d}, but violates the MPEC-NNAMCQ since the multipliers $\multG{_1}=1$ and $\multG{_2}=-1$ of the two copies of $\nabla G(\bar{\zzz})$, with all other multipliers zero, satisfy the equation in~\eqref{eq:MPCCMFCQC_multipliers} and have $\multG{_i}\multH{_i}=0$.
On the other hand, if the constraint $z_1+z_2\geq0$ is replaced by $z_1-z_2\geq 0$, then the \MPCCMFCQC{} holds at $\bar{\zzz}$, since $\ddd=[2,1]^\top$ (for $J=\emptyset$) and $\ddd=[-1,-2]^\top$ (for $J=I^{GH}(\bar{\zzz})$) satisfy~\eqref{eq:MPCCMFCQC_direction_d}, while the \MPCCMFCQT{} still fails for the same reason as above.
Hence, all implications in Theorem~\ref{theorem:cq_implications} that involve the \MPCCMFCQC{} or the MPEC-NNAMCQ are strict.
\end{example}

Another important example that satisfies
\MPCCMFCQC{} (and hence \MPCCMFCQR{}%
) but not \MPCCMFCQT{} is the main topic of Section~\ref{section:Application}. For further evidence of the very strong nature of the \MPCCMFCQT{}, we end on the following remark.

\begin{remark}
    \label{remark:MPCCMFCQT_boils_down_to_MPCCLICQ}
     For any problem of the form~(\MPCC{}) with $p=0$ (in other words with no inequality $g_i(\zzz)\geq0$ constraints but still possibly with equality $h_i(\zzz)=0$ and complementarity $0\leq G_i(\zzz)\bot H_i(\zzz)\geq 0$ constraints), the \MPCCMFCQT{} is equivalent to the \MPCCLICQ{}. 
\end{remark}

%% file: Article/article_22_stationarity.tex
\subsection{Stationarity}
\label{section:stationarity}
% ...
This subsection states some well-known stationarity concepts for problem \MPCC, which have been extensively studied in the literature; see, e.g.,~\cite{Scheel2000}, which seems to be one of the early papers on the subject. We begin by introducing a problem-specific Lagrangian function
\begin{equation}
\label{eq:MPCC_Lagragian}
    \mathcal{L}(\zzz,\multg{}, \multh{}, \multG{}, \multH{})  \coloneqq  
    f(\zzz)
    - \sum_{i=1}^{p} \multg{_i}g_i(\zzz)
    - \sum_{i=1}^{q} \multh{_i}h_i(\zzz)
    - \sum_{i=1}^{r} \multG{_i}G_i(\zzz)
    - \sum_{i=1}^{r} \multH{_i}H_i(\zzz).
\end{equation}
Observe that this function does not include the product terms $G_i(\zzz)H_i(\zzz)$ for $i=1, \ldots, r$. 

\begin{definition}
\label{definition:stationarity}
The feasible point $\bar{\zzz}$ of problem \MPCC{} is said to be Weakly (W)-stationary if there exist multipliers $(\multg{}{}, \multh{}{}, \multG{}{}, \multH{}{})$ such that the following conditions hold:
\begin{subequations}
\label{eq:stationarity}
\begin{align}
\label{eq:stationarity_a}
&\nabla_{\zzz}\mathcal{L}(\bar{\zzz},\multg{}, \multh{}, \multG{}, \multH{}) = 0
\\
\label{eq:stationarity_b}
&\multg{_i}\geq0
\text{ for } i\in I^g(\bar{\zzz}),\quad
\multg{_i} = 0
\text{ for } i\notin I^g(\bar{\zzz}),\quad
\multH{_i} = 0
\text{ for } i \in I^{G}(\bar{\zzz}),\quad
\multG{_i} = 0
\text{ for } i \in I^{H}(\bar{\zzz}).
\end{align}
\end{subequations}
The point $\bar{\zzz}$ will be said to be Alternative (A), Clarke (C), Mordukhovich (M) or Strongly (S)-stationary if additionally to \eqref{eq:stationarity_b}, we respectively have the conditions
\begin{align}
\label{eq:A-stationary}
\tag{\ref*{eq:stationarity}A}
&\text{A-stationary:} &
&\multG{_i} \geq 0 \text{  or  }   \multH{_i}{}\geq 0&
&\text{for } i \in I^{GH}(\bar{\zzz});&
\\
\label{eq:C-stationary}
\tag{\ref*{eq:stationarity}C}
&\text{C-stationary:}&
&\multG{_i} \multH{_i} \geq 0&
&\text{for } i \in I^{GH}(\bar{\zzz});&
\\
\label{eq:M-stationary}
\tag{\ref*{eq:stationarity}M}
&\text{M-stationary:}&
&\multG{_i}, \multH{_i} \geq 0 \text{  or  }   \multG{_i} \multH{_i} =0&
&\text{for } i \in I^{GH}(\bar{\zzz});&
\\
\label{eq:S-stationary}
\tag{\ref*{eq:stationarity}S}
&\text{S-stationary:}&
&\multG{_i}, \multH{_i} \geq 0&
&\text{for } i \in I^{GH}(\bar{\zzz}).&
\end{align}
\end{definition}

It is important to note a few properties of these stationarity conditions. They only differ by the nature of multipliers on the biactive index set $I^{GH}$. 
\ifsuplementaryincluded
Figure~\ref{fig:Stationary} highlights these differences.
\fi
So, if the problem is strictly complementary at $\bar{\zzz}$, that is
$G_i(\bar{\zzz})>0,H_i(\bar{\zzz})=0$ or $G_i(\bar{\zzz})=0,H_i(\bar{\zzz})>0$
for all $i=1,\dots,r$, then all five concepts are equivalent.  
S-stationary is the strongest and is equivalent to the KKT conditions of the \MPCC{} if the complementarity  constraints are viewed as ordinary inequality and equality constraints. The following relationships hold \cite[p. 60]{Flegel2005}: 
% ...
\newcommand*{\NEArrow}{\rotatebox[origin=c]{15}{\(\implies\)}}
\newcommand*{\SEArrow}{\rotatebox[origin=c]{-15}{\(\implies\)}}
\begin{equation*}
    \text{S-stationary }
    \implies
    \text{M-stationary}
    \quad
    \begin{matrix}
    \NEArrow & \text{A-stationary } & \SEArrow\\
    \SEArrow & \text{C-stationary } & \NEArrow
    \end{matrix}
    \quad
    \text{W-stationary}
\end{equation*}

%% file: Article/article_30_solution_methods.tex
\section{Solution methods}
\label{section:solution_methods}
% ...
\noindent In this section, we study the partial penalisation and Scholtes relaxation methods. 
There are at least two motivations to study these algorithms. 
First, they are among the most commonly used in the literature on MPCCs (see, e.g., \cite{Kim2020}).
Secondly, there is a nice symmetry between the two approaches. In the penalisation method, the constraint $G(\zzz)^\top H(\zzz) =0$ is moved from the feasible set to the objective function, while this same constraint is instead relaxed to form a larger feasible set in the context of the relaxation approach.

As will become clear in the next two subsections, a common point between the two methods, which does not seem to have been discovered before, is that they converge under the same CQ; i.e., the
\MPCCMFCQC{}, which is strictly weaker than
the \MPCCLICQ{} and \MPCCMFCQT{} used in \cite{Leyffer2006interior,schwartz2011mathematical} to prove the convergence of these algorithms.
In both cases, the limit point is C-stationary, and it is S-stationary under the weaker \MPCCMFCQR{} together with an additional assumption at the biactive indices.


%% file: Article/article_31_penalisation.tex

\subsection{Penalisation}
\label{section:Penalisation}
% ...
We start by recalling that the penalisation method works by moving the problematic product term $G(\zzz)^\top H(\zzz)$ from the constraints of problem (\MPCC{}) to its objective function. This leads to the penalised problem with $\penaltyparam{}\geq0$ being the penalty parameter:
\begin{equation}
\label{eq:penalisation}
\tag{\mbox{$\mathcal{P}_\penaltyparam{}$}} 
\begin{aligned}
\minimise_{\zzz \in \R^n} \quad
&f(\zzz)\ +\penaltyparam{} \ G(\zzz)^{\top} H(\zzz) 
\\
\subjectto \quad
& g_i(\zzz)\geq 0&  \text{for } i=1,...,p,\\
& h_i(\zzz) =   0&  \text{for } i=1,...,q,\\
& G_i(\zzz)\geq 0,\quad
  H_i(\zzz)\geq 0& \text{for } i=1,...,r.
\end{aligned}
\end{equation}
Interestingly, \penalisation{\penaltyparam{}} is more likely to satisfy standard constraint qualifications.  In fact, if \MPCC{} satisfies \MPCCMFCQR{} at a point $\zzz^*$, then \penalisation{\penaltyparam{}} will satisfy the \MFCQ{} at $\zzz^*$ \cite[Theorem 5.1]{Ralph2004}. This means that we can employ standard NLP solvers such as the interior-point method~(IPM) or sequential quadratic programming~(SQP) for \penalisation{\penaltyparam{}}. However, when employing such solvers, we should not expect them to find a global minimum but instead terminate within a certain tolerance $\epsilon\geq0$ of a KKT point. With this in mind, we define a point $\zzz\in\R^n$ to be an $\epsilon$-KKT point for \penalisation{\penaltyparam{}} if there exists a Lagrange multiplier vector $\left(\multg{}, \multh{}, \multG{}, \multH{}\right)\in\R^{(p+q+2r)}$ satisfying the following system:
\begin{subequations}
\label{eq:KKT_of_penalised}
\begin{align}
&\begin{aligned}
\Big\|
\nabla{}f(\zzz) 
 + \penaltyparam{} H(\zzz)\left[\nabla G(\zzz)\right] 
 + \penaltyparam{} G(\zzz)\left[\nabla H(\zzz)\right]
 - \sum_{i=1}^{p} \multg{_i} \nabla  g_i(\zzz)
 - \sum_{i=1}^{q} \multh{_i} \nabla  h_i(\zzz)&\\
 - \sum_{i=1}^{r} \multG{_i} \nabla  G_i(\zzz)
 - \sum_{i=1}^{r} \multH{_i} \nabla  H_i(\zzz)
&\Big\|\leq \epsilon,\label{eq:KKT_of_penalised_a}
\end{aligned}\\
&\begin{array}{r@{\quad} r@{\quad} r@{\quad} r}
g_i(\zzz) \geq 0, & \multg{_i} \geq 0,       & g_i(\zzz)\multg{_i}=0,  & \text{ for }i=1,\dots,p,\\
h_i(\zzz) =    0, & \multh{_i} \text{ free}, &                         & \text{ for }i=1,\dots,q,\\
G_i(\zzz) \geq 0, & \multG{_i} \geq 0,       & G_i(\zzz)\multG{_i}=0,  & \text{ for }i=1,\dots,r,\\
H_i(\zzz) \geq 0, & \multH{_i} \geq 0,       & H_i(\zzz)\multH{_i}=0,  & \text{ for }i=1,\dots,r.\label{eq:KKT_of_penalised_b}
\end{array}
\end{align}
\end{subequations}

\begin{figure}
    \centering
    \includegraphics[width=\linewidth]{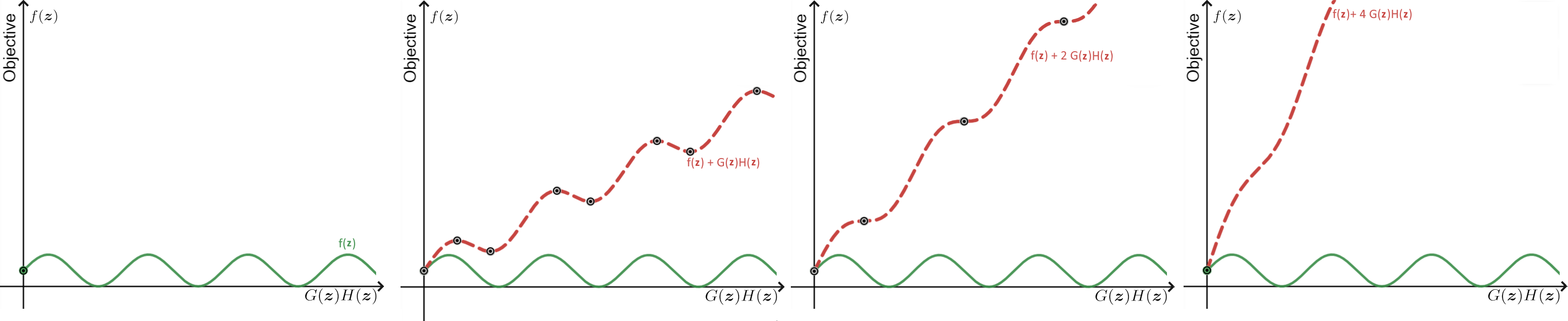}
    \caption{This figure visualises applying the penalisation sequence to a hypothetical MPCC with a single complementarity constraint. The vertical axis is the objective value, while the horizontal axis is the product of the two complementarity constraint functions $G(\zzz) H(\zzz)$. The green line shows the objective value $f(\zzz)$, and the red curve shows the penalised objective $f(\zzz)+\penaltyparam{} G(\zzz)H(\zzz)$ for increasing penalty parameters $\penaltyparam{}$. The grey dots show infeasible stationary points that are slowly ironed out. 
    }
    \label{fig:penalisation_sequence}
\end{figure}

This concept is widely used in the literature; see, e.g.,~\cite[Definition 3.1]{Dutta2013}.
For $\epsilon=0$ the system (\link{eq:KKT_of_penalised}{\ref*{eq:KKT_of_penalised}a-b}) is the classical KKT conditions for~\penalisation{\penaltyparam{}}.
The idea is then for a large enough penalty parameter~$\penaltyparam{}$ and small enough tolerance~$\epsilon$, finding a solution to (\link{eq:KKT_of_penalised}{\ref*{eq:KKT_of_penalised}a-b}) that will correspond to a stationary point of problem \MPCC{}.
If we know such values for $\penaltyparam{}$ and $\epsilon$, we only need to solve \penalisation{\penaltyparam{}} once.  
This is called \textit{exact penalisation}. 
However, in practice, we may not have this information, or in the case where $\penaltyparam{}$ is too large, \penalisation{\penaltyparam{}} may become ill-conditioned \cite[Example 3]{Leyffer2006interior}. 
Therefore, a more practical approach is to introduce an increasing sequence of non-negative penalty parameters $(\penaltyparam{}^{t})_{t\in\N}$ and a decreasing sequence of tolerances $(\epsilon^{t})_{t\in\N}$ and solve problem \penalisation{\penaltyparam{}^t} repeatedly using the previous solution as the starting point for the next iteration until our stopping criteria are met. We call this \textit{sequential penalisation}. 
A precise description of the method is provided in Algorithm~\ref{alg:penalisation} and its theoretical convergence is established in the result below. A graphical illustration, for increasing penalty parameter values, is given in Figure~\ref{fig:penalisation_sequence}.

\renewcommand{\algorithmicrequire}{\textbf{Input:}}
\renewcommand{\algorithmicensure}{\textbf{Procedure:}}
\begin{algorithm}[H]
\caption{Inexact sequential partial penalisation}
\label{alg:penalisation}
\begin{algorithmic}[1]

\Require{} Access to the \MPCC{}'s functions $f,g,h,G,H$ and their first order derivatives;
Starting point $\zzz^{0}\in \R^n$;
Initial penalty parameter $\penaltyparam{}^{1}>0$;
Initial KKT tolerance $\epsilon^1 > 0$;
Exit tolerance $\epsilon_{\text{exit}} > 0$;
\vspace{.5em}

\Ensure
\Repeat{ for $t=1,2,3,\dots$:}
\label{line:find_minimiser}
    % ...
    \State%
    \label{line:find_pen_kkt}%
    Find an $\epsilon^t$-KKT point $\zzz^{t}$ of \penalisation{\penaltyparam{}^{t}} with starting point $\zzz^{t-1}$;
    
    % ...
    \State%
    \label{line:increase_pen_param}%
    Increase the penalty parameter $\penaltyparam{}^{t+1} > \penaltyparam{}^t$%
    , ensuring $\penaltyparam{}^{t}\rightarrow\infty$;
    % ...
    \State%
    \label{line:decrease_pen_tol}%
    Decrease the KKT tolerance $\epsilon^{t+1} < \epsilon^t$, ensuring $\epsilon^{t}\rightarrow0$;
\Until{Stopping condition: $\max_{i}\left( \min(G_i(\zzz^t), \;H_i(\zzz^t)) \right)\leq \epsilon_{\text{exit}}$.}
\State Output $\zzz^t$.
\end{algorithmic}
\end{algorithm}

\begin{theorem}
\leavevmode\label{theorem:penalisation_convergence}%
Consider problem \MPCC{} and a sequence of iterates $\left(\zzz^t\right)_{t\in\N}$, generated by Algorithm~\ref{alg:penalisation}, which is assumed to converge to the point $\zzz^*$ where $G(\zzz^*)^\top H(\zzz^*)=0$, then
\begin{enumerate}[label=(\roman*),leftmargin=2em]
    \item If
    \MPCCMFCQC{} holds at $\zzz^*$, then $\zzz^*$ is a C-stationary point.
    \item If \MPCCMFCQR{} holds at $\zzz^*$ but additionally $\penaltyparam{}^tG_i(\zzz^t)\rightarrow0$, $\penaltyparam{}^tH_i(\zzz^t)\rightarrow0$, 
    for the biactive indices $i\in I^{GH}(\zzz^*)$
    then $\zzz^*$ is an S-stationary point.
\end{enumerate}
\end{theorem}

\begin{proof}
For each iterate $\zzz^t$, since it is an $\epsilon^t$-KKT point of \penalisation{\penaltyparam{}^t}, let $\left(\multg{^t}, \multh{^t}, \multG{^t}, \multH{^t}\right)\in\R^{(p+q+2r)}$ be a Lagrange multiplier vector satisfying (\link{eq:KKT_of_penalised}{\ref*{eq:KKT_of_penalised}a-b}).
We now define augmented multipliers:
\begin{align}
\label{eq:mu_definition}
\multG[hat]{_i^t}  \coloneqq  \multG{^t_i} - \penaltyparam{}^t H_i(\zzz^t)  \; \mbox{ and } \;
\multH[hat]{_i^t}  \coloneqq  \multH{^t_i} - \penaltyparam{}^t G_i(\zzz^t).
\end{align}
In this way we can rewrite (\ref{eq:KKT_of_penalised_a}) as
\begin{equation}
\label{eq:KKT_of_penalised_a_rearranged}
\left\|
\nabla f(\zzz^t)
 - \sum_{i=1}^{p} \multg{_i^t}      \nabla  g_i(\zzz^t)
 - \sum_{i=1}^{q} \multh{_i^t}      \nabla  h_i(\zzz^t)
 - \sum_{i=1}^{r} \multG[hat]{_i^t} \nabla  G_i(\zzz^t)
 - \sum_{i=1}^{r} \multH[hat]{_i^t} \nabla  H_i(\zzz^t)
 \right\|
 \leq \epsilon^t.
\end{equation}
We wish to prove that the sequence $\left(\aaa_t\right)_{t\in\N}$, where $\aaa_t \coloneqq \left(\multg{^t}, \multh{^t}, \multG[hat]{^t}, \multH[hat]{^t}\right)$, has a convergent subsequence. Suppose towards a contradiction that it does not. Then $\|\aaa_t\|\rightarrow\infty$; for otherwise some subsequence of $\left(\aaa_t\right)_{t\in\N}$ would be bounded and would, by the Bolzano--Weierstrass theorem, admit a convergent sub-subsequence.
Now, dividing equation (\ref{eq:KKT_of_penalised_a_rearranged}) through by $\left\|\aaa_t\right\|$ leads to
\begin{multline}
\label{eq:KKT_of_penalised_a_normed}
\left\|
\frac{\nabla f(\zzz^t)}{\|\aaa_t\|}
 - \sum_{i=1}^{p} \frac{\multg{_i^t}}{{\|\aaa_t\|}} \nabla  g_i(\zzz^t)
 - \sum_{i=1}^{q}  \frac{\multh{_i^t}}{{\|\aaa_t\|}} \nabla  h_i(\zzz^t)\right.
\left. - \sum_{i=1}^{r} \frac{\multG[hat]{_i^t}}{{\|\aaa_t\|}} \nabla  G_i(\zzz^t)
 - \sum_{i=1}^{r}  \frac{\multH[hat]{_i^t}}{{\|\aaa_t\|}} \nabla  H_i(\zzz^t)
\right\|
 \leq \frac{\epsilon^t}{\|\aaa_t\|}.
\end{multline}
We consider the normalised sequence of Lagrange multipliers $\aaa_t/\left\|\aaa_t\right\|$.
As it is bounded, by the Bolzano--Weierstrass theorem, there exists some infinite subsequence that converges.
Let $(\multg{'}, \multh{'}, \multG[hat]{'}, \multH[hat]{'})$ be the limit of that subsequence.
Since every term of the normalised sequence has unit norm, so does its limit, by continuity of the norm; in particular $(\multg{'}, \multh{'}, \multG[hat]{'}, \multH[hat]{'})\neq0$.
We next identify the components that vanish in the limit.
For $i\in \{1,\dots,p\} \setminus I^g$ we have $g_i(\zzz^*)>0$, hence $g_i(\zzz^t)>0$ for all large $t$, and the complementarity slackness in~\eqref{eq:KKT_of_penalised_b} gives $\multg{^t_i}=0$; therefore $\multg{'_i}=0$.
Now consider $i\in I^H$, so that $G_i(\zzz^*)>0$ and $H_i(\zzz^*)=0$. Then $G_i(\zzz^t)>0$ for all large $t$, so complementarity slackness gives $\multG{^t_i}=0$ and hence $\multG[hat]{_i^t}=-\penaltyparam{}^t H_i(\zzz^t)$ by~\eqref{eq:mu_definition}. If $H_i(\zzz^t)=0$, then $\multG[hat]{_i^t}=0$. Otherwise, complementarity slackness gives $\multH{^t_i}=0$, so that $\multH[hat]{_i^t}=-\penaltyparam{}^t G_i(\zzz^t)\neq0$, and since $|\multH[hat]{_i^t}|\leq\|\aaa_t\|$,
\begin{equation*}
\frac{\bigl|\multG[hat]{_i^t}\bigr|}{\|\aaa_t\|}
\;\leq\; \frac{\bigl|\multG[hat]{_i^t}\bigr|}{\bigl|\multH[hat]{_i^t}\bigr|}
\;=\; \frac{\penaltyparam{}^t H_i(\zzz^t)}{\penaltyparam{}^t G_i(\zzz^t)}
\;=\; \frac{H_i(\zzz^t)}{G_i(\zzz^t)}
\;\longrightarrow\; \frac{0}{G_i(\zzz^*)}=0 .
\end{equation*}
In either case $\multG[hat]{'_i}=0$ for $i\in I^H$; exchanging the roles of $G$ and $H$ gives $\multH[hat]{'_i}=0$ for $i\in I^G$. Using this indexing, we can say that it must be the $\left(\multg{'_{I^g}}, \multh{'}, \multG[hat]{_{I^G\cup I^{GH}}'}, \multH[hat]{_{I^H\cup I^{GH}}'}\right)$ components of the sequence that are non-zero.
Now, returning to \eqref{eq:KKT_of_penalised_a_normed}, as $t\rightarrow\infty$, the terms $\frac{\nabla f(\zzz^t)}{\|\aaa_t\|}$ and $\frac{\epsilon^t}{\|\aaa_t\|}$ go to zero, leaving the following linear combination equal to zero
\begin{equation}
\label{eq:p_linearly_dependence}
    \sum_{i\in I^g} \multg{'_i} \nabla  g_i(\zzz^*) +
    \sum_{i\in I^h} \multh{'_i} \nabla  h_i(\zzz^*)+
    \sum_{\mathclap{i\in I^{G}\cup I^{GH}}} \multG[hat]{'_i} \nabla  G_i(\zzz^*) +
    \sum_{\mathclap{i\in I^{H}\cup I^{GH}}} \multH[hat]{'_i} \nabla  H_i(\zzz^*)=0.
\end{equation}
Under the assumptions of part~(i), since $\multg{_i'}\geq0$
and, for $i\in I^{GH}$, $\multG[hat]{'_i}\multH[hat]{'_i}\geq0$ as the limit of $\multG[hat]{_i^t}\multH[hat]{_i^t}/\|\aaa_t\|^2$, which is non-negative for every $t$ by the identity shown below in the proof of C-stationarity, the 
non-zero vector $(\multg{'}, \multh{'}, \multG[hat]{'}, \multH[hat]{'})$ satisfies~\eqref{eq:MPCCMFCQC_multipliers} at $\zzz^*$, contradicting
\MPCCMFCQC{} by Lemma~\ref{lemma:positive-linearly-independent}(iv).
Under the assumptions of part~(ii), i.e. $\penaltyparam{}^tG_i(\zzz^t)\rightarrow0$ and $\penaltyparam{}^tH_i(\zzz^t)\rightarrow0$ for $i\in I^{GH}$, we have that $\multG[hat]{'_i}\geq0$ and $\multH[hat]{'_i}\geq0$ for $i\in I^{GH}$, so that the sets of gradients given in~\eqref{eq:MFCQ_dual_form_part_3} are positive-linearly dependent, contradicting \MPCCMFCQR{}.
By these contradictions, we conclude in both cases (i) and (ii) that $\left(\aaa_t\right)_{t\in\N}$ has a convergent subsequence. 
Let $\left(\multg{^*}, \multh{^*}, \multG[hat]{^*}, \multH[hat]{^*}\right)$ be the limit of that convergent subsequence.
By the continuous differentiability of $f, g, h, G, H$, and since the sequence $\left(\zzz^t\right)_{t\in\N}$ converges to $\zzz^*$ in their domain, these functions and their derivatives must converge in range. Applying this to inequality~(\ref{eq:KKT_of_penalised_a_rearranged}), and recalling that $\epsilon^t\rightarrow0$, we get
\begin{equation}
\label{eq:KKT_stationarity_of_MPCC}
\nabla f(\zzz^*)
 - \sum_{i=1}^{p} \multg{_i^*}      \nabla  g_i(\zzz^*)
 - \sum_{i=1}^{q} \multh{_i^*}      \nabla  h_i(\zzz^*)
 - \sum_{i=1}^{r} \multG[hat]{_i^*} \nabla  G_i(\zzz^*)
 - \sum_{i=1}^{r} \multH[hat]{_i^*} \nabla  H_i(\zzz^*)=0,
\end{equation}
satisfying the stationarity condition~\eqref{eq:stationarity_a} of Definition~\ref{definition:stationarity}.
To establish W-stationarity, it remains to be shown that condition~\eqref{eq:stationarity_b} is satisfied, i.e., $\multG[hat]{_i^*}=0$ for each $i\in I^H$ and $\multH[hat]{_i^*}=0$ for each $i\in I^G$.
Considering an index $i\in I^G$, by definition, we have $G_i(\zzz^*)=0$ and $H_i(\zzz^*)>0$. Since $H_i(\zzz^t)\rightarrow H_i(\zzz^*)>0$, complementarity slackness gives $\multH{^t_i}=0$ for all large $t$. To deduce $\multH[hat]{_i^*}=0$ from~\eqref{eq:mu_definition}, the penalty term $\penaltyparam{}^t G_i(\zzz^t)$ must vanish in the limit. For every $t$ with $G_i(\zzz^t)>0$, complementarity slackness forces $\multG{^t_i}=0$, so that $\multG[hat]{_i^t}=-\penaltyparam{}^t H_i(\zzz^t)$; since $\multG[hat]{_i^t}$ converges to the finite limit $\multG[hat]{_i^*}$ while $H_i(\zzz^t)$ is bounded away from zero, the penalty parameter $\penaltyparam{}^t$ remains bounded for all such $t$. Consequently $\penaltyparam{}^t G_i(\zzz^t)\rightarrow0$, because for each $t$ either $G_i(\zzz^t)=0$ or $\penaltyparam{}^t$ is bounded while $G_i(\zzz^t)\rightarrow0$; therefore $\multH[hat]{_i^*}=0$ for all $i\in I^G$. Exchanging the roles of $G$ and $H$, the same argument gives $\multG[hat]{_i^*}=0$ for all $i\in I^H$.

To establish C-stationarity~\eqref{eq:C-stationary} under the assumptions of Theorem~\ref{theorem:penalisation_convergence}~(i), consider a biactive index $i\in I^{GH}$ and the product $\multG[hat]{_i^t}\,\multH[hat]{_i^t}$. Expanding it and cancelling the two cross terms by the complementarity slackness $\multG{^t_i}G_i(\zzz^t)=0$ and $\multH{^t_i}H_i(\zzz^t)=0$ of~\eqref{eq:KKT_of_penalised_b}, we obtain, for every $t\in\N$,
\begin{equation*}
\multG[hat]{_i^t}\,\multH[hat]{_i^t}
= \bigl(\multG{^t_i}-\penaltyparam{}^t H_i(\zzz^t)\bigr)\bigl(\multH{^t_i}-\penaltyparam{}^t G_i(\zzz^t)\bigr)
= \multG{^t_i}\,\multH{^t_i} + \bigl(\penaltyparam{}^t\bigr)^{2} G_i(\zzz^t)\,H_i(\zzz^t)
\;\geq\; 0.
\end{equation*}

To establish S-stationarity~\eqref{eq:S-stationary} under the assumptions of Theorem~\ref{theorem:penalisation_convergence}~(ii), consider a biactive index $i\in I^{GH}$.
Since $\multG{^t_i},\multH{^t_i}\geq0$, \eqref{eq:mu_definition} gives $\multG[hat]{_i^t}\geq-\penaltyparam{}^tH_i(\zzz^t)$ and $\multH[hat]{_i^t}\geq-\penaltyparam{}^tG_i(\zzz^t)$. Under the assumption that $\penaltyparam{}^tG_i(\zzz^t)\rightarrow0$ and  $\penaltyparam{}^tH_i(\zzz^t)\rightarrow0$, we therefore have $\liminf_{t\rightarrow\infty}\multG[hat]{_i^t}\geq0$ and $\liminf_{t\rightarrow\infty}\multH[hat]{_i^t}\geq0$, and thus $\multG[hat]{_i^*}\geq 0$ and $\multH[hat]{_i^*}\geq 0$.
\end{proof}

Recall that Ralph and Wright \cite[Theorem 5.2]{Ralph2004} provide a nice proof for penalisation with a fixed parameter $\penaltyparam{}$.  However, in practice, a user would not know a good value for $\penaltyparam{}$ ahead of time. 
If the chosen $\penaltyparam{}$ is too small, the property $G(\zzz^*)^{\top} H(\zzz^*)=0$ may not hold. 
On the other hand, if it is too large, the program may become ill-conditioned \cite[Example 3]{Leyffer2006interior} as is common with penalisation-based algorithms \cite{Murray1971,Wright1994}. 
This provides a motivation for the sequential approach.
Also, observe that the paper by Leyffer et al.~\cite[Theorem 3.4]{Leyffer2006interior} already presents a sequential penalisation theory in the case where the subproblem~\eqref{eq:penalisation} is specifically solved with the interior-point algorithm.
Theorem~\ref{theorem:penalisation_convergence} provides two useful generalisations over their work. Firstly, it gives flexibility to the user in the selection of the algorithm for the penalised subproblem \eqref{eq:penalisation}.
Secondly, we only assume
the \MPCCMFCQC{} in part~(i) and the \MPCCMFCQR{}
in part~(ii), which
are far weaker CQs compared to the \MPCCLICQ{}, required in \cite{Leyffer2006interior}.
Table~\ref{table:comparison} compares the full assumptions of Theorem~\ref{theorem:penalisation_convergence}, and of Theorem~\ref{theorem:relaxation_convergence} below, with those of the closest results in the literature.
% ...
It is also important to note all the aforementioned convergence results for penalisation, including ours, require the fulfilment of firstly $G(\zzz^*)^{\top}  H(\zzz^*) = 0$, and secondly%
, for S-stationarity, $\pi^t G_i(\zzz^t)\rightarrow0$ and $\pi^t H_i(\zzz^t)\rightarrow0$ for biactive indices $i\in I^{GH}$. 
If
the penalty parameter is held fixed, $\pi^t\equiv\pi$, as in the exact penalisation of~\cite{Ralph2004}, then the latter assumption is automatically satisfied%
, since $\pi G_i(\zzz^t)\rightarrow\pi G_i(\zzz^*)=0$ and $\pi H_i(\zzz^t)\rightarrow0$ for $i\in I^{GH}$.
Note that the proof of Theorem~\ref{theorem:penalisation_convergence} does not use the monotonicity of $(\pi^t)_{t\in\N}$, and so applies verbatim in this case.
In this case, we have shown that one can weaken the
\MPCCMFCQC{} assumption
of part~(i) to \MPCCMFCQR{}%
, and obtain S-stationarity.
However, it is unclear how to ensure that $G(\zzz^*)^{\top}  H(\zzz^*) = 0$ in practice. In the remainder of this section, we provide Theorem~\ref{theorem:penalisation_global_iterates} and Example~\ref{example:example_poly} to establish when this assumption may and may not be reasonable.

\begin{table}[t]
% ...
\caption{%
Assumptions of Theorems~\ref{theorem:penalisation_convergence} and~\ref{theorem:relaxation_convergence} and of the closest convergence results in the literature.%
}
\label{table:comparison}
\setlength{\tabcolsep}{2.5pt}
\begin{tabular}{@{}>{
\raggedright\arraybackslash}p{3.05cm}           
*{3}{>{\centering\arraybackslash}p{1.25cm}}     
>{\centering\arraybackslash}p{1.65cm}           
*{2}{>{\centering\arraybackslash}p{1.25cm}}     
>{\centering\arraybackslash}p{1.4cm}            
>{\centering\arraybackslash}p{1.65cm}@{}        
}
% ...
\toprule
& \multicolumn{4}{c}{Penalisation} & \multicolumn{4}{c}{Relaxation}\\
\cmidrule(lr){2-5}\cmidrule(l){6-9}
& Ralph \& Wright \cite{Ralph2004}
& Hu \& Ralph \cite{Hu2004}
& Leyffer et~al. \cite{Leyffer2006interior}
& This paper, Theorem~\ref{theorem:penalisation_convergence}
& Scholtes \cite{Scholtes2001}
& Schwartz \cite{schwartz2011mathematical}
& Kanzow~\& Schwartz \cite{Kanzow2015}
& This paper, Theorem~\ref{theorem:relaxation_convergence}\\
\midrule
% ...
\multicolumn{9}{@{}l}{\textbf{Setting}}\\
\quad{}$\zzz^*$ is MPCC feasible:
& assumed & assumed & implied & assumed$^{b}$ & implied & implied & implied & implied\\
\quad{}Subproblem solver:
& any & any & IPM & any & any & any & any & any\\
\quad{}Subproblem inexact:
& -- & -- & \checkmark & \checkmark & -- & -- & \checkmark & \checkmark\\
\quad{}Tolerance condition:
& -- & -- & $\epsilon^t\rightarrow0$ & $\epsilon^t\rightarrow0$ & -- & -- & $\epsilon^t{=}\,o(\relaxparam{}^t)$ & $\epsilon^t\rightarrow0$\\

% ...
\addlinespace
\multicolumn{9}{@{}l}{\textbf{C-stationarity of $\zzz^*$ under}}\\
\quad{}\MPCCLICQ{}: & \multirow{4}{=}{\centering $^{a}$} & \checkmark & \checkmark & \checkmark & \checkmark & \checkmark & \checkmark & \checkmark\\
\quad{}\MPCCMFCQT{}:& & -- & -- & \checkmark & -- & \checkmark & \checkmark & \checkmark\\
\quad{}\MPCCMFCQC{}:& & -- & -- & \checkmark & -- & -- & -- & \checkmark\\
\quad{}\MPCCMFCQR{}:& & -- & -- & -- & -- & -- & -- & -- \\
\addlinespace
% ...
\addlinespace
\multicolumn{9}{@{}l}{\textbf{S-stationarity of $\zzz^*$ under}}\\
\quad{}\MPCCLICQ{}:& \multirow{5}{=}{\centering $^{a}$} & \checkmark & \checkmark & \checkmark & \checkmark & -- & -- & \checkmark\\
\quad{}\MPCCMFCQT{}:& & -- & -- & \checkmark & -- & -- & -- & \checkmark\\
\quad{}\MPCCMFCQC{}:& & -- & -- & \checkmark & -- & -- & -- & \checkmark\\
\quad{}\MPCCMFCQR{}:& & -- & -- & \checkmark & -- & -- & -- & \checkmark \\
\multicolumn{2}{@{}l}{\quad{}With further assumptions:}
& SONC, ULSC$^{c}$ & decay$^{d}$ & decay$^{d}$ & SONC, ULSC$^{c}$ & -- & -- & decay$^{d}$\\
\botrule
\end{tabular}
\footnotetext{\checkmark: the result holds in this setting or under this assumption. 
The ticks for the CQs follow from the implications of Theorem~\ref{theorem:cq_implications}. 
IPM stands for an interior-point method. 
The same result of~\cite{schwartz2011mathematical} also appears in~\cite{Hoheisel2013}.}
\footnotetext[a]{Ralph and Wright consider a fixed parameter $\pi$ and a single KKT point of~\penalisation{\penaltyparam{}} that is feasible for~\MPCC{}, rather than a limit of iterates, and show that it is S-stationary without any CQ. They use the \MPCCMFCQR{} to show that~\penalisation{\penaltyparam{}} satisfies the \MFCQ{}~\cite[Theorem~5.1]{Ralph2004}.}
\footnotetext[b]{Theorem~\ref{theorem:penalisation_global_iterates} gives a sufficient condition.}
\footnotetext[c]{A (weak) second-order necessary condition (SONC) at each $\zzz^t$, which gives M-stationarity, and upper-level strict complementarity (ULSC) at $\zzz^*$ to get S-stationarity.}
\footnotetext[d]{Decay refers to the assumption that for the indices $i\in I^{GH}(\zzz^*)$ we have $\penaltyparam{}^tG_i(\zzz^t)\rightarrow0$ and $\penaltyparam{}^tH_i(\zzz^t)\rightarrow0$ for the penalisation or similarly $\multGH{_i^t}G_i(\zzz^t)\rightarrow0$ and $\multGH{_i^t}H_i(\zzz^t)\rightarrow0$ for the relaxation.
Note that for a fixed penalty parameter, this holds automatically.}
\end{table}

\begin{theorem}
\leavevmode\label{theorem:penalisation_global_iterates}%
Let $\left(\zzz^t\right)_{t\in\N}$ be a sequence of iterates generated by Algorithm~\ref{alg:penalisation}
with $\pi^t\rightarrow\infty$ such that $\zzz^t$ is a global optimal solution of problem \penalisation{\pi^t} for $t\in \N$%
, and suppose that problem \MPCC{} has a feasible point.
If $\underset{t \rightarrow \infty}{\lim} \zzz^t = \zzz^*$, then $G(\zzz^*)^{\top} H(\zzz^*)=0$.
\end{theorem}
\ifsuplementaryincluded
The proof is given in supplementary material Subsection~\ref{SM:proof_global}.
\else
The proof is given in the supplementary material~\cite[Subsection~\ref{SM:proof_global}]{Ward2025}.  
\fi

Next, we provide an example, which shows when the assumption that $G(\zzz^*)^{\top} H(\zzz^*)=0$ fails.
In this example, the objective and constraint functions $f,G,H:\R^2\rightarrow\R$ are all continuously differentiable, and the  \MPCCLICQ{} is satisfied. 
Nevertheless, its corresponding penalisation formulation has a sequence of local minima that converge to an infeasible point for the corresponding \MPCC{}. 

\begin{example} Consider the following example of MPCC: 
\label{example:example_poly}
\begin{equation}
\label{eq:example_poly}
\begin{aligned}     
\minimise_{z_1, z_2} \quad
&
f(\zzz)  \coloneqq  (z_1 - 3)^2 + (z_2 - 3)^2
&
\\
\subjectto \quad
& G(\zzz) \coloneqq \left(\frac{1}{3}z_1^3-\frac{9}{4}z_1^2+\frac{9}{2}z_1\right) \geq 0,\\
& H(\zzz) \coloneqq \left(\frac{1}{3}z_2^3-\frac{9}{4}z_2^2+\frac{9}{2}z_2\right) \geq 0,\\
& G(\zzz)H(\zzz)=0.
\end{aligned}
\end{equation}
\end{example}

\begin{figure}
    \centering
    \includegraphics[width=.7\textwidth]{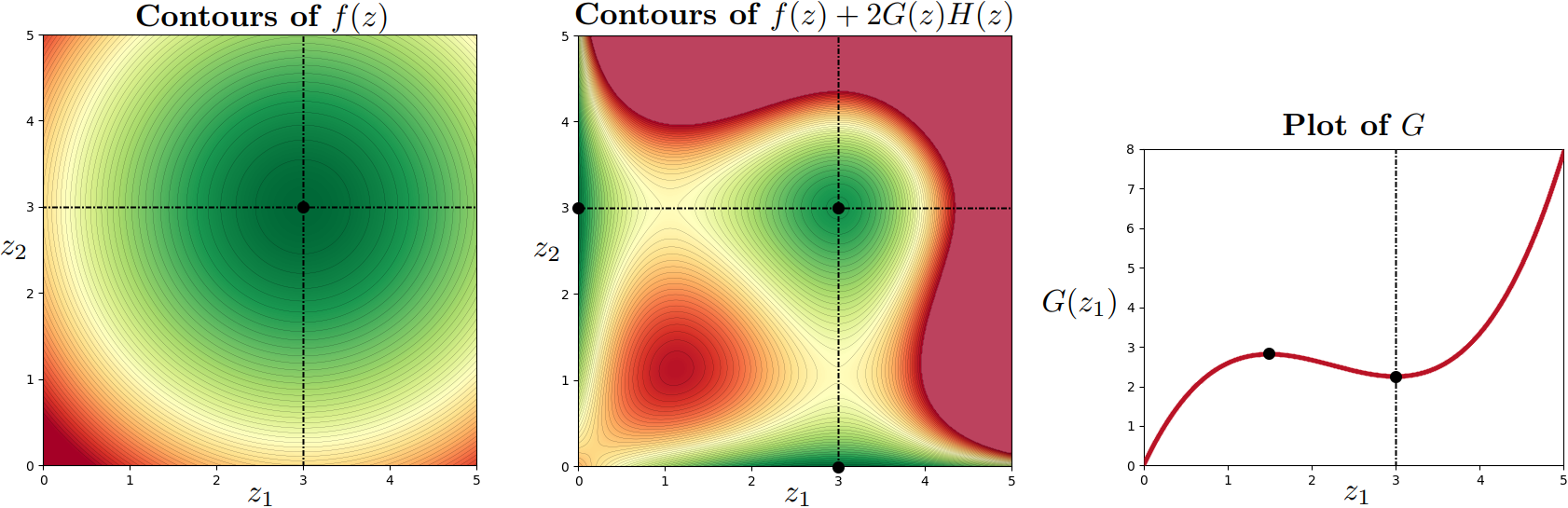}
    \caption{A graphical representation of Example~\ref{example:example_poly}. On the left are the contours of the objective function $f(\zzz)$. When the complementarity  constraints $0\leq G(\zzz)\bot H(\zzz) \geq 0$ are added, the feasible set becomes the positive vertical and horizontal axes. In the centre are the contours of the penalised objective function $f(\zzz)+\penaltyparam{} G(\zzz)H(\zzz)$. The two dark green wells at (0,3) and (3,0) are the two global minima. The light green well at (3,3) is a local minimum. On the right is the constraint $G(\zzz)$ plotted against $z_1$.} 
    \label{fig:example_poly}
\end{figure}
It is helpful to note that $G$ or $H$ are only equal to zero for $z_1=0$ or $z_2=0$, respectively. The feasible set is, therefore, the union of the two non-negative axes; see Figure~\ref{fig:example_poly}. 
We begin by writing the derivatives of all the functions involved in the above problem:
\begin{align*}
&\nabla f(\zzz) =
\begin{bmatrix}
    2z_1-6 \\
    2z_2-6 
\end{bmatrix},&
&\nabla G(\zzz) =
\begin{bmatrix}
    (z_1-\frac{3}{2})(z_1-3) \\
    0
\end{bmatrix},&
&\nabla H(\zzz) =
\begin{bmatrix}
    0 \\
    (z_2-\frac{3}{2})(z_2-3)
\end{bmatrix}.&
\end{align*}
Program \eqref{eq:example_poly} has two local minima at $[3,0]^{\top}$ and $[0,3]^{\top}$.
Furthermore, it satisfies \MPCCLICQ{} at all feasible points. This can be verified by checking the two strictly complementary  cases of 
$(G(\zzz)=0,H(\zzz)>0)$ and $(G(\zzz)>0,H(\zzz)=0)$ and the biactive case $(G(\zzz)=0,H(\zzz)=0)$. For each of these, in order, the set of gradients
of the active constraints are
$\left\{\begin{bmatrix}\frac{9}{2}, 0\end{bmatrix}^{\top}\right\}$
and
$\left\{\begin{bmatrix}0, \frac{9}{2}\end{bmatrix}^{\top}\right\}$
and 
$\left\{\begin{bmatrix}\frac{9}{2}, 0\end{bmatrix}^{\top}, \begin{bmatrix}0, \frac{9}{2}\end{bmatrix}^{\top}\right\}$ which are each linearly independent sets.

Now consider the penalisation formulation \eqref{eq:penalisation} of this example whose feasible set is the whole non-negative quadrant.
Note \eqref{eq:example_poly} has only complementarity  constraint functions $G,H$ and no inequality $g$ nor equality $h$ constraint functions. 
We introduce multipliers $\multG{_{}}$ and $\multH{_{}}$ corresponding to $G$ and $H$ respectively and write the KKT conditions of the corresponding problem \eqref{eq:penalisation}: 
\begin{equation}
    \label{eq:example_poly_kkt}
\begin{aligned}
&
    0=\nabla \left[
    f(\zzz) + \penaltyparam{}^t G(\zzz)H(\zzz)
    \right]
    - \multG{_{}}  \nabla  G(\zzz)
    - \multH{_{}} \nabla  H(\zzz), \\
&
    G(\zzz) \geq 0, \quad
    H(\zzz) \geq 0, \quad
    \multG{_{}} \geq 0, \quad
    \multH{_{}} \geq 0, \quad
    G(\zzz) \multG{_{}} = 0, \quad
    H(\zzz) \multH{_{}} = 0.
\end{aligned}
\end{equation}
Consider the point $(z_1^*,z_2^*)=(3,3)$, feasible for the corresponding version of problem \eqref{eq:penalisation} for our example but not for problem  \eqref{eq:example_poly} itself. At this point $G(\zzz^*)=H(\zzz^*)=2.25$ and all the derivatives are zero. $\nabla f(\zzz^*) = \nabla G(\zzz^*) = \nabla H(\zzz^*) = 0$.  Substituting these values into the KKT conditions for the penalisation problem (\ref{eq:example_poly_kkt}) we can find that $(z_1^*,z_2^*, \multG{_{}}, \multH{_{}})=(3,3,0,0)$ is a KKT point of the penalisation problem \penalisation{\penaltyparam{}^t} for all parameters $\penaltyparam{}^{t} \geq0$. Algorithm~\ref{alg:penalisation} could reasonably find the sequence of iterates $\zzz^t=(3,3)$ for each $t\in\N$. However, this converges to a point that is infeasible for the original MPCC. This highlights a key weakness of the penalisation method.

%% file: Article/article_32_relaxation.tex
\subsection{Relaxation}
\label{section:Relaxation}
% ...
Scholtes \cite{Scholtes2001} proposed the first regularisation scheme where problem \MPCC{} is approximated by a sequence of parametrised nonlinear programs. Since then, many alternative relaxation schemes have been proposed (e.g., \cite{Demiguel2005,Steffensen2010}), and a good survey on the subject is given by Hoheisel et al. \cite{Hoheisel2013}. Typically, in relaxation schemes for MPCCs, the constraint involving the problematic product term $G(\zzz)^\top H(\zzz)$ of the complementarity system is relaxed.  This enlarges the feasible set and can make the resulting problem more tractable. In the specific context of the Scholtes relaxation, the following subproblem is considered: 
% ...
\begin{equation}
\label{eq:relaxation}
\tag{\mbox{$\mathcal{R}_\relaxparam{}$}} 
\begin{aligned}
\minimise_{\zzz} \quad
&f(\zzz)&
\\
\subjectto \quad
&g_i(\zzz)\geq 0&
&\text{ for }i=1,\dots,p,&\\
&h_i(\zzz) = 0&
&\text{ for }i=1,\dots,q,&\\
&G_i(\zzz)\geq 0,\
H_i(\zzz)\geq 0, \
G_i(\zzz) H_i(\zzz) \leq \relaxparam{}&&\text{ for }i=1,\dots,r.&
\end{aligned}
\end{equation}
Note that for $\relaxparam{}$ close to zero, the feasible set of subproblem \relaxation{\relaxparam{}} closely approximates that of  \MPCC{}.
\ifsuplementaryincluded
This is illustrated in Figure~\ref{fig:relaxation_sequence}.
\fi
This forms the basis of the relaxation method.

% ...
In a similar fashion to our construction for penalisation, we shall not expect the NLP solvers to find the exact global minimum but instead terminate within a certain tolerance $\epsilon\geq0$ of a KKT point. We say a point $\zzz\in\R^n$ is an $\epsilon$-KKT point for \relaxation{\relaxparam{}} if there exists a Lagrange multiplier vector  $\left( \multg{},\multh{}, \multG{}, \multH{}, \multGH{} \right)\in\R^{(p+q+3r)}$  such that 
\begin{subequations}\label{eq:KKT_of_relaxed}
\begin{align}
&\begin{aligned}\label{eq:KKT_of_relaxed_a}
\Big\|\nabla f(\zzz)
 - \sum_{i=1}^{p} \multg{_i} \nabla  g_i(\zzz)
 - \sum_{i=1}^{q} \multh{_i} \nabla  h_i(\zzz)
 - \sum_{i=1}^{r} \multG{_i} \nabla  G_i(\zzz)
- \sum_{i=1}^{r} \multH{_i} \nabla  H_i(\zzz)&\\
 + \sum_{i=1}^{r} \multGH{_i} \left[G_i(\zzz)\nabla H_i(\zzz) + H_i(\zzz)\nabla G_i(\zzz)\right]&\Big\|\leq \epsilon,
\end{aligned}\\
&\begin{array}{r@{\quad} r@{\quad} r@{\quad} r}
g_i(\zzz) \geq 0, & \multg{_i} \geq 0,       & g_i(\zzz) \multg{_i}=0, & \text{ for }i=1,\dots,p,\\
h_i(\zzz) =    0, & \multh{_i} \text{ free}, &                      & \text{ for }i=1,\dots,q,\\
G_i(\zzz) \geq 0, & \multG{_i} \geq 0,       & G_i(\zzz) \multG{_i} =0, & \text{ for }i=1,\dots,r,\\
H_i(\zzz) \geq 0, & \multH{_i} \geq 0,       & H_i(\zzz) \multH{_i}=0, & \text{ for }i=1,\dots,r,\\
\left[\relaxparam{}-G_i(\zzz)H_i(\zzz)\right]\geq 0, &\multGH{_i} \geq 0, &\left[\relaxparam{}-G_i(\zzz)H_i(\zzz)\right]\multGH{_i}= 0, &\text{ for }i=1,\dots,r.
\label{eq:KKT_of_relaxed_b}
\end{array}
\end{align}
\end{subequations}

In the following we prove that a sequence of such $\epsilon$-KKT points generated by Algorithm~\ref{alg:relxation}, under certain assumptions may converge to a C-stationary point%
, or to an S-stationary point, of the original \MPCC{} problem.

\renewcommand{\algorithmicrequire}{\textbf{Input:}}
\renewcommand{\algorithmicensure}{\textbf{Procedure:}}
\begin{algorithm}[H]
\caption{Inexact sequential Scholtes relaxation}
\label{alg:relxation}
\begin{algorithmic}[1]

\Require{} Access to the \MPCC{}'s functions $f,g,h,G,H$ and their first order derivatives; Starting point $\zzz^{0}\in \R^n$;
Initial relaxation parameter $\relaxparam{}^{1}>0$;
Initial KKT tolerance $\epsilon^1 > 0$;
Exit tolerance $\epsilon_{\text{exit}} > 0$;
\vspace{.5em}

\Ensure
\Repeat{ for $t=1,2,3,\dots$:}
    \State Find an $\epsilon^t$-KKT point $\zzz^{t}$ of \relaxation{\relaxparam{}^{t}} with starting point $\zzz^{t-1}$;
    \State Decrease the relaxation parameter $\relaxparam{}^{t+1} < \relaxparam{}^t$%
    , ensuring $\relaxparam{}^{t}\rightarrow0$;
    \State Decrease the KKT tolerance $\epsilon^{t+1} < \epsilon^t$%
    , ensuring $\epsilon^{t}\rightarrow0$;
\Until{Stopping condition: $\max_{i}\left( \min(G_i(\zzz^t), \;H_i(\zzz^t)) \right)\leq \epsilon_{\text{exit}}$.}
\State Output $\zzz^t$.
\end{algorithmic}
\end{algorithm}

\begin{theorem}
\leavevmode\label{theorem:relaxation_convergence}%
Consider problem \MPCC{}. Let $\left(\zzz^{t}\right)_{t\in\N}$ be a sequence of $\epsilon^t$-KKT points of \relaxation{\relaxparam{}^t} parametrised by a strictly decreasing sequence of positive relaxation parameters $\left(\relaxparam{}^{t}\right)_{t\in\N}$ and tolerances $\left(\epsilon^{t}\right)_{t\in\N}$ that both tend to zero%
, and let $\multGH{^t}$ denote the corresponding multipliers of the constraints $G_i(\zzz)H_i(\zzz)\leq\relaxparam{}^t$. Let $\zzz^*:=\lim_{t\rightarrow\infty}\left(\zzz^{t}\right)$ exist.
Then
\begin{enumerate}[label=(\roman*),leftmargin=2em]
    \item If \MPCCMFCQC{} holds at $\zzz^*$, then $\zzz^*$ is a C-stationary point for problem \MPCC{}.
    \item If \MPCCMFCQR{} holds at $\zzz^*$ but additionally $\multGH{_i^t}G_i(\zzz^t)\rightarrow0$ and $\multGH{_i^t}H_i(\zzz^t)\rightarrow0$ for the biactive indices $i\in I^{GH}(\zzz^*)$, then $\zzz^*$ is an S-stationary point for problem \MPCC{}.
\end{enumerate}
\end{theorem}
% ...
\begin{proof}
For each iterate $\zzz^{t}$, since it is an $\epsilon^t$-KKT point of \relaxation{\relaxparam{}^t}, there exists a vector of Lagrange multipliers   
$\left( \multg{^t},\multh{^t}, \multG{^t}, \multH{^t}, \multGH{^t} \right)\in\R^{(p+q+3r)}$ 
satisfying (\link{eq:KKT_of_relaxed}{\ref*{eq:KKT_of_relaxed}a-b}). Consider the following constructions:
\begin{align*}
    &\multG[hat]{_i^t} \coloneqq 
    \begin{cases}
        \multG{_i^t}                 & \text{ if } i \in I^{G}(\zzz^t) \cup I^{GH}(\zzz^t),\\
        -\multGH{_i^t} H_i(\zzz^t)    & \text{ if } i \in \left\{i: G_i(\zzz^t)H_i(\zzz^t)=\relaxparam{}^t\right\},\\
        0                       & \text{ otherwise},
    \end{cases}&\\
    &\multH[hat]{_i^t} \coloneqq 
    \begin{cases}
        \multH{_i^t}                 & \text{ if } i \in I^{H}(\zzz^t) \cup I^{GH}(\zzz^t),\\
        -\multGH{_i^t} G_i(\zzz^t)    & \text{ if } i \in \left\{i: G_i(\zzz^t)H_i(\zzz^t)=\relaxparam{}^t\right\},\\
        0    & \text{ otherwise}.
    \end{cases}&
\end{align*}
This allows us to rewrite equation (\ref{eq:KKT_of_relaxed_a}) as
\begin{equation}
\label{eq:KKT_of_relaxed_a_rewritten}
\left\|\nabla f(\zzz^t)
 - \sum_{i=1}^{p} \multg{_i^t} \nabla  g_i(\zzz^t)
 - \sum_{i=1}^{q} \multh{_i^t} \nabla  h_i(\zzz^t)
 - \sum_{i=1}^{r} \multG[hat]{_i^t} \nabla  G_i(\zzz^t)
 - \sum_{i=1}^{r} \multH[hat]{_i^t} \nabla  H_i(\zzz^t)\right\|\leq \epsilon^t.
\end{equation}
We also note that for each $i\in \left\{i: G_i(\zzz^t)H_i(\zzz^t)=\relaxparam{}^t\right\}$, both $\multG[hat]{_i^t}$ and $\multH[hat]{_i^t}$ are non-positive, while for each $i \in I^{G}(\zzz^t) \cup I^{H}(\zzz^t) \cup{} I^{GH}(\zzz^t)$, both $\multG[hat]{_i^t}$ and $\multH[hat]{_i^t}$ are non-negative. For all other indices, they are zero.  Therefore, we can conclude that 
\begin{equation}
    \label{eq:r_C_stationarity}
    \multG[hat]{_i^t}\multH[hat]{_i^t}\geq 0 \quad \forall i\in \{1,...,r\}.
\end{equation}

We now wish to prove that the sequence $\aaa_t \coloneqq \left( \multg{^t},\multh{^t}, \multG[hat]{^t}, \multH[hat]{^t} \right)_{t\in\N}$ has a convergent subsequence. Suppose towards a contradiction that it does not. It must be that $\aaa_t$ is unbounded and $\|\aaa_t\|\rightarrow\infty$. Dividing equation  \eqref{eq:KKT_of_relaxed_a_rewritten} through by $\|\aaa_t\|$ leads to 
\begin{multline}
\label{eq:KKT_a_divided_by_at}
\left\|
\frac{\nabla f(\zzz^t)}{\|\aaa_t\|}
 - \sum_{i=1}^{p} \frac{\multg{_i^t}}{\|\aaa_t\|}      \nabla  g_i(\zzz^t)
 - \sum_{i=1}^{q} \frac{\multh{_i^t}}{\|\aaa_t\|}      \nabla  h_i(\zzz^t)
 - \sum_{i=1}^{r} \frac{\multG[hat]{_i^t}}{\|\aaa_t\|} \nabla  G_i(\zzz^t)
 - \sum_{i=1}^{r} \frac{\multH[hat]{_i^t}}{\|\aaa_t\|} \nabla  H_i(\zzz^t)
 \right\| \leq
 \frac{\epsilon^t}{\|\aaa_t\|}.
\end{multline}

Now consider the normed sequence 
$\frac{
\left(\multg{^t},\multh{^t}, \multG[hat]{^t}, \multH[hat]{^t}\right)}
{\left\|\left(\multg{^t},\multh{^t}, \multG[hat]{^t}, \multH[hat]{^t}\right)\right\|}$. Clearly, it is bounded. It follows from the Bolzano--Weierstrass theorem that there exists some subsequence that converges. Let the point $\left(\multg{'}, \multh{'}, \multG[hat]{'}, \multH[hat]{'}\right)$ be the limit of that subsequence.
Since every term of the normed sequence has unit norm, so does its limit; in particular, the limit is non-zero. From the complementarity slackness conditions in \eqref{eq:KKT_of_relaxed_b} we conclude that: $\multg{_i'}$ is non-zero only for $i\in I^g$ and $\multG[hat]{'_i}$ is non-zero only for $i\in I^G\cup I^{GH}$ and $\multH[hat]{'_i}$ is non-zero only for $i\in I^H\cup I^{GH}$.
While the claim on $\multg{'}$, and also $\multg{'}\geq0$, follows directly, the claims on $\multG[hat]{'}$ and $\multH[hat]{'}$ require a short argument, since the constraint $G_i(\zzz^t)H_i(\zzz^t)\leq\relaxparam{}^t$ may be active. Consider $i\in I^H$, so that $G_i(\zzz^*)>0$ and $H_i(\zzz^*)=0$. Then $G_i(\zzz^t)>0$ for all large $t$, so that $i\notin I^{G}(\zzz^t) \cup I^{GH}(\zzz^t)$ and, by construction, $\multG[hat]{_i^t}$ is either zero or equal to $-\multGH{_i^t}H_i(\zzz^t)$, the latter only if $G_i(\zzz^t)H_i(\zzz^t)=\relaxparam{}^t$. In the latter case $H_i(\zzz^t)>0$, so that also $\multH[hat]{_i^t}=-\multGH{_i^t}G_i(\zzz^t)$, and if $\multGH{_i^t}>0$, then, since $|\multH[hat]{_i^t}|\leq\|\aaa_t\|$,
\begin{equation*}
\frac{\bigl|\multG[hat]{_i^t}\bigr|}{\|\aaa_t\|}
\;\leq\; \frac{\bigl|\multG[hat]{_i^t}\bigr|}{\bigl|\multH[hat]{_i^t}\bigr|}
\;=\; \frac{\multGH{_i^t} H_i(\zzz^t)}{\multGH{_i^t} G_i(\zzz^t)}
\;=\; \frac{H_i(\zzz^t)}{G_i(\zzz^t)}
\;\longrightarrow\; \frac{0}{G_i(\zzz^*)}=0 .
\end{equation*}
In either case $\multG[hat]{'_i}=0$ for $i\in I^H$; exchanging the roles of $G$ and $H$ gives $\multH[hat]{'_i}=0$ for $i\in I^G$.
Using this indexing we can say that it must be the
$\left(\multg{'_{I^g}}, \multh{'}, \multG[hat]{_{I^G\cup I^{GH}}'}, \multH[hat]{_{I^H\cup I^{GH}}'}\right)$
components of the sequence that are non-zero.

Returning to \eqref{eq:KKT_a_divided_by_at}, as $t\rightarrow\infty$, the terms $\frac{\nabla f(\zzz^t)}{\|\aaa_t\|}$ and $ \frac{\epsilon^t}{\|\aaa_t\|}$ go to zero, and we are left with
\begin{equation}
\label{eq:r_linearly_dependence}
\begin{aligned}
\sum_{i\in I^g} \multg{_i'} \nabla  g_i(\zzz^*)
 + \sum_{i\in I^h} \multh{_i'} \nabla  h_i(\zzz^*)
 + \sum_{{i\in I^G\cup I^{GH}}} \multG[hat]{'_i} \nabla  G_i(\zzz^*)
 + \sum_{{i\in I^H\cup I^{GH}}} \multH[hat]{'_i} \nabla  H_i(\zzz^*)=\veczero{}.
\end{aligned}
\end{equation}
Under the assumptions of part~(i), since $\multg{'}\geq0$ and, by~\eqref{eq:r_C_stationarity}, $\multG[hat]{'_i}\multH[hat]{'_i}\geq 0$ for $i\in I^{GH}$, the non-zero vector $\left(\multg{'}, \multh{'}, \multG[hat]{'}, \multH[hat]{'}\right)$ satisfies~\eqref{eq:MPCCMFCQC_multipliers} at $\zzz^*$, which contradicts the \MPCCMFCQC{} by Lemma~\ref{lemma:positive-linearly-independent}(iv).
Under the assumptions of part~(ii), consider $i\in I^{GH}$. By construction, $\multG[hat]{_i^t}$ is equal to $\multG{_i^t}\geq 0$, to $-\multGH{_i^t}H_i(\zzz^t)$ or to zero, so that $\multG[hat]{_i^t}\geq-\multGH{_i^t}H_i(\zzz^t)\rightarrow 0$. As $\|\aaa_t\|\rightarrow\infty$, this gives $\multG[hat]{'_i}\geq 0$ and, in the same way, $\multH[hat]{'_i}\geq 0$ for $i\in I^{GH}$.
From the assumptions of
part~(ii), we have that $\zzz^*$ satisfies \MPCCMFCQR{}.  By Lemma~\ref{lemma:positive-linearly-independent}(iii), the set of gradient vectors in (\ref{eq:MFCQ_dual_form_part_3}) should be positive-linearly independent. But that is contradicted by (\ref{eq:r_linearly_dependence}).  Therefore,
in both cases, we conclude  $\aaa_t$ must have a convergent subsequence.  Let $\left(\multg{^*}, \multh{^*}, \multG[hat]{^*}, \multH[hat]{^*}\right)$ be the limit of that convergent subsequence.
By the continuous differentiability of $f, g, h, G, H$, and since the sequence $\left(\zzz^t\right)_{t\in\N}$ converges in their domain, these functions and their derivatives must converge in range. We can apply all this to
\eqref{eq:KKT_of_relaxed_a_rewritten}, recalling that $\epsilon^t\rightarrow0$, to get the following:
\begin{equation}
\label{eq:r_final_stationarity}
\begin{aligned}
\nabla f(\zzz^*)
 - \sum_{i=1}^{p} \multg{_i^*} \nabla  g_i(\zzz^*)
 - \sum_{i=1}^{q} \multh{_i^*} \nabla  h_i(\zzz^*)
 - \sum_{i=1}^{r} \multG[hat]{_i^*} \nabla  G_i(\zzz^*)
 - \sum_{i=1}^{r} \multH[hat]{_i^*} \nabla  H_i(\zzz^*) = 0.
\end{aligned}
\end{equation}
The conditions in~\eqref{eq:stationarity_b} on $\multg{^*}$ follow directly from~\eqref{eq:KKT_of_relaxed_b}. For those on $\multG[hat]{^*}$ and $\multH[hat]{^*}$, consider $i\in I^H$. As shown above, for all large $t$, $\multG[hat]{_i^t}$ is either zero or equal to $-\multGH{_i^t}H_i(\zzz^t)$, and in the latter case $\multH[hat]{_i^t}=-\multGH{_i^t}G_i(\zzz^t)$. Since $\multH[hat]{_i^t}$ converges to the finite limit $\multH[hat]{_i^*}$ along the subsequence while $G_i(\zzz^t)$ is bounded away from zero, $\multGH{_i^t}$ remains bounded for all such $t$. Consequently $\multG[hat]{_i^t}\rightarrow0$, because $H_i(\zzz^t)\rightarrow0$, and therefore $\multG[hat]{_i^*}=0$ for all $i\in I^H$. Exchanging the roles of $G$ and $H$ gives $\multH[hat]{_i^*}=0$ for all $i\in I^G$.
Therefore, under the assumptions of part~(i), from equations~\eqref{eq:r_final_stationarity}, \eqref{eq:KKT_of_relaxed_b}
together with the argument above, and \eqref{eq:r_C_stationarity}, we can conclude that~\eqref{eq:stationarity_a}, \eqref{eq:stationarity_b} and \eqref{eq:C-stationary} hold, respectively. These are exactly the C-stationary conditions.
Under the assumptions of part~(ii), we showed above that $\multG[hat]{_i^t}\geq-\multGH{_i^t}H_i(\zzz^t)\rightarrow 0$ and, in the same way, $\multH[hat]{_i^t}\geq-\multGH{_i^t}G_i(\zzz^t)\rightarrow 0$ for $i\in I^{GH}$. Hence $\multG[hat]{_i^*}\geq 0$ and $\multH[hat]{_i^*}\geq 0$ for $i\in I^{GH}$, which together with~\eqref{eq:stationarity_a} and~\eqref{eq:stationarity_b} gives~\eqref{eq:S-stationary}; that is, $\zzz^*$ is S-stationary.
\end{proof}

A similar version of Theorem~\ref{theorem:relaxation_convergence} was initially proven by Scholtes \cite[p. 922, Section 3]{Scholtes2001} but required the stronger  \MPCCLICQ{}. An updated proof is given by Schwartz \cite[p. 112, Theorem 7.3]{schwartz2011mathematical}, which weakened the requirement to \MPCCMFCQT{}.
We present a similar proof but weakening the CQ further to
\MPCCMFCQC{} in part~(i), and to \MPCCMFCQR{} under the additional assumption on the multipliers in part~(ii).
This is a significant weakening of assumptions (cf. Theorem~\ref{theorem:cq_implications} and Example~\ref{eq:Example01}).
The \MPCCMFCQR{} alone, however, is not sufficient; see
\ifsuplementaryincluded
the supplementary materials, Subsection~\ref{SM:counterexample_relaxation},
\else
the supplementary materials~\mbox{\cite[Subsection~\ref{SM:counterexample_relaxation}]{Ward2025}},  
\fi
for an example in which a sequence of global minimisers of the relaxed problems converges to a point that satisfies the \MPCCMFCQR{} but is not even W-stationary.
Another feature of our result (Theorem~\ref{theorem:relaxation_convergence}) is that we only require an inexact KKT point to be computed at each iteration, unlike in
\cite{Scholtes2001,schwartz2011mathematical}, where exact KKT points were used for the
relaxed problem.
Kanzow and Schwartz~\cite{Kanzow2015} also consider inexact KKT points of the relaxed problems, with errors in all the KKT conditions, including the signs of the multipliers, but they require the \MPCCMFCQT{} and $\epsilon^t=o(\relaxparam{}^t)$. Our $\epsilon^t$-KKT points only allow an error in the stationarity condition, but any $\epsilon^t\rightarrow0$ suffices, and the required qualification \MPCCMFCQC{} is weaker; see Table~\ref{table:comparison}.  
This allows the theorem to be applied to a much broader class of problems and for the numerical implementation to reflect the theory more closely.  

%% file: Article/article_40_application.tex
\section{Application}
\label{section:Application}
% ...
In this section, we construct an \MPCC{} model to optimally select the hyperparameters of a Support Vector Machine (SVM).
We then show, amongst other theoretical properties, that the model satisfies
\MPCCMFCQC{}, and hence \MPCCMFCQR{}%
, for any feasible point with positive hyperparameters%
, under two further mild assumptions, but that it systematically fails \MPCCLICQ{} and \MPCCMFCQT{} in general.

%% file: Article/article_41_model.tex
\subsection{Modelling hyperparameter tuning}
\label{section:model}
% ...
A major aspect of machine learning is judging if a model is overfitted to its training data \cite[p. 28]{Geron2022}, \cite[Chapter 7]{Hastie2009}.
For classifiers, this could mean unreasonably contorting the decision boundary to include outliers. 
By doing so, the classifier achieves near-zero empirical training error but misses the overarching pattern. 
SVMs' power comes from their ability to accept some level of empirical error in favour of a simpler boundary that generalises well to unseen data. 
These two objectives are balanced by the hyperparameters. 
This leads to the question of how to evaluate the accuracy of a model or the choice of hyperparameters if training error does not provide the whole picture.

The cross-validation procedure (see literature \cite[Section 7.10.1]{Hastie2009}, \cite[p. 77]{Geron2022}, \cite{Stone1974}) partitions the data indices into $K$ folds. It then runs $K$ experiment(s) such that at the $k$-th experiment ($k=1, \ldots, K$), new model parameters are learned using all the examples but withholding those from fold $k$. Then, the remaining data points from fold $k$ can be used as unseen validation data to evaluate the model's ability to generalise. This typically results in a less biased evaluation of the model.

To denote the cross-validation procedure applied to a classification dataset, let $K$ be the total number of folds.
For each experiment $k=1,\hdots, K$, the feature vectors~$\xxx^k_i\in\R^d$ and labels~$y^k_i\in\{-1,+1\}$ relating to that fold will be identified with the superscript ${k}$.
More precisely, we use the following notation for experiment $k$:
% ...
\begin{equation*}
\begin{aligned}
&&
&\text{\MakeUppercase{Training}}&
&\text{\MakeUppercase{Validation}}&\notag\\
&\text{Cardinality:}&
&n^{k}&
&\bar{n}^{k}&\\
&\text{Data/examples:}&
&\left(\xxx_i^{k}, y_i^{k}\right)_{i=1,\dots,n^k}&
&\left(\bar{\xxx}_i^{k},\bar{y}_i^{k}\right)_{i=1,\dots,\bar{n}^k}&\\
&\text{Kernel matrix:}&
&Q(\gamma)^k&
&\bar{Q}(\gamma)^k&
\end{aligned}
\end{equation*}

To represent a non-linear relationship we use the RBF kernel though our theory holds for any positive semi-definite kernel%
, except for Proposition~\ref{proposition:satisfy_mpccmfcqr}, which requires the kernel matrices to be positive definite; for the RBF kernel, this holds whenever the training inputs are distinct (Lemma~\ref{lemma:psd_kernel}).
Note that in training, the kernel 
$Q(\gamma)^k_{ij}\coloneqq{}y_i^k y_j^k\exp{\left(-\gamma\|\xxx_i^k-\xxx_j^k\|^2\right)}$ 
is calculated between two training examples $(\xxx_i^k,y_i^k)$ and $(\xxx_j^k,y_j^k)$ but, in validation, the kernel 
$\bar{Q}(\gamma)^k_{ij}\coloneqq \bar{y}_i^k y_j^k \exp{\left(-\gamma\|\bar{\xxx}_i^k-\xxx_j^k\|^2\right)}$ 
is between one validation $(\bar{\xxx}_i^k,\bar{y}^k_i)$ and one training example $(\xxx_j^k,y^k_j)$. 
The scalar $\gamma>0$ is a hyperparameter.
With this notation in mind we proceed to define the SVM hyperparameter tuning optimisation problem.

\begin{equation}
\label{eq:final-model}
\tag{\mbox{$\mathcal{M}$}} 
\begin{aligned}
\minimise_{C, \gamma, \zeta, \ablue{}, \mlower{}, \mupper{}, \uu{}} \quad
&\frac{1}{K}\sum_{k=1}^{K} \frac{1}{\bar{n}^k} \sum_{i=1}^{\bar{n}^k} \zeta_i^k&
\\
\subjectto \quad
&\left. C \geq 0, \right.& \\
&\left. \gamma \geq 0, \right.& \\
&\left.
\begin{aligned}
&\zeta_i^k \geq 0 \quad \quad \quad \\
&\zeta_i^k \geq 1-\sum_{j=1}^{n^k} \ablue{_j^k} \bar{Q}(\gamma)^k_{ij} - \bar{y}_i^k \uu{^k}\\
\end{aligned}
\hskip39px
\right\}
\begin{aligned}
&\texttt{(Validation)}\\
&\text{for }k=1,\hdots,K,\\
&\quad\text{for }{i=1,\hdots,\bar{n}^k},
\end{aligned}
\\
&
\left.
\begin{aligned}
&\sum_{j=1}^{n^k}   \ablue{^k_j} Q(\gamma)^k_{ij}  -1 - \mlower{^k_i} + \mupper{^k_i} + \uu{^k} y_i^k=0
\\
& 0 \leq \ablue{_i^k} \bot \mlower{_i^k} \geq 0
\\
& 0 \leq  (C - \ablue{_i^k}) \bot \mupper{_i^k} \geq 0 
\\
& \sum_{j=1}^{n^k} \ablue{_j^k} y_j^k = 0
\end{aligned}
\right\}
\begin{aligned}
&\texttt{(Training)}\\
&\text{for }k=1,\hdots,K,\\
&\quad\text{for } {i=1,\hdots,n^k}.
\end{aligned}
\end{aligned}
\end{equation}

\newcommand{\SVMdual}{\ifsuplementaryincluded\eqref{eq:SVMdual}\else{dual SVM}\fi}

The set of constraints tagged \texttt{(Training)} are exactly the KKT conditions of the  \SVMdual{} for hyperparameters~$C$ and~$\gamma$~\cite{Fan2005}. 
The decision variable~$\ablue{_i^k}$ is non-zero if training data point~$i$ from experiment~$k$ is a support vector and~$\mlower{_i^k}$, $\mupper{_i^k}$ and~$\uu{^k}$ are the Lagrange multipliers corresponding to its
lower bound, upper bound and bias respectively. 
Since the 
\SVMdual{}%
\ is convex with affine constraints, these constraints are only satisfied for an optimal choice of support vectors.
\ifsuplementaryincluded
For a full expos\'{e} on this system see supplementary material Subsection~\ref{SM:KKT_of_dual_SVM} and~\cite[p. 1891]{Fan2005}.
\else
For the full theory of the SVM and this KKT representation see~\cite[p. 1891]{Fan2005}.
\fi

The variable~$\zeta^k_i$ and the two constraints tagged \texttt{(Validation)} calculate the validation loss of each data point~$i$ from fold~$k$. 
The objective function then seeks to minimise the average validation loss.
In the language of bilevel optimisation this is equivalent to the KKT-reformulation of a bilevel program with lower-level programs corresponding to training $K$ dual SVMs over training datasets $k=1,\dots,K$ and an upper-level program minimising the validation loss over the choice of hyperparameters.
For full details on the link to bilevel programming see 
\ifsuplementaryincluded
supplementary material Section~\ref{SM:relationship_to_bilevel}.
\else
the supplementary material~\cite[Section~\ref{SM:relationship_to_bilevel}]{Ward2025}.  
\fi

Problem  \MODEL{} is explicit but not concise. It uses a two-dimensional indexing scheme: First, it indexes the fold $k=1,\dots, K$, then the individual examples $i$ within that fold. We now wish to flatten the indices into a one-dimensional sequence. Hence, we introduce the index sets
% ...
\begin{equation}
    \label{eq:flat_index}
    I  \coloneqq  \left\{(k,i):k=1,\dots,K, i=1,\dots,n^k\right\} \,\mbox{ and }\,
    \bar{I}  \coloneqq  \left\{(k,i):k=1,\dots,K, i=1,\dots,\bar{n}^k\right\}.
\end{equation}
% ...
This allows us to stack all the constraint functions into vectors via the notations
\begin{subequations}
\label{eq:constraint_function_definitions}
\begin{align}
    \label{eq:Z}
    &Z_i^k(\zzz)
     \coloneqq \zeta_i^k-1+\sum_{j=1}^{n^k} \ablue{_j^k} \bar{Q}(\gamma)_{ij}^k + \bar{y}_i^k \uu{^k}&
    &\mbox{ for }\;\, (k, i)\in \bar{I},&
    \\
    \label{eq:theta}
    &\theta_i^k(\zzz)
     \coloneqq \sum_{j=1}^{n^k}   \ablue{^k_j} Q(\gamma)^k_{ij} - 1 - \mlower{^k_i} + \mupper{^k_i} + \uu{^k} y_i^k&
    &\mbox{ for }\;\,{(k, i)\in I,}&
    \\
    \label{eq:phi}
    &\varphi^k(\zzz)
     \coloneqq \sum_{j=1}^{n^k} \ablue{_j^k} y_j^k&
    &\mbox{ for }\;\,{k=1,\dots,K,}&
\end{align}
\end{subequations}
where we have also introduced  $\zzz  \coloneqq  \left[C, \gamma, \zeta, \ablue{}, \mlower{}, \mupper{}, \uu{}\right]^{\top} \in \R^{N}$ to collect all variables, which has the dimension $N \coloneqq 3(K-1)n+n+K+2$%
, where $n$ denotes the total number of data points. This leads to the representation  $f(\zzz) \coloneqq  \frac{1}{K}\sum_{k=1}^{K} \frac{1}{\bar{n}^k} \sum_{i=1}^{\bar{n}^k} \zeta_i^k$ for the upper-level objective function of  \MODEL{}, as well as its constraint functions 
\begin{equation}
\label{eq:ghGH}
\begin{aligned}
&g(\zzz)  \coloneqq 
\begin{bmatrix}
    C\\
    \gamma\\
    \zeta_1^{1}\\
    \vdots \\
    \zeta_{\bar{n}^K}^K\\
    Z_1^{1}(\zzz)\\
    \vdots \\
    Z_{\bar{n}^K}^K(\zzz)\\
\end{bmatrix},&
&h(\zzz)  \coloneqq  
\begin{bmatrix}
    \theta_1^{1}(\zzz)\\
    \vdots \\
    \theta_{n^K}^K(\zzz)\\
    \varphi^1(\zzz)\\
    \vdots\\
    \varphi^K(\zzz)
\end{bmatrix},&
&G(\zzz)  \coloneqq 
\begin{bmatrix}
    \ablue{_1^1} \\
    \vdots \\
    \ablue{_{n^K}^K} \\
    \left(C - \ablue{_1^1}\right)\\
    \vdots \\
    \left(C - \ablue{_{n^K}^K}\right)
\end{bmatrix},&
&H(\zzz)  \coloneqq 
\begin{bmatrix}
    \mlower{_1^{1}} \\
    \vdots \\
    \mlower{_{n^K}^K} \\
    \mupper{_1^1}\\
    \vdots \\
    \mupper{_{n^K}^K}
    \end{bmatrix},&
\end{aligned}
\end{equation}
where we obviously have $g(\zzz)\in\R^{2+2n}, h(\zzz)\in\R^{(K-1)n+K}$ and $G(\zzz)$, $H(\zzz)\in \R^{2(K-1)n}$. 
% ...
Subsequently, using the above constructions for the functions $f$, $g$, $h$, $G$ and $H$, we can re-write problem~\MODEL{} more compactly in the general form~\MPCC{}.

%% file: Article/article_42_properties.tex

\subsection{Theoretical properties of the model}
\label{section:theoretical_properties_of_the_model}

We begin by making an important assumption that for each fold, there exists at least one dual variable in the interior of its bounds. This
implies that at least one data point lies on the margin.
\begin{equation}
    \label{eq:assumption_on_alpha}
    \forall k=1,\dots,K,\quad \exists i\in \{1,\dots,n^k\}\quad \text{such that } 0<\ablue{_i^k}<C.
\end{equation}
In practice,~\eqref{eq:assumption_on_alpha} typically holds.
Furthermore, for a fixed $\gamma>0$ and distinct training inputs, it has been proven that there exists a threshold $\bar{C}$ such that \eqref{eq:assumption_on_alpha} holds for all $C\geq\bar{C}$; see~\cite[Section~2, Case~2]{KeerthiLin2003asymptotic}.
In the proof of Proposition~\ref{proposition:satisfy_mpccmfcqr} below, assumption~\eqref{eq:assumption_on_alpha} is used to show that the multiplier associated with the bias constraint $\varphi^k(\zzz)=0$ vanishes.

We now give the derivatives of the constraint functions of \MODEL{}. This will be useful later in studying which constraint qualifications problem \MODEL{} satisfies. 
\begin{proposition}
The transposed Jacobian matrix of the constraint functions of \MODEL{} is
\begin{equation}
\label{eq:Jacobian}
\setlength{\arraycolsep}{5pt}
\begin{array}{c c | c c c c c c c c c c}
& \multicolumn{1}{c}{}& \multicolumn{10}{c}{\text{Constraint functions}}  \\
& \multicolumn{1}{c}{}
&\multicolumn{4}{c}{\overbrace{\hphantom{\begin{array}{c c c c}C&\gamma&\zeta&Z(\zzz)\end{array}}}^{g(\zzz)}}
&\multicolumn{2}{c}{\overbrace{\hphantom{\begin{array}{c c c c}\nabla_\gamma \theta(\zzz)&\varphi{}(\zzz)\end{array}}}^{h(\zzz)}}
&\multicolumn{2}{c}{\overbrace{\hphantom{\begin{array}{c c c c}\ablue{}&C-\ablue{}\end{array}}}^{G(\zzz)}}
&\multicolumn{2}{c}{\overbrace{\hphantom{\begin{array}{c c c c}\mlower{}&\mupper{}\end{array}}}^{H(\zzz)}}
\\
& & C &\gamma &\zeta & Z(\zzz) & \theta{(\zzz)} & \varphi{(\zzz)} & \ablue{} & (C-\ablue{}) & \mlower{} & \mupper{}
\\\cline{2-12}
\parbox[t]{1em}{\multirow{7 }{*}{\rotatebox[origin=c]{90}{\text{Derivative w.r.t.}}}}   
&\nabla_{C}        \bullet&1&0&0  &\nabla_C Z(\zzz)     &0                     &0 &0  &1   &0  &0  \\
&\nabla_{\gamma}   \bullet&0&1&0  &\nabla_\gamma Z(\zzz)&\nabla_\gamma\theta(\zzz)&0 &0  &0   &0  &0  \\
&\nabla_{\zeta}    \bullet&0&0&\id&\id               &0                     &0 &0  &0   &0  &0  \\
&\nabla_{\ablue{}} \bullet&0&0&0  &\nabla_\alpha Z(\zzz)&\Q(\gamma)            &\Y&\id&-\id&0  &0  \\
&\nabla_{\mlower{}}\bullet&0&0&0  &0                 &-\id                  &0 &0  &0   &\id&0  \\
&\nabla_{\mupper{}}\bullet&0&0&0  &0                 &\id                   &0 &0  &0   &0  &\id\\
&\nabla_{\uu{}}    \bullet&0&0&0  &\bar{\mathcal{Y}}^{\top}&\mathcal{Y}^{\top}         &0 &0  &0   &0  &0  \\
\end{array}
\end{equation}
with the column headers denoting which constraint function (refer to \eqref{eq:ghGH}) the column corresponds to; the row headers $\nabla_v\bullet$ denote that their row corresponds to the derivative with respect to decision variable $v$; $\id$ is the identity matrix with inferred dimension; $\Q\in\R^{|I|\times|I|}$ and $\mathcal{Y}\in\{-1,0,1\}^{|I|\times K}$ are block diagonal matrices formed by 
\begin{align*}
&\Q(\gamma)  \coloneqq  
    \begin{bmatrix}
        Q(\gamma)^1 & & 0\\
        & \ddots & \\
        0 &  & Q(\gamma)^K\\
\end{bmatrix},\;&
&\mathcal{Y} \coloneqq 
\begin{bmatrix}
    \boldsymbol{y}^1 & & \veczero{} \\
     & \ddots & \\
    \veczero{} & & \boldsymbol{y}^K
\end{bmatrix}
\text{ and }&
&\boldsymbol{y}^k  \coloneqq 
\begin{bmatrix}
    y_1^k \\
    \vdots \\
    y_{n^k}^k
\end{bmatrix} 
\text{ for }
k=1,\dots,K.&
\end{align*}
\end{proposition}

\begin{figure}
    \centering
    \includegraphics[width=.6\textwidth]{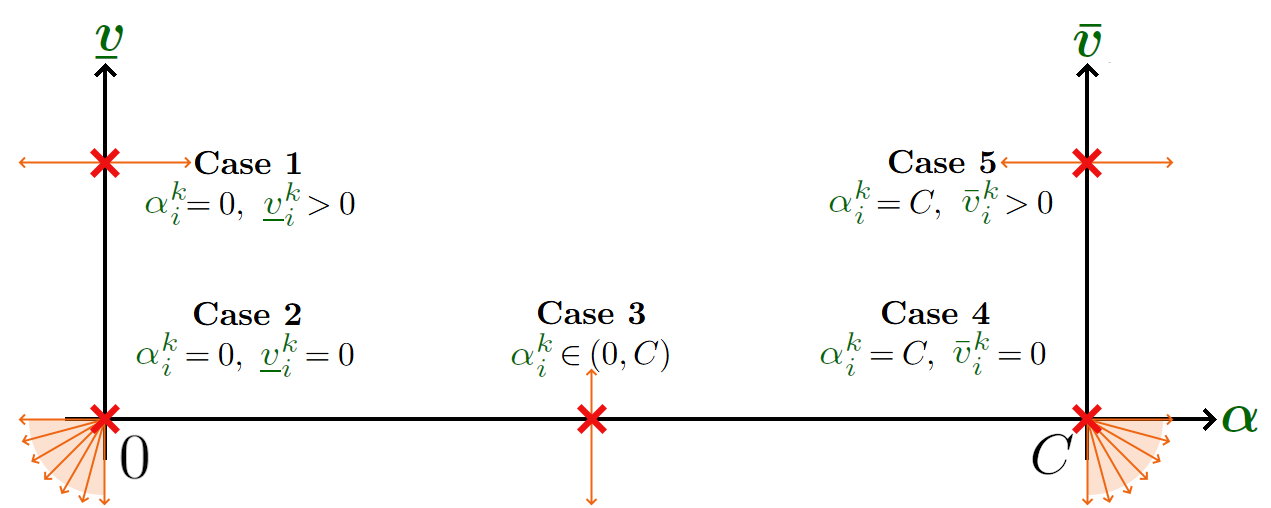}
    \caption{This figure illustrates the geometry of our problem by showing example feasible points for each of the five index sets $C_1$, $C_2$, $C_3$, $C_4$ and $C_5$. The horizontal axis represents values of $\ablue{_i^k}$. The vertical axis represents non-negative values of {\color{dark_green}\underbar{$v$}${_i^k}$} and $\mupper{_i^k}$ when $\ablue{_i^k}=0$ and $\ablue{_i^k}=C$, respectively. The orange arrows show directions in the normal cone of the feasible set of the problem.}
    \label{fig:five_cases}
\end{figure}

Using the standard MPCC index sets $I^g, I^G, I^H, I^{GH}$ that were defined in (\ref{eq:index_sets}) becomes laborious and hides the structure of our problem. Recalling the indexing scheme of \eqref{eq:flat_index}, we define the problem-specific index sets
\[  \bar{I}^\zeta  \coloneqq  \left\{(k,i)\in \bar{I}:\zeta_i^k=0 \right\} \;  \mbox{ and } \;
\bar{I}^Z  \coloneqq      \left\{(k,i)\in \bar{I}:Z_i^k(\zzz)=0 \right\}. 
\]
We then partition the index set $I$ into five cases that better describe \MODEL{}:
\begin{equation*}
\begin{aligned}
    &C_1  \coloneqq  \left\{(k,i)\in I:\quad \ablue{_i^k}=0,\quad\quad \mlower{_i^{k}}>0,\quad (C-\ablue{_i^k})=C,\quad\quad \mupper{_i^k}=0\right\};\\
    &C_2  \coloneqq  \left\{(k,i)\in I:\quad \ablue{_i^k}=0,\quad\quad \mlower{_i^{k}}=0,\quad (C-\ablue{_i^k})=C,\quad\quad \mupper{_i^k}=0\right\};\\
    &C_3  \coloneqq  \left\{(k,i)\in I:\quad \ablue{_i^k}\in(0,C),\    \mlower{_i^{k}}=0,\quad   (C-\ablue{_i^k})\in(0,C),\ \mupper{_i^k}=0\right\};\\
    &C_4  \coloneqq  \left\{(k,i)\in I:\quad \ablue{_i^k}=C,\quad\quad \mlower{_i^{k}}=0,\quad (C-\ablue{_i^k})=0,\quad\quad \mupper{_i^k}=0\right\};\\
    &C_5  \coloneqq  \left\{(k,i)\in I:\quad \ablue{_i^k}=C,\quad\quad \mlower{_i^{k}}=0,\quad (C-\ablue{_i^k})=0,\quad\quad \mupper{_i^k}>0\right\}.
\end{aligned}
\end{equation*}
We denote the union of multiple cases in the subscript. For example, $C_{123} \coloneqq C_1 \cup C_2 \cup C_3$. It is easy to see that these sets are disjoint and cover all possible cases. That is, $C_{12345}=I$.  Figure~\ref{fig:five_cases} illustrates the geometry of these index sets.
Next, we show which constraint qualification of Section~\ref{section:constraint_qualification} \MODEL{} satisfies.
\begin{proposition}
\leavevmode\label{proposition:satisfy_mpccmfcqr}%
Let $\zzz \coloneqq [C, \gamma, \zeta, \ablue{}, \mlower{}, \mupper{}, \uu{}]^\top$ be any feasible point of problem \MODEL{} such that $C>0$, $\gamma>0$ and \eqref{eq:assumption_on_alpha} holds%
, and suppose that, for each fold $k$, the training inputs $\xxx_1^k,\dots,\xxx_{n^k}^k$ are distinct. Then
\MPCCMFCQC{}, and hence \MPCCMFCQR{}%
, holds at this point.
\end{proposition}

\begin{proof}

Recall from Lemma~\ref{lemma:positive-linearly-independent}%
(iv) that \MPCC{} satisfies
\MPCCMFCQC{} if and only if there are no scalars, not all of them zero, satisfying \eqref{eq:MPCCMFCQC_multipliers}, which involves the family of gradients $(A,B)$ defined in \eqref{eq:MFCQ_dual_form_part_3}. We write these gradients explicitly for \MODEL{} using the problem specific variables, index sets $\bar{I}^\zeta$, $\bar{I}^Z$, and partition $C_l$ for $l=1, \ldots, 5$: 
\begin{subequations}
\begin{multline}
\label{eq:set_A}
A  \coloneqq  
     \Big\{\derivative{\zzz} \ablue{_i^k},     \derivative{\zzz} \mlower{_i^k}: {(k,i)\in{C_2}}\Big\}
\cup \Big\{\derivative{\zzz} (C-\ablue{_i^k}), \derivative{\zzz} \mupper{_i^k}: {(k,i)\in{C_4}} \Big\}\\
\cup \Big\{\derivative{\zzz} \zeta_i^k:     {(k,i)\in \bar{I}^\zeta}  \Big\} 
\cup \Big\{\derivative{\zzz} Z_i^k(\zzz):         {(k,i)\in \bar{I}^Z}      \Big\};
\end{multline}
\vspace{-.5em}
\begin{multline}
\label{eq:set_B}
B \coloneqq  
    \Big\{\derivative{\zzz} \ablue{_i^k}: {(k,i)\in{C_1}}\Big\}
\cup\Big\{\derivative{\zzz} \mlower{_i^k}: {(k,i)\in{C_{345}}}\Big\} 
\cup\Big\{\derivative{\zzz} (C-\ablue{_i^k}): {(k,i)\in{C_5}}\Big\}\\
\cup\Big\{\derivative{\zzz} \mupper{_i^k}: {(k,i)\in{C_{123}}}\Big\}
\cup\Big\{\derivative{\zzz} \theta{^k_i}(\zzz): {(k,i)\in{I}} \Big\}
\cup\Big\{\derivative{\zzz} \varphi^k(\zzz): {k=1,...,K}\Big\}.
\end{multline}
\end{subequations}
Suppose there exist
multipliers $\lambda^1,\lambda^2,\lambda^3,\lambda^4,\lambda^5, \lambda^6$
corresponding to the vectors in set~$A$%
, with $\lambda^5,\lambda^6\geq0$, $\lambda^1_{ki}\lambda^2_{ki}\geq0$ for $(k,i)\in C_2$ and $\lambda^3_{ki}\lambda^4_{ki}\geq0$ for $(k,i)\in C_4$,
and free multipliers $\mu^1,\mu^2,\mu^3,\mu^4,\mu^5,\mu^6$
corresponding to the vectors in set~$B$ such that together they
give rise to a linear combination whose sum is zero.  More precisely, 
\begin{multline}
\label{eq:positive_linear_dependance_system_simple}
\sum_{{(k,i)\in{C_2}}} \left(
    \lambda_{ki}^{1} \derivative{\zzz} \ablue{_i^k}+
    \lambda_{ki}^{2} \derivative{\zzz} \mlower{_i^k}
\right)
+ \sum_{{(k,i)\in{C_4}}} \left(
    \lambda_{ki}^{3} \derivative{\zzz} (C-\ablue{_i^k})+
    \lambda_{ki}^{4} \derivative{\zzz} \mupper{_i^k}
\right)
\\
 + \sum_{{(k, i)\in \bar{I}^\zeta}}
    \lambda_{ki}^{5}    \derivative{\zzz} \zeta_i^k
+ \sum_{{(k, i)\in \bar{I}^Z}}
    \lambda_{ki}^{6}   \derivative{\zzz} Z_i^k(\zzz)
+ \sum_{(k,i)\in{C_1}}      \mu_{ki}^{1} \derivative{\zzz} \ablue{_i^k}
+ \sum_{(k,i)\in{C_{345}}}  \mu_{ki}^{2} \derivative{\zzz} \mlower{_i^k}
\\
+ \sum_{(k,i)\in{C_5}}     \mu_{ki}^{3}  \derivative{\zzz} (C-\ablue{_i^k})
+ \sum_{(k,i)\in{C_{123}}} \mu_{ki}^{4}  \derivative{\zzz} \mupper{_i^k}
+ \sum_{(k,i)\in I}
\mu_{ki}^{6}  \derivative{\zzz} \theta{^k_i}(\zzz)
+ \sum_{k=1}^K
\mu_{k}^{5}  \derivative{\zzz} \varphi^k(\zzz)
= \veczero{}.   
\end{multline}
% ...
where we have written $\derivative{\zzz}\bullet$ to denote the derivative with respect to the full decision variable $\zzz \coloneqq [C, \gamma, \zeta, \ablue{}, \mlower{}, \mupper{}, \uu{}]^\top$. Moving forward, we will ignore the derivatives with respect to $C$ and $\gamma$ and write the remaining components of the derivative as vectors of five blocks corresponding to the derivatives with respect to $[\zeta, \ablue{}, \mlower{}, \mupper{}, \uu{}]$.
For example $\nabla_{[\zeta, \ablue{}, \mlower{}, \mupper{}, \uu{}]}(\ablue{_i^k}) = [0,e_i^k,0,0,0]^\top$
where $e_i^k$ is the standard basis vector of zeros with a $1$ at index $(k,i)$. Thus we take \eqref{eq:positive_linear_dependance_system_simple} and rewrite explicitly
{\scriptsize
\begin{multline}
\label{eq:positive_linear_dependance_system}
\sum_{{(k,i)\in{C_2}}} \left(
    \lambda_{ki}^{1} \begin{bmatrix} 0 \\ e_i^k \\ 0     \\ 0 \\ 0 \end{bmatrix}+
    \lambda_{ki}^{2} \begin{bmatrix} 0 \\ 0     \\ e_i^k \\ 0 \\ 0 \end{bmatrix}
\right)
+ \sum_{{(k,i)\in{C_4}}} \left(
    \lambda_{ki}^{3}    \begin{bmatrix} 0 \\ -e_i^k \\ 0 \\ 0      \\ 0 \end{bmatrix}+
    \lambda_{ki}^{4}    \begin{bmatrix} 0 \\ 0      \\ 0 \\ e_i^k  \\ 0 \end{bmatrix}
\right)
\\
\begin{aligned}
 + \sum_{{(k, i)\in \bar{I}^\zeta}}
    \lambda_{ki}^{5}    \begin{bmatrix} e_i^k \\ 0 \\ 0 \\ 0 \\ 0 \end{bmatrix}
+ \sum_{{(k, i)\in \bar{I}^Z}}
    \lambda_{ki}^{6}    \begin{bmatrix} e_i^k \\ \nabla_\alpha Z_i^k(\zzz) \\ 0 \\ 0 \\ e^k\bar{y}_i^k  \end{bmatrix}
+ \sum_{(k,i)\in{C_1}}      \mu_{ki}^{1} \begin{bmatrix} 0 \\ e_i^k \\ 0     \\ 0 \\ 0 \end{bmatrix}
+ \sum_{(k,i)\in{C_{345}}}  \mu_{ki}^{2} \begin{bmatrix} 0 \\ 0     \\ e_i^k \\ 0 \\ 0 \end{bmatrix}&
\\
+ \sum_{(k,i)\in{C_5}}      \mu_{ki}^{3}    \begin{bmatrix} 0 \\ -e_i^k \\ 0 \\ 0 \\ 0 \end{bmatrix}
+ \sum_{(k,i)\in{C_{123}}}  \mu_{ki}^{4}    \begin{bmatrix} 0 \\ 0 \\ 0 \\ e_i^k \\ 0  \end{bmatrix}
+ \sum_{k=1}^K
\mu_{k}^{5}    \begin{bmatrix} 0 \\ Y^k \\ 0 \\ 0 \\ 0 \end{bmatrix}
+ \sum_{(k,i)\in I}  \mu_{ki}^{6}    \begin{bmatrix} 0 \\ Q(\gamma)_i^k \\ -e_i^k \\ e_i^k \\
e^k y_i^k  \end{bmatrix} &
= \veczero{}.   
\end{aligned}%
\end{multline}%
}%
Inspecting the first row of the vector derivatives, which represents the derivative with respect to $\zeta$, we have that $\lambda^5+\lambda^6=\veczero{}$. By their non-negativity we conclude $\lambda^5=\veczero{}$ and $\lambda^6=\veczero{}$.

We continue by inspecting subvectors of the multipliers based on $C_{12345}$. For a vector $\mu\in\R^{|I|}$, denote its subvector according to the index set $C_l\subseteq{}I$ as
\begin{equation*}
    \mu_{C_l} \coloneqq \mu_{[(k,i)\in C_l]}\in\R^{|C_l|}
    \quad\text{for } l=1, \ldots, 5.
\end{equation*}
Note that the matrix $\Q=\Q(\gamma)\in\R^{|I|\times |I|}$ is positive definite by Lemma~\ref{lemma:psd_kernel}%
, which applies since $\gamma>0$ and the training inputs are distinct. It is still a function of $\gamma$, but we drop the argument for brevity.

Consider the derivative with respect to
$\mlower{_i^k}$ for $(k,i)\in C_2$; i.e., the third row in our block notation in \eqref{eq:positive_linear_dependance_system}. We have that 
$\lambda^{2}-\mu^{6}_{C_2}=\veczero{}$,
which implies that
$\mu^{6}_{C_2}$
is equal to $\lambda^2$. For $(k,i)\in C_1$, the same row block gives $\mu^{6}_{C_1}=\veczero{}$, since $\mlower{_i^k}>0$ and the gradient of $\mlower{_i^k}$ does not appear in~\eqref{eq:positive_linear_dependance_system_simple}.

In a similar manner, we now consider the derivative with respect to
$\mupper{_i^k}$ for $(k,i)\in C_4$; i.e., the fourth row block in \eqref{eq:positive_linear_dependance_system}. We have that
$\lambda^{4}+\mu^{6}_{C_4}=\veczero{}$,
which implies that
$\mu^{6}_{C_4}$
is equal to $-\lambda^4$. For $(k,i)\in C_5$, the same row block gives $\mu^{6}_{C_5}=\veczero{}$. Moreover, since $\lambda^6=\veczero{}$, the fifth row block, i.e. the derivative with respect to $\uu{}$, gives $\Y^\top\mu^6=\veczero{}$.
Consider now the derivative with respect to $\ablue{_i^k}$ for all $(k,i)\in I$, i.e. the second row block in \eqref{eq:positive_linear_dependance_system}. Since $\lambda^6=\veczero{}$, it reads
\begin{equation*}
\lambda^{1}-\lambda^{3}+\mu^{1}-\mu^{3}+\Q\mu^{6}+\Y\mu^{5}=\veczero{},
\end{equation*}
where $\lambda^{1}$, $\lambda^{3}$, $\mu^{1}$ and $\mu^{3}$ are extended by zeros from $C_2$, $C_4$, $C_1$ and $C_5$, respectively, to $I$. Taking the inner product with $\mu^{6}$ and using $\mu^{6}_{C_1}=\veczero{}$, $\mu^{6}_{C_5}=\veczero{}$ and $\Y^\top\mu^6=\veczero{}$ gives
\begin{equation*}
    \left(\mu^{6}_{C_2}\right)^\top\lambda^{1}
    -\left(\mu^{6}_{C_4}\right)^\top\lambda^{3}
    +\left(\mu^{6}\right)^\top\Q\,\mu^{6}
    =\left(\lambda^{2}\right)^\top\lambda^{1}
    +\left(\lambda^{4}\right)^\top\lambda^{3}
    +\left(\mu^{6}\right)^\top\Q\,\mu^{6}
    =0.
\end{equation*}
Since $\lambda^1_{ki}\lambda^2_{ki}\geq0$ for $(k,i)\in C_2$ and $\lambda^3_{ki}\lambda^4_{ki}\geq0$ for $(k,i)\in C_4$, all three terms are non-negative and hence zero. By the positive definiteness of $\Q$ (Lemma~\ref{lemma:psd_kernel}), we must have $\mu^{6}=\veczero{}$ and therefore $\lambda^{2}=\mu^{6}_{C_2}=\veczero{}$ and $\lambda^{4}=-\mu^{6}_{C_4}=\veczero{}$.

At this point, we have shown that $\lambda^{2}$, $\lambda^{4}$, $\lambda^{5}$, $\lambda^{6}$ and
$\mu^{6}$ are all zero. Now, consider the remaining vectors involved in~\eqref{eq:positive_linear_dependance_system} whose multipliers we have not yet concluded are zero.
The third and fourth row blocks of~\eqref{eq:positive_linear_dependance_system} give $\mu^{2}=\mu^{6}_{C_{345}}=\veczero{}$ and $\mu^{4}=-\mu^{6}_{C_{123}}=\veczero{}$, and the second row block reduces to $\lambda^{1}-\lambda^{3}+\mu^{1}-\mu^{3}+\Y\mu^{5}=\veczero{}$. Its components for $(k,i)\in C_3$ read $y_i^k\mu^{5}_{k}=0$. For each fold $k$, assumption~\eqref{eq:assumption_on_alpha} provides an index $i$ with $0<\ablue{_i^k}<C$, for which the complementarity constraints force $\mlower{_i^k}=\mupper{_i^k}=0$, so that $(k,i)\in C_3$. Since $y_i^k\in\{-1,1\}$, we obtain $\mu^{5}=\veczero{}$. The remaining components for $C_1$, $C_2$, $C_4$ and $C_5$ then give $\mu^{1}=\veczero{}$, $\lambda^{1}=\veczero{}$, $\lambda^{3}=\veczero{}$ and $\mu^{3}=\veczero{}$, respectively.
Therefore, equation~\eqref{eq:positive_linear_dependance_system_simple} is satisfied for
multipliers with the above sign restrictions only when $\lambda=\veczero{}, \mu=\veczero{}$, giving us that
\MPCCMFCQC{} holds by Lemma~\ref{lemma:positive-linearly-independent}(iv), and hence \MPCCMFCQR{} holds by Theorem~\ref{theorem:cq_implications}.
\end{proof}

\begin{remark}
The above proof uses the dual form of
\MPCCMFCQC{}.
It would be equivalent to show for the sets defined in~\link{eq:set_A}{(\ref*{eq:set_A},b)} that $B$ is linearly independent and%
, for every subset $J$ of the biactive indices, there exists a direction $\ddd$ such that $B\ddd=\veczero{}$ and
$A_J\ddd>\veczero{}$, where $A_J$ is obtained from $A$ by changing the sign of the gradients of the biactive pairs in $J$ (cf. Definition~\ref{def:MPCC-MFCQ-C})%
, and where $A$, $A_J$ and $B$ also denote the matrices whose rows are the corresponding gradient vectors.
For $J=\emptyset$, i.e. for the \MPCCMFCQR{}, this direction $\ddd$ is illustrated geometrically by the orange arrows in Figure~\ref{fig:five_cases}. Also, as a key consequence of Proposition~\ref{proposition:satisfy_mpccmfcqr}, the
constraint qualifications assumed in Theorem~\ref{theorem:penalisation_convergence} and Theorem~\ref{theorem:relaxation_convergence} hold accordingly. 
\end{remark}

We now show that problem \MODEL{} in general fails the stronger constraint qualification \MPCCMFCQT{}, and therefore displaying that the second implication in Theorem~\ref{theorem:cq_implications} is strict for the application problem under consideration in this section.

\begin{proposition}
\label{proposition:fail_mpcc_licq}
The \MPCCMFCQT{} fails at all feasible points $\zzz  \coloneqq  (C, \gamma, \zeta, \ablue{}, \mlower{}, \mupper{}, \uu{})$ of the program \MODEL{} with more than two biactive indices; i.e., with  $|I^{GH}(\zzz)|>2$. 
\end{proposition}
\begin{proof}
\MPCCMFCQT{} requires that the following family of vectors be linearly independent:
\begin{equation}
\begin{aligned}
	\label{eq:set_of_LI_vectors}
    &\Big\{\nabla h_i(\zzz):i\in I^h\Big\} \cup
    \Big\{\nabla G_i(\zzz):i\in {I^{G}(\zzz)\cup I^{GH}(\zzz)}\Big\} \cup
    \Big\{\nabla H_i(\zzz):i\in{I^{H}(\zzz) \cup I^{GH}(\zzz)}\Big\}.
\end{aligned}
\end{equation}
The constraint functions $h$, $G$ and $H$ are in the variable $(\ablue{}, \mlower{}, \mupper{}, \uu{}, C, \gamma)$, which has a dimension of $3(K-1)n+K+2$, where $n$ is the
total number of data points and $K$ is the number of folds in the cross-validation.  Observe that all the vectors in $\eqref{eq:set_of_LI_vectors}$ are constant with respect to $\zeta$, so their derivatives are zero with respect to $\zeta$. There are a total of $3(K-1)n+K+|I^{GH}|$ vectors in this set (\ref{eq:set_of_LI_vectors}). Therefore, if $|I^{GH}|>2$, then there are more vectors than non-zero dimensions, and so they cannot be linearly independent.
\end{proof}

By Theorem~\ref{theorem:cq_implications}, we conclude that \MPCCLICQ{} also fails for program~\MODEL{}. 
These results are conditional on the hyperparameters $C$ and $\gamma$ being strictly positive. This is a very sensible restriction, since on any dataset for which $C=0$ or $\gamma=0$ is optimal, the best classifier ignores the data entirely.
Proposition~\ref{proposition:satisfy_mpccmfcqr} therefore makes no claim at such boundary points $C=0$, $\gamma=0$, which are admitted by~\MODEL{} but are of no practical interest.

We close this subsection by relating these results to earlier work. 
Firstly, the CQs of an MPCC reformulation of the \emph{linear} SVM hyperparameter tuning problem (i.e., when $Q$ is a linear kernel) are analysed by Li et al.~\cite{Li2022linear}, in 2022, with the training problem in primal form.
Secondly, a model very similar to~\MODEL{} has already been introduced by Coniglio et al.~\cite{Coniglio2022}.
They also use a reformulation for an RBF kernel that is derived from the dual SVM, but with the multipliers $\mlower{_i^k}$ eliminated using $\theta_i^k(\zzz)=0$ (see Section~\ref{section:behaviour}).
Another notable difference is that they keep the max of the hinge loss in the upper-level objective, smoothing it only in their experiments, rather than introducing slack variables.
In Theorem~20 of this preprint~\cite{Coniglio2022}, they assert that the MPEC-MFCQ, i.e., the \MPCCMFCQT{}, holds at every feasible point, but they prove something much weaker.
This stems from a misreading of \emph{positive linear independence}: they rule out only combinations in which all multipliers are non-negative, rather than splitting the gradients into the sets $A$ and $B$ of~\eqref{eq:MFCQ_dual_form_part_2}, where the multipliers relating to $B$ are free.
The proof carries this mistake through to \cite[Equations~(41) and~(44)]{Coniglio2022}, where they conclude that the coefficients of $\nabla G$ and $\nabla H$ vanish only because they are non-negative.
This motivates the care of Subsection~\ref{section:constraint_qualification} in establishing the correct dual form of each qualification.
Hence the \MPCCMFCQC{}, and with it the \MPCCMFCQR{}, is established for the first time in Proposition~\ref{proposition:satisfy_mpccmfcqr}, whereas the \MPCCMFCQT{}, claimed by~\cite{Coniglio2022}, in fact fails by Proposition~\ref{proposition:fail_mpcc_licq}.

%% file: Article/article_50_experiments.tex
\section{Numerical experiments}
\label{section:experiments}
% ...
In this section, we present our numerical experiments on hyperparameter tuning for support vector machines where we solve large instances of the highly non-convex MPCC problem \MODEL{} using the solution methods presented in Section~\ref{section:solution_methods} under various settings.

%% file: Article/article_51_implementation.tex
\subsection{Implementation details}
\label{section:implementation}
% ...
Each program and system of equations we wish to solve, 
\ifsuplementaryincluded
such as~\eqref{eq:SVMdual}, \MODEL{}, \eqref{eq:relaxation}, \eqref{eq:penalisation} and~\eqref{eq:system}, 
\else
such as~\MODEL{}, \eqref{eq:relaxation}, and \eqref{eq:penalisation},
\fi
is written in the algebraic modelling language AMPL~\cite{Fourer1990}.
This allows us to robustly compute the derivatives and interface with the solvers we wish to use.
% ...
The algorithms from Section~\ref{section:solution_methods} are implemented in Python~\cite{Python2009} with use of the library NumPy~\cite{Harris2020} for efficient linear algebra computation.
% ...
For solving the subproblems~\eqref{eq:relaxation} and~\eqref{eq:penalisation}, we trialled two different solvers: Interior Point Optimizer~(IPOPT)~\cite{Wachter2006}, an open source solver that implements a primal-dual interior point method with line search, and Artelys' Knitro~\cite{Byrd2006}, a commercial solver that implements an interior point trust region method. 
Both converged reliably in our context though our results are stated for Knitro as it shows faster runtimes.

Regarding cross-fold validation, it is common in the literature to set $K=3$ (see \cite[p. 91]{Geron2022}) and $K=5$ (see \cite[Section 7.10.1]{Hastie2009}) as the number of folds. 
We set $K=3$ because it gives a more stable validation accuracy for small datasets.
% ...
The stopping tolerance for the subprograms is gradually decreased to $10^{-6}$. The overarching \MPCC{} loop in Algorithms~\ref{alg:penalisation} and~\ref{alg:relxation} is terminated when one of each pair of complementarity constraint functions is at most $10^{-6}$, i.e.,
\begin{equation}
    \label{eq:stopping_criterion}
    \max_{i}\left( \min(|G_i(\zzz)|, \;|H_i(\zzz)|) \right)\leq 10^{-6}.
\end{equation}

For the penalisation method, we select the initial penalty parameter to be $100$ and sequentially increase it by multiples of $10$.  
For the relaxation method, we select the initial relaxation parameter to be $0.01$ and sequentially decrease it by multiples of $10^{-1}$.
Many publications have suggested using the so-called \textit{exact penalisation} (e.g.,~\cite{Ralph2004}) and \textit{exact relaxation} (e.g.,~\cite{Hoheisel2013}) approaches.
Here, the parameter $\penaltyparam{}$ (resp. $\relaxparam{}$) is fixed at some chosen value.
If chosen correctly, only a single instance of the subproblem \penalisation{\penaltyparam{}}  (resp. \relaxation{\relaxparam{}}) is solved resulting in an MPCC feasible solution.  
In our experiments, we implement this approach but note that our designation of fixed parameter comes from an intimate knowledge of the problem and is, in general, hard to find.
Regarding the theory of Section~\ref{section:solution_methods}, the above geometric updates satisfy $\penaltyparam{}^t\rightarrow\infty$ and $\relaxparam{}^t\rightarrow0$, so that the sequential methods meet the parameter requirements of Theorems~\ref{theorem:penalisation_convergence}, \ref{theorem:penalisation_global_iterates} and~\ref{theorem:relaxation_convergence}.
For the exact approaches, Theorem~\ref{theorem:penalisation_global_iterates} does not apply, and the MPCC feasibility of the computed point is not guaranteed.
However, if a KKT point of \penalisation{\penaltyparam{}} with a fixed parameter $\penaltyparam{}$ is MPCC feasible, then it is S-stationary~\cite[Theorem~5.2]{Ralph2004}.
The same holds for a KKT point of \relaxation{\relaxparam{}} with a fixed $\relaxparam{}>0$, since the constraints $G_i(\zzz)H_i(\zzz)\leq\relaxparam{}$ are then inactive, so that their multipliers vanish.


We use 18 datasets, 6 synthetically generated and 12 from real-world observations.
\ifsuplementaryincluded
For full details see the supplementary material Subsection~\ref{SM:datasets}.
\else
For full details of these datasets see~\cite[Subsection~\ref{SM:datasets}]{Ward2025}.  
\fi
When \MODEL{} is parametrised with these data sets, it is large in comparison to other non-convex MPCCs that have been solved in the literature so far. The median problem size is 1773 variables, whereas for example, see Leyffer's famous MacMPEC collection~\cite{MacMPEC}, which has a median problem size of 49 variables.

%% file: Article/article_52_behaviour.tex
\subsection{Numerical behaviour of penalisation and relaxation methods}
\label{section:behaviour}

\begin{figure}
\centering
\includegraphics[width=.4\linewidth]{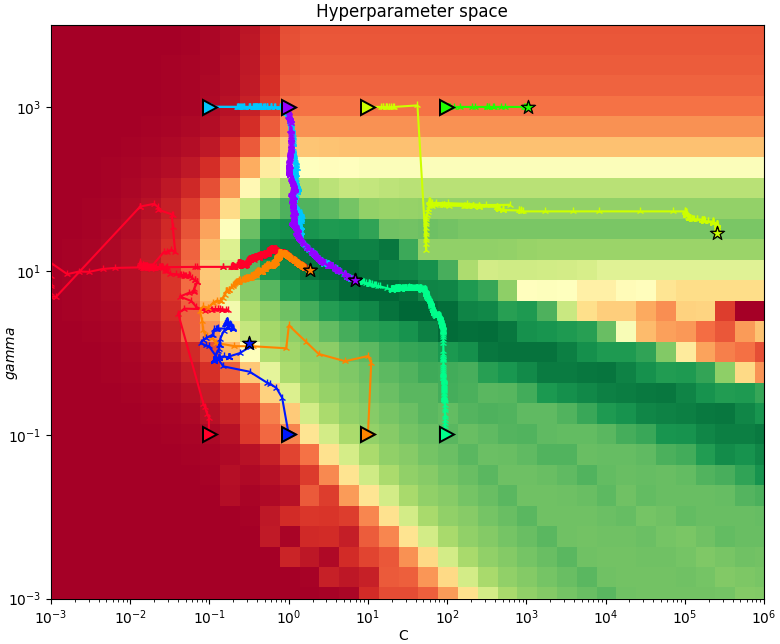}
\caption{This heat map shows the hyperparameter space of the Moons dataset. Green regions correspond to $(C,\gamma)$ combinations, which allow the SVM to achieve high validation accuracy. The eight coloured lines represent the path of iterates produced by the Relaxation method from eight different starting points (triangles) to various stationary points (stars).}
\label{fig:multistart}
\end{figure}

Our first observation is that due to the highly non-convex nature of \MODEL{}, the stationary point to which the algorithms converge is quite variable depending on the initialisation of the decision variables; see Figure~\ref{fig:multistart}. 
Throughout the experiments a multi-start strategy is used which is fully described in 
\ifsuplementaryincluded
supplementary material Subsection~\ref{SM:starting_point}.
\else
supplementary material~\cite[Subsection~\ref{SM:starting_point}]{Ward2025}.
\fi

Secondly, we observe that for the sequential relaxation method, parameter $\relaxparam{}$ almost always decreases to $10^{-6}$.  
That is, at each iteration $t$ with $\relaxparam{}^t = 10^{-2}, 10^{-3}, 10^{-4}, 10^{-5}$, the solutions found by Algorithm~\ref{alg:relxation} for the subproblem \relaxation{\relaxparam{}^t} are not MPCC feasible (i.e. \eqref{eq:stopping_criterion} does not hold) and the algorithm only terminates when the relaxation parameter reaches $10^{-6}$ and the complementarity constraints are forced within the exit tolerance.
On the other hand, the penalisation method terminates with a range of different penalty parameters $\penaltyparam{}^t$ between $10^3$ and $10^6$ producing MPCC feasible solutions.
For more details see
\ifsuplementaryincluded
Table~\ref{table:pi_and_feasibility}.
\else
\cite[Section~\ref{SM:futher_numerics}]{Ward2025}.  
\fi
Further research should investigate the structure of \MODEL{} to quantify why the penalisation approach is able to find an MPCC feasible solution before the penalty parameter~$\penaltyparam{}^t$ potentially gets too large and causes \penalisation{\penaltyparam{}^t} to become ill-conditioned (e.g.,~\cite[Example 3]{Leyffer2006interior}), as is usually observed for penalty-based algorithms~\cite{Murray1971,Wright1994}.

Next, we note that the complementarity constraint functions of \MODEL{} are all linear.  
Also, observe that for all indices $k=1,\dots,K, i=1,\dots,n^k$, the equality constraint~\eqref{eq:theta} could be rearranged and used to substitute out either $\mlower{_i^k}$ or $\mupper{_i^k}$ from the complementarity  constraints as is done by Coniglio et al.~\cite{Coniglio2022}.  
This would reduce the number of variables and constraints by $(K-1)n$ but introduce a nonlinear term in the complementarity  constraints. 
In line with the advice of Fletcher and Leyffer~\cite{fletcher2002numerical,Leyffer2006complementarity}, we find that keeping the complementarity constraints linear is more beneficial than reducing the number of variables via substitution.

Finally, we remark that the MPCC optimisation approach for hyperparameter optimisation, studied in this paper, sometimes overfits to the validation data. 
To monitor this, on top of the training and validation data that parametrises problem~\MODEL{}, we require further data that is entirely unseen by the model; this will be called~\emph{test} data and is formed by putting aside 10\% of examples and parametrising with the remaining 90\%.
\ifsuplementaryincluded
Refer to Figure~\ref{fig:k_fold_cross_val}.
\fi
The overfitting can then be observed when the classifiers, represented by the solutions of problem~\MODEL{}, achieve higher validation accuracy than test accuracy. 
Evidence of this is presented in our later experiments, in particular, see the Ionosphere and Moons datasets in 
\ifsuplementaryincluded
Table~\ref{table:experiment_results}.
\else
\cite[Table~\ref{table:experiment_results}]{Ward2025}.
\fi

%% file: Article/article_53_performance.tex
\subsection{Performance evaluation in comparison to derivative-free optimisation}
\label{section:performance}
% ...
Across industry and academia, Derivative-Free Optimisation (DFO) remains by far the most used tool for hyperparameter tuning~\cite{Marsland2011}.  
It is simple and effective.  
The defining characteristic of DFO is that it treats training as a black box, meaning it does not use any information from the training process other than the final evaluation of the model in optimising the hyperparameters.  
We compare the MPCC optimisation to four such methods: \textit{grid search}, \textit{random search}, \textit{Bayesian optimisation}, and \textit{pattern search}. 
For the exact derivative-free parametrisations and procedures used see
\ifsuplementaryincluded
Subsection~\ref{SM:derivative_free}.
\else
supplementary material~\cite[Subsection~\ref{SM:derivative_free}]{Ward2025}. 
\fi

\begin{figure}
    \centering
    \includegraphics[width=1\linewidth]{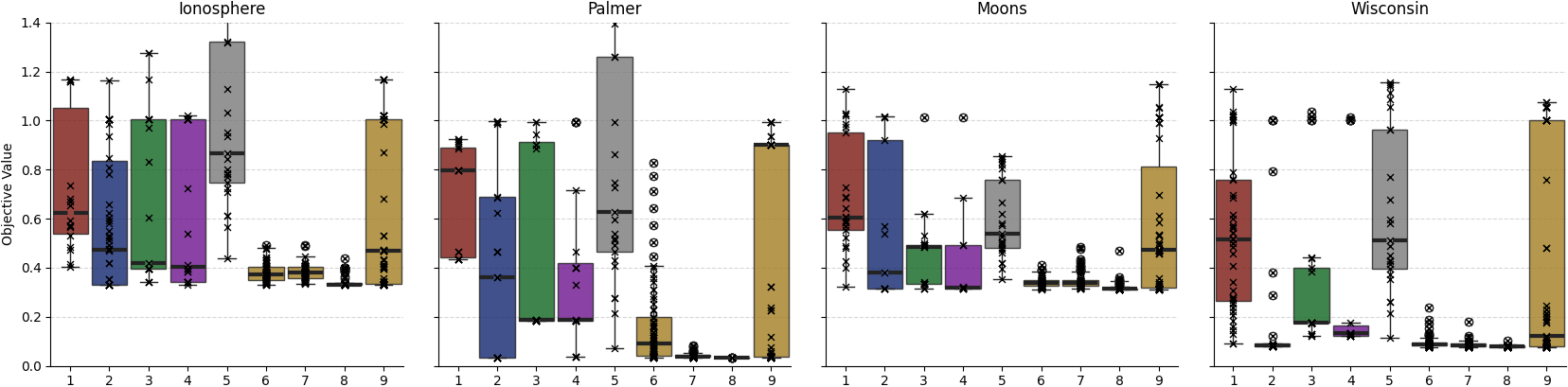}
    \caption{These box plots show the distribution of objective values achieved by initialising from a range of different starting points.  The four axes represent the datasets (left to right):  Ionosphere, Palmer, Moons and Wisconsin.  Within each axis, the box plots represent the solution method (left to right):
    \textbf{\color{c1}1.~Penalisation (Exact),}
    \textbf{\color{c2}2.~Penalisation (Sequential),}
    \textbf{\color{c3}3.~Relaxation (Exact),}
    \textbf{\color{c4}4.~Relaxation (Sequential),}
    \textbf{\color{c5}5.~Semi-Smooth Newton,}
    \textbf{\color{c6}6.~Grid Search,}
    \textbf{\color{c6}7.~Random Search,}
    \textbf{\color{c6}8.~Bayesian Optimisation,}
    \textbf{\color{c6}9.~Pattern Search.}
    }
    \label{fig:box_plots}
\end{figure}

Note that
\ifsuplementaryincluded
\eqref{eq:SVMdual}
\else
the dual SVM as is solved by most libraries~\cite[Equation (1)]{Fan2005}
\fi
has%
, for each fold $k$, only $n^k<n$ variables and $2n^k+1$ constraints.
On the other hand,
with $K=3$ folds, the MPCC model~\MODEL{} has $7n+5$ variables,
$2n+2$ inequality constraints, $2n+3$ equality constraints and $4n$ complementarity constraints.
Therefore, the comparison must be made as to whether better performance can be achieved by minimising the much smaller problem 
\ifsuplementaryincluded
\eqref{eq:SVMdual}
\else
dual SVM problem~\cite[Equation (1)]{Fan2005}
\fi
many times employing a DFO method or by minimising the much larger \MODEL{} with the solution methods from Section~\ref{section:solution_methods}.

We also add the semi-smooth Newton method to the experiment suite.  
This is an entirely different approach to solving MPCCs and converges to an M-stationary point, making it an interesting comparison, given that the sequential penalisation and Scholtes relaxation algorithms that we primarily study in this paper have been shown to instead converge to
C-stationary points, and to S-stationary points under an additional assumption.
For the full algorithm see 
\ifsuplementaryincluded
supplementary material Subsection~\ref{SM:m_stationarity_system}.
\else
supplementary material~\cite[Subsection~\ref{SM:m_stationarity_system}]{Ward2025}.  
\fi

\begin{figure}
    \centering
    \includegraphics[width=\linewidth]{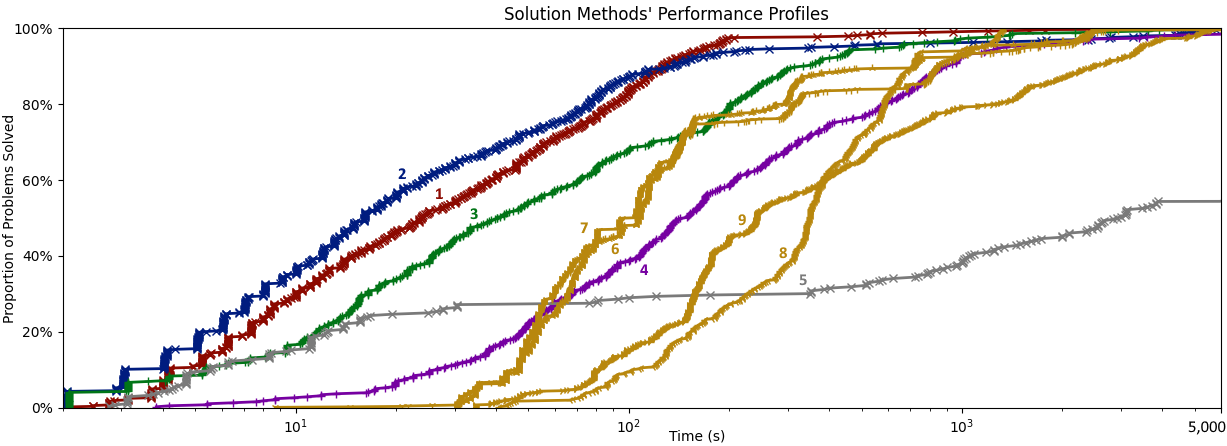}
    \caption{This performance profile graph plots the proportion of instances that were solved to the target test accuracy (vertical axis) within the given time (horizontal axis).
    The lines are coloured and enumerated for clarity as follows:
    \textbf{\color{c1}1.~Penalisation (Exact),}
    \textbf{\color{c2}2.~Penalisation (Sequential),}
    \textbf{\color{c3}3.~Relaxation (Exact),}
    \textbf{\color{c4}4.~Relaxation (Sequential),}
    \textbf{\color{c5}5.~Semi-Smooth Newton,}
    \textbf{\color{c6}6.~Grid Search,}
    \textbf{\color{c6}7.~Random Search,}
    \textbf{\color{c6}8.~Bayesian Optimisation,}
    \textbf{\color{c6}9.~Pattern Search.}
    }
    \label{fig:perfomance_profile}
\end{figure}

For each of the eighteen datasets
\ifsuplementaryincluded
detailed in supplementary material Subsection~\ref{SM:datasets}, 
\fi
each solution method is run 100 times from different initialisations.  
We document the range of objective values achieved; 
\ifsuplementaryincluded
see Table~\ref{table:experiment_results} and box-plots in Figure~\ref{fig:box_plots}.
\else
see Figure~\ref{fig:box_plots}.
\fi
In general, the MPCC methods run much faster, but, on the other hand, converge to a large variety of stationary points, many of which have poor objective function values.
To gain an understanding of their overall performance, we attribute a target test accuracy to each of the datasets.  
We then study the average runtime to achieve the target test accuracy.  
These results are plotted as a performance profile in Figure~\ref{fig:perfomance_profile}. 

Overall, the obtained results make a convincing argument that our model~\MODEL{}, being solved with the sequential penalisation method, can consistently outperform derivative-free methods in the tuning of SVM hyperparameters. 
Additionally, the sequential approach, when compared to the exact approach, often achieves a smaller range of objective function values (see Figure~\ref{fig:box_plots}).  
The exact relaxation outperforms the derivative-free methods on a majority of instances.  
The semi-smooth Newton method performs well on small problems but fails to converge for larger problems -- it seems to get stuck at the non-differentiable corners of the stationarity function.
Among the derivative-free methods, grid and random search perform equally well on all but the very largest problems for which Bayesian optimisation is preferable.
The penalisation methods clearly come out on top overall.
Sequential penalisation is the fastest method to achieve the target test accuracy in 88\% of instances.
For the other 12\%, exact penalisation is the fastest.
However, it must be noted that, in practice, a user may prefer to always use sequential penalisation as it does not require the selection of a penalty parameter a priori.

%% file: Article/article_60_conclusion.tex
\section{Conclusion}
\label{section:conclusion}
% ...
We began by highlighting the distinction between the definitions of \MPCCMFCQT{} and \MPCCMFCQR{}, which is often overlooked in the literature.
We also introduced the \MPCCMFCQC{}, which lies strictly between them and is associated with the sign pattern of C-stationarity.
From there, we presented two solution methods for the \MPCC{}.
For
both the penalisation and relaxation methods, we provided new convergence proofs that an inexact sequence of KKT points converges
to a C-stationary point under the assumption of \MPCCMFCQC{}, and to an S-stationary point under the assumption of \MPCCMFCQR{} together with an additional assumption at the biactive indices.
This provides a generalisation over the existing literature and allows the methods to be applied to a broader range of problems.

In the paper's second half, we presented our model \MODEL{} for the hyperparameter tuning of support vector machines.  
We then proved that this satisfies the
\MPCCMFCQC{}, and hence the \MPCCMFCQR{}%
, under mild assumptions, thus satisfying the
constraint qualifications assumed in the theorems proved earlier in the paper.
Finally, we ran an extensive set of experiments on real-world data sets.
In this context, we showed the practical use of our methods in beating the derivative-free search that is most commonly used. 
We conclude that the sequential penalisation method both converges faster and is more likely to converge to a higher quality objective value than the other methods.

For future work, we will develop first-order gradient descent methods to solve \MODEL{} such as those being studied by Bengio~\cite{Bengio2000} and Chapelle et al.~\cite{Chapelle2001}.  
These have a distinct advantage in machine learning applications as they can be efficiently parallelised on modern graphics processing units.

%% file: Article/article_declarations.tex

\backmatter

\ifsuplementaryincluded
\nocite{Hunter2007,Waskom2021}
\else
\input{Article/article_acknowledgements}
\fi

\ifsuplementaryincluded\else
\subsection*{Statements and Declarations}
The work of the first author is jointly funded by \href{https://www.das-ltd.co.uk/}{DAS Ltd} and the School of Mathematical Sciences at the University of Southampton through a PhD studentship grant.
The second author's research is partly funded by an EPSRC grant with reference EP/X040909/1.
The authors have no competing interests to declare that are relevant to the content of this article.
All datasets analysed during this study are publicly available; the individual sources are cited in
\ifsuplementaryincluded
supplementary material Subsection~\ref{SM:datasets}.
\else
the supplementary material~\cite[Subsection~\ref{SM:datasets}]{Ward2025}.
\fi
The AMPL models and Python implementations of the algorithms reported in this paper are available from the corresponding author on reasonable request.
All authors contributed to theory established in this paper.
All authors read and approved the final manuscript.
This article is accompanied by supplementary material containing further technical details, proofs and numerical results; it is available at~\cite{Ward2025}.
\fi

%% file: Article/article_acknowledgements.tex

\subsection*{Acknowledgements}
All graphics are produced using Matplotlib Pyplot~\cite{Hunter2007} (version 3.8.2) and seaborn~\cite{Waskom2021} (version 0.13.2).
Thank you to Qingna Li from the Beijing Institute of Technology for her correction regarding assumption~\eqref{eq:assumption_on_alpha}.
We are grateful to an anonymous reviewer for the counterexample in
\ifsuplementaryincluded
supplementary material Subsection~\ref{SM:counterexample_relaxation},
\else
the supplementary materials~\mbox{\cite[Subsection~\ref{SM:counterexample_relaxation}]{Ward2025}},
\fi
which revealed an error in an earlier version of Theorem~\ref{theorem:relaxation_convergence}.

%% file: Supplementary/supplementary__sections.tex
\input{Supplementary/supplementary_10_notation}
\newpage
\input{Supplementary/supplementary_20_relationship_to_bilevel}
\input{Supplementary/supplementary_21_model_as_bilevel}
\input{Supplementary/supplementary_22_positive_definite}
\input{Supplementary/supplementary_30_futher_solution_methods}
\input{Supplementary/Supplementary_31_proof_global}
\input{Supplementary/supplementary_33_proof_relaxation}

\input{Supplementary/supplementary_33b_counterexample_relaxation}
\input{Supplementary/supplementary_34_m_stationarity}

\input{Supplementary/supplementary_35_lagrangian}
\input{Supplementary/supplementary_40_svm}
\input{Supplementary/supplementary_41_primal_svm}

\input{Supplementary/supplementary_42_dual_svm}
\input{Supplementary/supplementary_43_kkt_svm}
\input{Supplementary/supplementary_44_reconstruction}
\input{Supplementary/supplementary_50_futher_numerics}
\input{Supplementary/supplementary_51_starting_point}

\input{Supplementary/supplementary_52_datasets}
\input{Supplementary/supplementary_53_derivative_free}
\input{Supplementary/supplementary_54_results}

\begin{figure}[H]
    \centering
    \includegraphics[width=\textwidth]{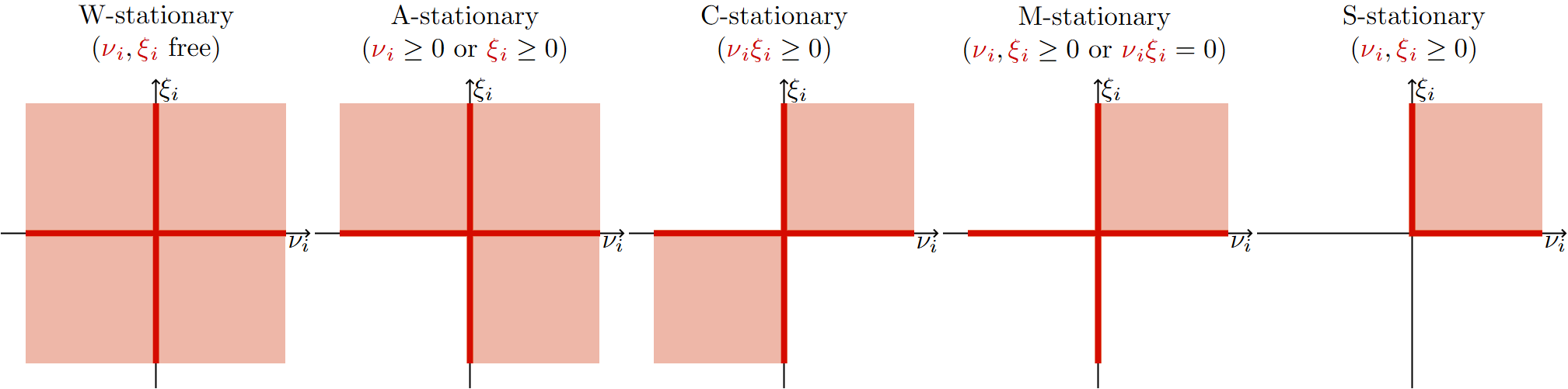}
    \caption{This diagram highlights the differences in the W, A, C, M and S stationary conditions. Areas coloured in red show possible values the multipliers $\multG{_i}$ and $\multH{_i}$ may take for a stationary point $\bar{\zzz}$ where $G_i(\bar{\zzz}) = H_i(\bar{\zzz})=0$}
    \label{fig:Stationary}
\end{figure}

%% file: Supplementary/supplementary_10_notation.tex
\section{Table of notation}%
\label{SM:notation}%
% ...
\begin{center}
\begin{tabular}{l c p{10cm} }
\multicolumn{3}{c}{\underline{Parameters}}\\
$\gamma, \delta$   & $\coloneqq$ & Hyperparameters that determine the kernel.\\
$C$                & $\coloneqq$ & Regularisation hyperparameter for \eqref{eq:soft_margin}.\\
$\penaltyparam{}^t$& $\coloneqq$ & Penalty parameter at iteration $t$.\\
$\relaxparam{}^t$  & $\coloneqq$ & Relaxation parameter at iteration $t$.\\
\multicolumn{3}{c}{}\\
\multicolumn{3}{c}{\underline{Training and Validation Data}}\\
$x_i$            & $\coloneqq$ & Data point $i$.\\
$x_i^k$          & $\coloneqq$ & Training data point $i$ in cross fold experiment $k$.\\
$\bar{x}_i^k$    & $\coloneqq$ & Validation data point $i$ in cross fold experiment $k$.\\
\multicolumn{3}{c}{}\\
$y_i$               & $\coloneqq$ & Label for data point $i$.\\
$y_i^k$          & $\coloneqq$ & Label for training data point $i$ in cross fold experiment $k$.\\
$\bar{y}_i^k$    & $\coloneqq$ & Label for validation data point $i$ in cross fold experiment $k$.\\
\multicolumn{3}{c}{}\\
$n$            & $\coloneqq$ & Total number of data points.\\
$n^k$          & $\coloneqq$ & Number of training data points in cross fold experiment $k$.\\
$\bar{n}^k$    & $\coloneqq$ & Number of validation data points in cross fold experiment $k$.\\
\multicolumn{3}{c}{}\\
$Q(\gamma)_{ij}$        & $\coloneqq$ & $y_i y_j \exp\left(-\gamma\|x_i-x_j\|^2\right)$ the full RBF kernel matrix.\\
$Q(\gamma)_{ij}^k$      & $\coloneqq$ & $y^k_i y^k_j \exp\left(-\gamma\|x^k_i-x^k_j\|^2\right)$ the $k$\textsuperscript{th} training RBF kernel matrix.\\
$\bar{Q}(\gamma)_{ij}^k$& $\coloneqq$ & $\bar{y}_i^k y_j^k \exp\left(-\gamma\|\bar{x}^k_i-x^k_j\|^2\right)$ the $k$\textsuperscript{th} validation RBF kernel matrix.\\
\multicolumn{3}{c}{}\\
\multicolumn{3}{c}{\underline{Decision Variables}}\\
$w, b$                          & $\coloneqq$ & Weights and bias for the primal program \eqref{eq:soft_margin}.\\
$\ablue{}, \etablue{}$        & $\coloneqq$ & Dual variables for the dual program \eqref{eq:SVMdual}.\\
$\mlower{}, \mupper{}, \uu{}$       & $\coloneqq$ & Multipliers in the KKT conditions of the dual SVM \eqref{eq:lowerLevelKKT}.\\
$\multg{}, \multh{}, \multG{}, \multH{}{}$       & $\coloneqq$ & Lagrange multipliers for the MPCC; see \eqref{eq:MPCC_Lagragian}.\\
\multicolumn{3}{c}{}\\
\multicolumn{3}{c}{\underline{Functions}}\\
$f  $             & $\coloneqq$ & Objective function.\\
$g  $             & $\coloneqq$ & Inequality constraint function i.e. $g(\zzz)\geq 0$.\\
$h  $             & $\coloneqq$ & Equality constraint function i.e. $h(\zzz)= 0$.\\
$G, H$            & $\coloneqq$ & Complementarity constraint functions i.e. $0\leq G(\zzz)\bot H(\zzz) \geq 0$.\\
$Z, \theta, \varphi$ & $\coloneqq$ & Specific constraint functions of model \MODEL{} defined in \eqref{eq:constraint_function_definitions}.\\
$F$               & $\coloneqq$ & System of nonlinear equations corresponding to M-stationary points \eqref{eq:system}.\\
\multicolumn{3}{c}{}\\
\multicolumn{3}{c}{\underline{Index Sets}}\\
$I^g, I^h  $      & $\coloneqq$ & Index sets of the active constraints; see Section~\ref{section:constraint_qualification}.\\
$I^G,I^H,I^{GH}$  & $\coloneqq$ & Index sets of the complementarity constraints; see equation~\eqref{eq:index_sets}.\\
$I, \bar{I}$      & $\coloneqq$ & Sets of indices $(k,i)$ for fold k, example i in training and validation respectively; see~equation~\eqref{eq:flat_index}.\\
$C_{12345}$       & $\coloneqq$ & The five index cases used for the proof of Proposition~\ref{proposition:satisfy_mpccmfcqr}.\\
\end{tabular}
\end{center}

%% file: Supplementary/supplementary_20_relationship_to_bilevel.tex
\section{Relationship to bilevel optimisation}
\label{SM:relationship_to_bilevel}

As our support vector machine (SVM) application can be viewed as a bilevel optimisation problem, it is important to introduce the definition of a bilevel optimisation program and present its relationship to the MPCC. 

First, recall that bilevel optimisation dates back to Stackelberg's economic leader-follower game theory \cite{Stackelberg}, where the leader first gets to select their decision variable,  that we denote here by $\zone\in\R^n$, to optimise their upper-level objective function $\fupper(\zzz_1,\zzz_2)$ subject to some upper-level constraints $\gupper(\zzz_1)\geq 0$ and $\hupper(\zzz_1)=0$. 
Whatever selection the leader makes for $\zzz_1$, the follower will then get to choose their own variable, represented here by 
$\ztwo\in\R^m$, 
to optimise their lower-level objective function $\flower(\zzz_1,\zzz_2)$ subject to some lower-level constraints $\glower(\zzz_1,\zzz_2)\geq 0$ and $\hlower(\zzz_1,\zzz_2)=0$. However, as the leader's choices depend on the follower's variable $\zzz_2$, the upper-level player has to carefully choose their strategy. 
More precisely, given dimensions $n,m,p,q,r,s \in\N$ and differentiable functions $\fupper$, $\flower$, $\gupper$, $\hupper$, $\glower$, and $\hlower$ from domain $\R^{n+m}$ to range $\R$, $\R$, $\R^p$, $\R^q$, $\R^r$, and $\R^s$, respectively, the described \emph{Bilevel Optimisation Program} can be written as

\begin{equation}
\label{eq:bilevel}
\tag{BP}
\begin{aligned}
\minimise_{\zzz_1,\zzz_2} \quad 
&\fupper(\zzz_1,\zzz_2)\\
\subjectto \quad
    & \gupper(\zzz_1) \geq 0,\\ 
    & \hupper(\zzz_1) = 0,\\ 
    & \zzz_2 \in S(\zzz_1)\coloneqq\argmin_{\zzz_2'}
\left\{
\flower(\zzz_1,\zzz_2'):\quad \glower(\zzz_1,\zzz_2') \geq 0,\;\; \hlower(\zzz_1,\zzz_2') = 0
\right\},
\end{aligned}
\end{equation}
where the lower-level problem is the following program parametrised by $\zzz_1$:
\begin{align*}
    \minimise_{\zzz_2'}\ \flower(\zzz_1,\zzz_2')\quad 
    \subjectto\ \glower(\zzz_1,\zzz_2')\geq 0,\; \hlower(\zzz_1,\zzz_2')=0.
\end{align*}
Problem \eqref{eq:bilevel} has a wide range of applications (see, e.g., \cite{dempe2020bilevel} for more details) and our SVM hyperparameter tuning problem in Section~\ref{section:Application} is one of its practical applications. 

General NLP assumptions and CQs can be applied to the lower-level program.  For instance, we say that problem \eqref{eq:bilevel} has a convex lower-level program if and only if for all $\zzz_1\in \mathbb{R}^n$, the function $\flower(\zzz_1, \cdot)$ is convex,
the functions $g_{\text{low},i}(\zzz_1, \cdot)$, for $i=1, \ldots, r$, are concave (so that the constraints $g_{\text{low},i}(\zzz_1, \cdot)\geq0$ describe a convex set) and $\hlower(\zzz_1, \cdot)$ is affine linear. Additionally, we say that the \emph{Lower-Level Linear Independence Constraint Qualification} (LLICQ) and respectively the \emph{Lower-Level Mangasarian--Fromovitz Constraint Qualification} (LMFCQ) are satisfied at a point $(\zzz_1,\zzz_2)$ if and only if (\LICQ) and respectively (\MFCQ{}) are satisfied at the point $\zzz_2$ for the lower-level problem parametrised by $\zzz_1$.

A common technique towards solving problem \eqref{eq:bilevel} is to replace the lower-level problem with its Karush–Kuhn–Tucker (KKT) conditions, leading to the problem
\begin{equation}
\label{eq:KKT_reformation}
\tag{KKTR}
\begin{aligned}
\minimise_{\zzz_1,\zzz_2,\multg{}, \multh{}} \quad 
&\fupper(\zzz_1,\zzz_2)\\
\subjectto \quad   
    & \gupper(\zzz_1) \geq 0,\\
    & \hupper(\zzz_1) = 0,\\ 
    &
\begin{cases}
    \nabla_{\zzz_2} \flower(\zzz_1,\zzz_2) - \multg{}^{\top} \nabla_{\zzz_2} \glower(\zzz_1,\zzz_2) - \multh{}^{\top} \nabla_{\zzz_2} \hlower(\zzz_1,\zzz_2) = 0,\\
    0 \leq \multg{} \bot \glower(\zzz_1,\zzz_2)\geq 0,\\
    \hlower(\zzz_1,\zzz_2)= 0,\\
\end{cases} 
\end{aligned}
\end{equation}
where $\multg{}\in\R^r$ and $\multh{}\in\R^s$ correspond to the Lagrange multiplier associated with the lower-level constraint functions $g(\zzz_1, \zzz_2)$ and $h(\zzz_1, \zzz_2)$ at the point $(\zzz_1, \zzz_2)$ such that $\zzz_1$ is the upper-level selection and $\zzz_2$ is the corresponding lower-level optimal solution. 

Clearly, problem \eqref{eq:KKT_reformation} is a special class of problem (\MPCC{}). However, problem \eqref{eq:KKT_reformation} is not necessarily equivalent to the original problem \eqref{eq:bilevel}, as shown in \cite{DempeDutta2012}. Nevertheless, under suitable assumptions, the following global and local relationships are well-known.

\begin{theorem}[\cite{DempeDutta2012}]
\label{theorem:bilevel_kkt_reform_bijection}
Let the lower-level problem in \eqref{eq:bilevel} be convex.
\begin{description}
\item[(i)] 
Let $(\zzz_1^*, \zzz_2^*)$ be a global (resp. local) optimal solution of problem \eqref{eq:bilevel}
and assume that the LMFCQ holds at $(\zzz_1^*, \zzz_2^*)$. 
Then,
for each choice of Lagrange multipliers 
\begin{align*}
\label{eq:bilevel_multiplier_set}
(\multg{^*}, \multh{^*}) \text{ such that } &
\begin{cases}
    \nabla_{\zzz_2^*} \flower(\zzz_1^*,\zzz_2^*) - \multg{^*}^{\top}  \nabla_{\zzz_2^*} \glower(\zzz_1^*,\zzz_2^*) - \multh{^*}^{\top} \nabla_{\zzz_2^*} \hlower(\zzz_1^*,\zzz_2^*) = 0,\\
    0 \leq \multg{^*}  \bot \glower(\zzz_1^*,\zzz_2^*)\geq 0,\\
    \hlower(\zzz_1^*,\zzz_2^*)= 0,\\
\end{cases}
\end{align*}
the point $(\zzz_1^*, \zzz_2^*, \multg{^*}, \multh{^*})$ is a global (resp. local) optimal solution of problem \eqref{eq:KKT_reformation}. 
\item[(ii)] 
Let $(\zzz_1^*, \zzz_2^*, \multg{^*}, \multh{^*})$ be a global (resp. local) optimal solution of problem \eqref{eq:KKT_reformation}
and assume the LLICQ holds at $(\zzz_1^*, \zzz_2^*)$  (resp. LLICQ holds at $(\zzz_1, \zzz_2)$, for all $\zzz_2 \in S(\zzz_1)$, for all $\zzz_1: \gupper(\zzz_1)\geq0, \hupper(\zzz_1)=0$).
Then,
$(\zzz_1^*, \zzz_2^*)$ is a global (resp. local) optimal solution of problem \eqref{eq:bilevel}. 
\end{description}
\end{theorem}
This theorem establishes global and local relationships between problem \eqref{eq:bilevel} and its KKT/MPCC reformulation \eqref{eq:KKT_reformation}. This will be useful in Section~\ref{section:Application} when we construct a model for hyperparameter selection in support vector machines.

%% file: Supplementary/supplementary_21_model_as_bilevel.tex
\subsection{SVM hyperparameter tuning as a bilevel program}
\label{SM:model_as_bilevel}
% ...
\begin{figure}
\centering
    \includegraphics[width=0.66\textwidth]{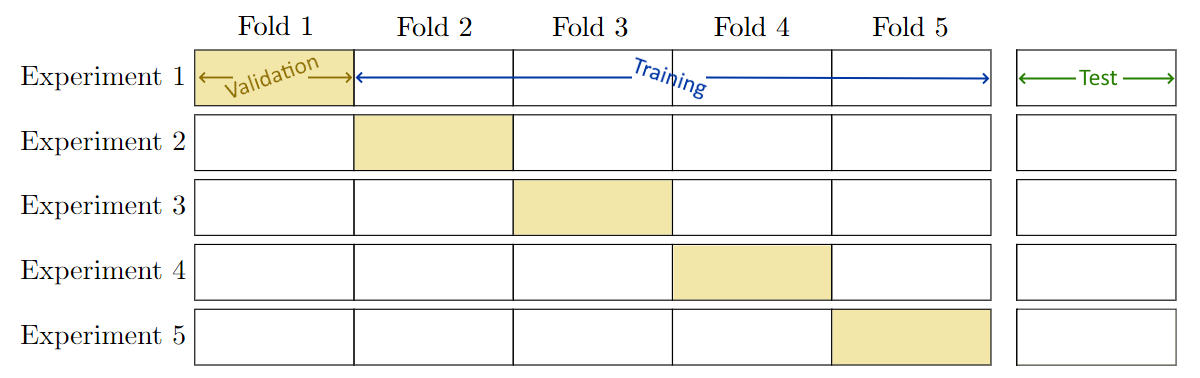}
    \caption{A visualisation of 5-fold cross-validation with dedicated test data kept aside.}
    \label{fig:k_fold_cross_val}
\end{figure}
% ...
With the above theory in mind, we can now model our SVM hyperparameter tuning problem as a bilevel program.
We begin in the primal space using the notation of Subsection~\ref{SM:primal_svm}.
The leader is to select hyperparameters $C$ and a feature mapping $\phi$ that minimise the average $K$-fold cross-validation error measured by the hinge loss.
The follower wishes to select optimal model parameters  $w^k$ and $b^k$ according to \eqref{eq:soft_margin}, for each fold $k$, that minimise the training error.
\begin{equation}
\tag{BP-primal}
\label{eq:bilvel_svm_primal}
\begin{aligned}
\minimise_{C,\phi,w^1,\dots,w^K,b^1,\dots,b^K}\quad
&\frac{1}{K}\sum_{k=1}^K \frac{1}{\bar{n}^k} \sum_{i=1}^{\bar{n}^k}
\max\left(0, 1-\bar{y}_i^k \left(w^{k\top} \phi \left(\bar{\xxx}_i^k \right)+b^k\right)\right)
\\
\subjectto \qquad\quad & C \geq 0,\\
&   \text{For each }  k=1,...,K, \text{ lower-level variables } w^k,b^k \text{ optimise }\\
&   \text{ \eqref{eq:soft_margin} parametrised by  } (x_i^k, y_i^k)_{i=1,\dots,n^k}.
\end{aligned}
\end{equation}

For the problem to be well-posed, as before, we will convert the lower-level variables into dual ones, by means of the kernel trick referring to the notation of Subsection~\ref{SM:dual_svm}.
This can be done by substituting $w^k=\sum_{j} \ablue{_j^k} y_j^k \phi\left(\xxx_j^k\right)$  into the objective function and applying the kernel trick $\bar{y}^k_i y^k_j \phi\left(\bar{\xxx}_i^k \right)^{\top}\phi\left(\xxx_j^k \right)=\bar{Q}(\gamma)^k_{ij}$. A few other technical arrangements with the bias $b^k$ will lead to the dual formulation.
\begin{equation}
\tag{BP-dual}
\label{eq:bilvel_svm_dual}
\begin{aligned}
\minimise_{C,\gamma, \ablue{^1},\dots,\ablue{^K}}\quad
&\frac{1}{K}\sum_{k=1}^K \frac{1}{\bar{n}^k} \sum_{i=1}^{\bar{n}^k}    \max\left(0, 1-\sum_{j=1}^{n^k} \ablue{_j^k} \bar{Q}(\gamma)^k_{ij} - \bar{y}_i^k b^k\left(\ablue{^k}, \gamma, C\right)\right)
\\
\subjectto \quad & C \geq 0,\gamma\geq0,\\
&   \text{For each }  k=1,...,K, \text{ lower-level variables } \ablue{^k} \text{ optimise }\\
&   \text{ \eqref{eq:SVMdual} parametrised by  } (x_i^k, y_i^k)_{i=1,\dots,n^k}.
\end{aligned}
\end{equation}

To obtain a smooth upper-level objective function, we remove the max operator by introducing slack variables $\zeta_i^k$ and pairs of constraints $\zeta_i^k\geq 0$ and $\zeta_i^k \geq 1-\sum_{j} \ablue{_j^k} \bar{Q}(\gamma)^k_{ij} - \bar{y}_i^k b^k\left(\ablue{^k}, \gamma, C\right)$. This must be done for each fold $k=1,...,K$  and for each validation example in those folds $i=1,\hdots,\bar{n}^k$.  To remove the embedded lower-level optimisation problem, we perform a KKT reformulation.  
This results in the MPCC model~\MODEL{} given in Subsection~\ref{section:model}.

We now show that the lower-level problem (\ref{eq:SVMdual}) satisfies the assumptions of Theorem~\ref{theorem:bilevel_kkt_reform_bijection}, and so, in the provided framework, problems \MODEL{} and~\eqref{eq:bilvel_svm_dual} will be locally equivalent.
\begin{proposition}
\label{proposition:convex_licq}
Problem~\eqref{eq:bilvel_svm_dual} has a convex lower level program, and satisfies the LLICQ at every point $\zzz \coloneqq [C, \gamma, \zeta, \ablue{}, \mlower{}, \mupper{}, \uu{}]^\top$ where $C>0$, $\gamma> 0$ and \eqref{eq:assumption_on_alpha} holds.
\end{proposition}
\begin{proof}
By Lemma~\ref{lemma:psd_kernel}%
, and more generally since every kernel matrix is positive semi-definite, the objective function of (\ref{eq:SVMdual}) is a
convex quadratic function. Each of the constraints is linear in $\ablue{}$. Therefore, we have convexity. 
For each index $i$, only one of the constraints $0\leq\ablue{_i}$ or $\ablue{_i}\leq C$ can be active. 
Therefore, each of our active inequality constraints is in a unique variable. 
Furthermore, by assumption \eqref{eq:assumption_on_alpha}, there exists some index $i$ such that the variable $\ablue{_i}$ is involved in the equality constraint $\sum\ablue{_i}y_i=0$ but none of the active inequality constraints. Thus, their derivatives are linearly independent, and problem~\eqref{eq:bilvel_svm_dual} fulfils the LLICQ.
\end{proof}

%% file: Supplementary/supplementary_22_positive_definite.tex
\subsection{Kernel matrix Q is positive definite}
\label{SM:Q_is_positive_definite}
% ...
\begin{lemma}
\label{lemma:psd_kernel}
Given distinct data points $x_1^k,\dots,x_{n^k}^k\in\R^d$ and labels $y_1^k,\dots,y_{n^k}^k\in\{-1,1\}$ for each $k=1,\dots,K$ and $\gamma> 0$,  we have:
\begin{enumerate}[leftmargin=*,label=(\roman*)]
    \item Each RBF kernel matrix $Q(\gamma)^k\in\R^{n^k\times n^k}$, defined by $Q(\gamma)^k_{ij} \coloneqq y_i^ky_j^k\exp\left(-\gamma\norm{x_i^k-x_j^k}^2\right)$ for~$i=1, \ldots, n^k$ and~$j=1, \ldots, n^k$,  is positive definite.
    \item The following block diagonal matrix $\Q(\gamma)$ formed from blocks $Q(\gamma)^1,\dots,Q(\gamma)^K$ is also positive definite:%
    \begin{equation}
    \label{eq:Q}
    \Q(\gamma)  \coloneqq  
    \begin{bmatrix}
        Q(\gamma)^1 & & 0\\
        & \ddots & \\
        0 &  & Q(\gamma)^K\\
    \end{bmatrix}.
    \end{equation}
\end{enumerate}
\end{lemma}
\begin{proof}
    The Gaussian $x\mapsto\exp(-\gamma\|x\|^2)$ is a positive definite function on $\R^d$ \cite[Theorem 6.10]{Wendland2004}. Therefore, 
    $P^k_{ij} \coloneqq \exp\left(-\gamma\norm{x_i^k-x_j^k}^2\right)$ is a positive definite matrix \cite[Definition 6.1]{Wendland2004}. Hence, for any non-zero vector $v \in\R^{n^k}$, by the positive definiteness of $P^k$, it holds that
    \begin{equation*}
        v^{\top} Q(\gamma)^k v = \sum_{i=1}^{n^k}\sum_{j=1}^{n^k} v_i v_j y_i^k y_j^k \exp\left(-\gamma\norm{x_i^k-x_j^k}^2\right) = (vy)^{\top} P^k (vy)>0,
    \end{equation*}
     where $vy\in\R^{n^k}$ is the non-zero vector $(vy)_i=v_iy_i^k$ for $i=1, \ldots, n^k$.  Thus, $Q(\gamma)^k$ is positive definite and so has only positive eigenvalues. Applying either \cite[Definition 1.59, Property a]{George2024} or Schur determinant formula \cite{Zhang2006} to the characteristic polynomial of $\Q(\gamma)$ we get:
\begin{equation*}
    \det\left(\Q(\gamma)-\lambda\id_{|I|}\right)=\prod_{k=1}^{K}\det\left(Q(\gamma)^k-\lambda\id_{n^k}\right).
\end{equation*}
 Hence,  $\Q(\gamma)$ has all positive eigenvalues and is therefore positive definite.
\end{proof}

%% file: Supplementary/supplementary_30_futher_solution_methods.tex
\section{Further details of the solution methods}
\label{SM:futher_solution_methods}
Here we provide the equations, algorithms and proofs that were excluded from Section~\ref{section:solution_methods} of the main article for brevity. 

%% file: Supplementary/Supplementary_31_proof_global.tex
\subsection{Proof of Theorem~\ref{theorem:penalisation_global_iterates}}
\label{SM:proof_global}
% ...
\begin{proof} 
We first want to show that $\left(G(\zzz^t)^\top H(\zzz^t)\right)_{t\in\N}$ is a non-increasing sequence.
We will apply the classical trick of taking a difference of terms at two indices.
Let $\overline{t}>t$. Notice both $\zzz^{\overline{t}}$ and $\zzz^t$ are feasible for both
\penalisation{\penaltyparam{}^{\overline{t}}} and \penalisation{\penaltyparam{}^t}.
Since $\zzz^t$ minimises $f(\zzz)  + \penaltyparam{}^{t} G(\zzz)H(\zzz)$ and $\zzz^{\overline{t}}$ minimises $f(\zzz)  + \penaltyparam{}^{\overline{t}} G(\zzz)H(\zzz)$ we can conclude that 
\begin{subequations}
\begin{align}
\label{eq:penalty_t}
    f(\zzz^{\overline{t}}) + \penaltyparam{}^{t} G(\zzz^{\overline{t}})H(\zzz^{\overline{t}}) &\geq f(\zzz^{t})  + \penaltyparam{}^{t} G(\zzz^{t})H(\zzz^{t}),
    \\
\label{eq:penalty_t_prime}
    f(\zzz^{t}) + \penaltyparam{}^{\overline{t}} G(\zzz^{t})H(\zzz^{t})  &\geq f(\zzz^{\overline{t}})  + \penaltyparam{}^{\overline{t}} G(\zzz^{\overline{t}})H(\zzz^{\overline{t}}).
\end{align}
\end{subequations}
Adding the two inequalities (\ref{eq:penalty_t}) and (\ref{eq:penalty_t_prime}) gives
\begin{align*}
    (\penaltyparam{}^{\overline{t}}-\penaltyparam{}^{t})\left(G(\zzz^{t})H(\zzz^{t})-G(\zzz^{\overline{t}})H(\zzz^{\overline{t}}) \right)\geq 0.
\end{align*}
Since $(\penaltyparam{}^{\overline{t}}-\penaltyparam{}^{t})>0$, it follows that $G(\zzz^{t})H(\zzz^{t})\geq G(\zzz^{\overline{t}})H(\zzz^{\overline{t}})$.

Next, we show that $\lim_{t\rightarrow\infty}\left( G(\zzz^t)^\top H(\zzz^t) \right)=0$. Let $\bar{\zzz}$ be any feasible point of \MPCC{}%
, which exists by assumption. Let $\zzz^1$ be the minimiser of the penalised problem for
$\penaltyparam{}^1$, i.e. the first iterate. Towards a contradiction, suppose that there exists an $\epsilon>0$ such that for all $T\in\N$ we have $G(\zzz^t)^\top H(\zzz^t)>\epsilon$ for some $t>T$. Under that premise,
and since $\penaltyparam{}^t\rightarrow\infty$, we can choose an arbitrarily large $t$ such that $\penaltyparam{}^t \geq \frac{1}{\epsilon}| f(\bar{\zzz})-f(\zzz^1)|+2$ and $G(\zzz^t)H(\zzz^t)>\epsilon$. Consider the infimum over the objective values of problem \MPCC{}:
\begin{subequations}
\begin{align}
\label{eq:inf_mpcc}
&\inf\left\{f(\zzz):\ g(\zzz)\geq0,\ h(\zzz)=0,\ 0\leq G(\zzz) \bot H(\zzz) \geq 0 \right\}\\
\label{eq:inf_penal}
&\geq\inf\left\{f(\zzz)+\penaltyparam{}^t G(\zzz)H(\zzz):\ g(\zzz)\geq 0,\ h(\zzz)=0,\ G(\zzz)\geq 0,\ H(\zzz) \geq 0 \right\}\\
\notag
&=f(\zzz^t)+\penaltyparam{}^t G(\zzz^t)H(\zzz^t)\\
\label{eq:ineq_second}
&\geq f(\zzz^1)+\penaltyparam{}^t G(\zzz^t)H(\zzz^t)\\
\notag
&\geq f(\zzz^1)+\penaltyparam{}^t \epsilon \\
\notag
&\geq f(\zzz^1)+| f(\bar{\zzz})-f(\zzz^1)|+2\epsilon \\
\label{eq:ineq_conclusion}
&>f(\bar{\zzz}).
\end{align}
\end{subequations}

 It is helpful to remember that $t$ here is chosen according to $\epsilon$. The first inequality \eqref{eq:inf_penal} comes from the fact that any feasible point of \MPCC{} is also feasible for \penalisation{\penaltyparam{}^t}, with the same objective value since its penalty term vanishes. The
 equality comes from the fact that $\zzz^t$ is a global minimum of $f(\zzz)+\penaltyparam{}^t G(\zzz)H(\zzz)$%
 , and the second inequality~\eqref{eq:ineq_second} from $f(\zzz^t)\geq f(\zzz^1)$, which follows from~\eqref{eq:penalty_t}, with $t$ and $\overline{t}$ replaced by $1$ and $t$, respectively, since $G(\zzz^t)^\top H(\zzz^t)$ is non-increasing.

The conclusion of \eqref{eq:ineq_conclusion} implies that $\bar{\zzz}$ achieves a smaller objective value than the infimum over objectives of the MPCC feasible points in problem \eqref{eq:inf_mpcc}. But this contradicts the MPCC feasibility of $\bar{\zzz}$.  Therefore we conclude that $\lim_{t\rightarrow\infty}\left( G(\zzz^t)^\top H(\zzz^t) \right)=0$. By the continuity of the functions $G$ and $H$, we have $ G(\zzz^*)^\top H(\zzz^*)=0$.
\end{proof}

%% file: Supplementary/supplementary_33_proof_relaxation.tex

\begin{figure}
    \centering
    \includegraphics[width=\linewidth]{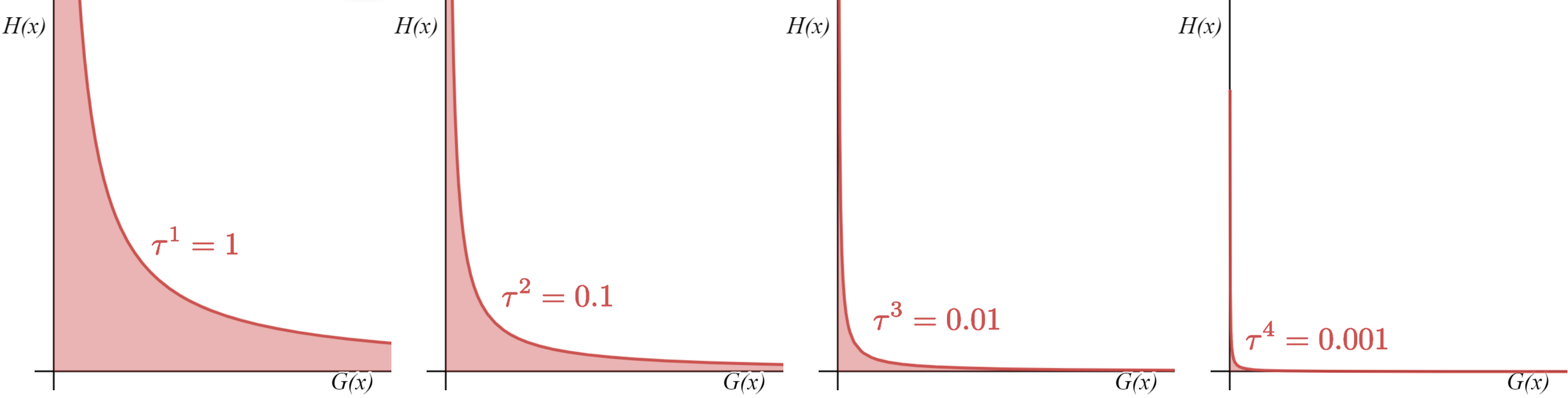}
    \caption{The feasible set for the relaxation subproblems with (left to right) decreasing relaxation parameters.}
    \label{fig:relaxation_sequence}
\end{figure}

%% file: Supplementary/supplementary_33b_counterexample_relaxation.tex
\subsection{Scholtes relaxation under MPCC-MFCQ-R alone}
\label{SM:counterexample_relaxation}
% ...
The following example, which is due to an anonymous reviewer, shows that the \MPCCMFCQR{} alone is not sufficient for the conclusion of Theorem~\ref{theorem:relaxation_convergence}. A sequence of global minimisers, and hence exact KKT points, of the relaxed problems converges to a point at which the \MPCCMFCQR{} holds, but which is not even W-stationary. The example also illustrates the roles of the \MPCCMFCQC{} in part~(i) and of the additional assumption in part~(ii) of the theorem.

\begin{example}
\label{example:counterexample_relaxation}
Consider the following MPCC in the variables $\zzz=(x,y,w,s)\in\R^4$:
\begin{equation}
\label{eq:counterexample_relaxation}
\begin{aligned}
\minimise_{x,y,w,s} \quad
& -w
\\
\subjectto \quad
& s=0,\\
& 0\leq x\bot y\geq 0,\\
& 0\leq x-w^2\bot y-w^2\geq 0.
\end{aligned}
\end{equation}
That is, $p=0$, $q=1$ and $r=2$ with $f(\zzz)\coloneqq-w$, $h_1(\zzz)\coloneqq s$, $G_1(\zzz)\coloneqq x$, $H_1(\zzz)\coloneqq y$, $G_2(\zzz)\coloneqq x-w^2$ and $H_2(\zzz)\coloneqq y-w^2$.
\end{example}

\noindent\textit{The limit point.}
If $\zzz$ is feasible for~\eqref{eq:counterexample_relaxation}, then $x=0$ or $y=0$ by the first complementarity constraint, and so $w^2\leq\min\{x,y\}=0$ by the second. Hence, the feasible set is $\{(x,y,0,0): 0\leq x\bot y\geq 0\}$, on which $f\equiv0$, so that $\zzz^*\coloneqq(0,0,0,0)$ is a global minimiser. Both complementarity pairs are biactive at $\zzz^*$; i.e., $I^h=\{1\}$, $I^G(\zzz^*)=I^H(\zzz^*)=\emptyset$ and $I^{GH}(\zzz^*)=\{1,2\}$. Denoting by $e_x,e_y,e_w,e_s$ the unit vectors of $\R^4$ corresponding to the four variables, the gradients at $\zzz^*$ are
\begin{equation*}
\nabla f(\zzz^*)=-e_w,\quad
\nabla h_1(\zzz^*)=e_s,\quad
\nabla G_1(\zzz^*)=\nabla G_2(\zzz^*)=e_x,\quad
\nabla H_1(\zzz^*)=\nabla H_2(\zzz^*)=e_y.
\end{equation*}
The \MPCCMFCQR{} holds at $\zzz^*$: the family in~\eqref{eq:MPCCMFCQR_linearly_independent_vectors} reduces to $\{e_s\}$, which is linearly independent, and $\ddd\coloneqq(1,1,0,0)$ satisfies~\eqref{eq:MPCCMFCQR_direction_d}, since $\nabla h_1(\zzz^*)^\top\ddd=0$ and $\nabla G_i(\zzz^*)^\top\ddd=\nabla H_i(\zzz^*)^\top\ddd=1>0$ for $i=1,2$.
However, $\nabla f(\zzz^*)=-e_w$ does not lie in the span of $e_x,e_y,e_s$, which contains all the constraint gradients at $\zzz^*$. Therefore, condition~\eqref{eq:stationarity_a} cannot hold, and $\zzz^*$ is not W-stationary and, in particular, not C-stationary.
In line with Theorem~\ref{theorem:relaxation_convergence}(i), the \MPCCMFCQC{} fails at $\zzz^*$: the scalars $\multG{_1}=1$ and $\multG{_2}=-1$, with all others zero, satisfy~\eqref{eq:MPCCMFCQC_multipliers}, since $\multG{_i}\multH{_i}=0$ for $i=1,2$. Note also that, as $\zzz^*$ is a global minimiser, the \MPCCMFCQR{}, unlike the \MPCCMFCQC{}, does not ensure that local minimisers of \MPCC{} are W-stationary.

\medskip
\noindent\textit{The relaxed problems.}
For $\delta>0$, consider the relaxed problem~\relaxation{\relaxparam{}} with $\relaxparam{}\coloneqq\delta^4$ and the point
\begin{equation*}
\zzz(\delta)\coloneqq(\delta^2,\delta^2,\delta,0).
\end{equation*}
This point is feasible for~\relaxation{\relaxparam{}}, since $G_1(\zzz(\delta))=H_1(\zzz(\delta))=\delta^2>0$ with $G_1(\zzz(\delta))H_1(\zzz(\delta))=\delta^4=\relaxparam{}$, and $G_2(\zzz(\delta))=H_2(\zzz(\delta))=0$. It is in fact a global minimiser of~\relaxation{\relaxparam{}}: every feasible point of~\relaxation{\relaxparam{}} satisfies $x\geq w^2$ and $y\geq w^2$, as well as $xy\leq\delta^4$, so that $w^4\leq xy\leq\delta^4$ and hence $f(\zzz)=-w\geq-\delta=f(\zzz(\delta))$.
Moreover, $\zzz(\delta)$ satisfies~\eqref{eq:KKT_of_relaxed} with $\epsilon=0$ for the multipliers
\begin{equation*}
\multh{_1}=\multG{_1}=\multH{_1}=\multGH{_2}=0,\qquad
\multG{_2}=\multH{_2}=\frac{1}{4\delta},\qquad
\multGH{_1}=\frac{1}{4\delta^3}.
\end{equation*}
Indeed, all multipliers are non-negative, the complementarity slackness conditions in~\eqref{eq:KKT_of_relaxed_b} hold since $G_2(\zzz(\delta))=H_2(\zzz(\delta))=0$ and the constraint $G_1(\zzz)H_1(\zzz)\leq\relaxparam{}$ is active, and, with $\nabla G_2(\zzz(\delta))=e_x-2\delta e_w$, $\nabla H_2(\zzz(\delta))=e_y-2\delta e_w$ and $G_1\nabla H_1+H_1\nabla G_1=\delta^2(e_x+e_y)$ at $\zzz(\delta)$, the left-hand side of~\eqref{eq:KKT_of_relaxed_a} is the norm of
\begin{equation*}
-e_w
-\frac{1}{4\delta}\left(e_x-2\delta e_w\right)
-\frac{1}{4\delta}\left(e_y-2\delta e_w\right)
+\frac{1}{4\delta^3}\,\delta^2\left(e_x+e_y\right)
=\veczero{}.
\end{equation*}

\medskip
\noindent\textit{The sequence of iterates.}
Let $(\delta^t)_{t\in\N}$ be any strictly decreasing sequence of positive numbers tending to zero, and let $\relaxparam{}^t\coloneqq(\delta^t)^4$ and $\zzz^t\coloneqq\zzz(\delta^t)$. Since each $\zzz^t$ is an exact KKT point of~\relaxation{\relaxparam{}^t}, it is an $\epsilon^t$-KKT point for any sequence of tolerances $\epsilon^t\downarrow0$, and $\zzz^t\rightarrow\zzz^*$. Thus, all the assumptions of Theorem~\ref{theorem:relaxation_convergence} other than the constraint qualifications are satisfied, and the \MPCCMFCQR{} holds at the limit $\zzz^*$, yet $\zzz^*$ is not C-stationary.
The additional assumption of Theorem~\ref{theorem:relaxation_convergence}(ii) fails, since $\multGH{_1^t}G_1(\zzz^t)=\multGH{_1^t}H_1(\zzz^t)=1/(4\delta^t)\rightarrow\infty$.

It is instructive to follow the proof of Theorem~\ref{theorem:relaxation_convergence} for this sequence. The augmented multipliers are $\multG[hat]{_1^t}=-\multGH{_1^t}H_1(\zzz^t)=-1/(4\delta^t)$ and $\multH[hat]{_1^t}=-\multGH{_1^t}G_1(\zzz^t)=-1/(4\delta^t)$, since the constraint $G_1(\zzz)H_1(\zzz)\leq\relaxparam{}^t$ is active, and $\multG[hat]{_2^t}=\multH[hat]{_2^t}=1/(4\delta^t)$, since $G_2(\zzz^t)=H_2(\zzz^t)=0$. Hence, $\|\aaa_t\|=1/(2\delta^t)\rightarrow\infty$, and the normalised multipliers converge to
\begin{equation*}
\left(\multG[hat]{'_1},\multG[hat]{'_2}\right)=\left(\multH[hat]{'_1},\multH[hat]{'_2}\right)=\left(-\tfrac{1}{2},\tfrac{1}{2}\right),
\end{equation*}
for which~\eqref{eq:r_linearly_dependence} reads $-\tfrac12e_x+\tfrac12e_x-\tfrac12e_y+\tfrac12e_y=\veczero{}$. The coefficients of the biactive index $i=1$ are negative, so this does not contradict the \MPCCMFCQR{}, whereas $\multG[hat]{'_i}\multH[hat]{'_i}=\tfrac14\geq0$ for $i=1,2$, so that it does contradict the \MPCCMFCQC{}, as in the proof of part~(i).

%% file: Supplementary/supplementary_34_m_stationarity.tex
\subsection{Semismooth Newton method}
\label{SM:m_stationarity_system}
% ...
During our experiments, we compared to the semismooth Newton method. This method writes the M-stationary conditions of problem \MPCC{} as a system of semismooth equations, which can be solved with a Newton-style method. This framework is proposed in \cite{Harder2021}.
It requires, in the first instance, the use of an NCP function that, given inputs $a,b$ returns zero if and only if $a,b\geq 0$ and $ab=0$. Two examples of such functions are the minimum and Fischer--Burmeister functions defined by
\begin{align*}
&\min(a,b) \coloneqq  \begin{cases}
    a \text{ if } a\leq b,\\
    b \text{ otherwise};
\end{cases}&
&\FischerBurmeister(a,b)  \coloneqq  \sqrt{a^2 + b^2} - a - b.&
\end{align*}

The second step of building the framework consists of introducing the M-stationarity function $\MStationaryFunc:\R^4 \rightarrow \R^2$, which has the property that its zero level set $\{(G_i,H_i,\multG{_i},\multH{_i}): \MStationaryFunc(G_i,H_i,\multG{_i},\multH{_i})=0\}$ is equal to the set of points that satisfy the M-stationary conditions; i.e., 
\begin{equation}
    \label{eq:MStationarityFunction}
     \MStationaryFunc(G_i, H_i, \multG{_i}, \multH{_i}) = 0 \quad \iff \quad
     \begin{cases}
          0 \leq G_i \bot H_i \geq 0,\\
          G_i \multG{_i} =0,\\
          H_i \multH{_i} =0,\\
          \left( \multG{_i},\multH{_i} \geq 0  \text{ or }  \multG{_i}\multH{_i} =0 \right).
     \end{cases}
\end{equation} 
 For an example of such a function, see \cite[p. 1470, Section 3.2]{Harder2021}. The full system of equations to be solved to compute M-stationary points for problem \MPCC{} becomes  
\begin{align}
    \label{eq:system}
    F(\zzz, \multg{}, \multh{}, \multG{}, \multH{}) \; \coloneqq  \;
    \begin{bmatrix}
    \nabla_{\zzz}\mathcal{L}(\zzz, \multg{}, \multh{}, \multG{}, \multH{})\\
    \left[\FischerBurmeister(g_i(\zzz),\multg{_i})\right]_{i=1,\dots,p}\\
    \left[h_i(\zzz)\right]_{i=1,\dots,q}\\
    \left[\MStationaryFunc( G_i(\zzz), H_i(\zzz), \multG{_i}, \multH{_i})\right]_{i=1,\dots,r}
    \end{bmatrix}=0
\end{align}
with the function $F$ being semismooth (see \cite{Harder2021} and references therein for the definition), while $\mathcal{L}$ is the Lagrangian function as defined in (\ref{eq:MPCC_Lagragian}). Note that here, the $\FischerBurmeister$ function may be replaced by $\min$ or any other NCP function.

It can easily be shown that $(\zzz^*, \multg{^*}, \multh{^*}, \multG{^*}, \multH{^*})$ solves equation \eqref{eq:system} if and only if $\zzz^*$ is an M-stationary point of the \MPCC{} with Lagrange multipliers ($ \multg{^*}, \multh{^*}, \multG{^*}, \multH{^*}$). To solve the system, we use a semismooth Newton method with an Armijo line search and the Newton derivative to address the involved non-smoothness, according to \cite[Algorithm 5.1]{Harder2021}.

The system of equations, when applied to our SVM hyperparameter tuning model \MODEL{}, is
\begin{align}
\label{eq:FunctionInFull}
F_{\mathcal{M}}(\zzz,\multip{})=
\begin{bmatrix}
\nabla_{\zzz}          \mathcal{L}(\zzz,\multip{})\\
\FischerBurmeister(C,\multip{^{C}_{}})\\
\FischerBurmeister(\gamma,\multip{^{\gamma}_{}})\\
\left[
\FischerBurmeister\left(\zeta_i^{k}, \multip{^{\zeta}_{(k,i)}}\right)
\right]_{(k,i)\in\bar{I}}
\\
\left[
\FischerBurmeister\left(Z(\zzz)_i^{k}, \multip{^{Z}_{(k,i)}}\right)
\right]_{(k,i)\in\bar{I}}
\\
\left[
\theta(\zzz)^k_i
\right]_{(k,i)\in I}
\\
\left[
\varphi(\zzz)^k
\right]_{k=1,\dots,K}
\\
\left[
\MStationaryFunc\left(
\ablue{^k_i},\mlower{^k_i},\multip{^{G_1}_{(k,i)}}, \multip{^{H_1}_{(k,i)}}
\right)\right]_{(k,i)\in I}
\\
\left[
\MStationaryFunc\left(
C-\ablue{^k_i},\mupper{^k_i},\multip{^{G_2}_{(k,i)}}, \multip{^{H_2}_{(k,i)}}
\right)\right]_{(k,i)\in I}
\end{bmatrix}
\end{align}
where $\FischerBurmeister:\R^2\rightarrow\R$ is the Fischer--Burmeister function and $\MStationaryFunc:\R^4\rightarrow\R^2$ is the M-stationarity function, the derivatives of the Lagrangian are written explicitly in Subsection~\ref{SM:Lagrangian} and the functions $Z, \theta, \varphi$ are written in (\link{eq:constraint_function_definitions}{\ref*{eq:constraint_function_definitions}a,b,c}).

%% file: Supplementary/supplementary_35_lagrangian.tex
\subsection{Explicit formula of the Lagrangian and its derivatives}
\label{SM:Lagrangian}
% ...
Early in the paper, we dealt with a compact form of the MPCC problem \MODEL{}. However, when it comes to implementing solution methods, we need accurate formulas for the Lagrangian and its first and second-order derivatives. We use the notation defined in equations \eqref{eq:flat_index} and \eqref{eq:constraint_function_definitions}. In this supplementary material, we present those explicit formulas in full, beginning with the Lagrangian
\begin{multline}
\mathcal{L}(\zzz, \multip{}) =
\sum_{(k,i)\in\bar{I}} \zeta_i^k - 
\multip{^{C}_{}} C - 
\multip{^{\gamma}_{}} \gamma - 
\sum_{(k,i)\in\bar{I}} \multip{^{\zeta}_{(k,i)}} \zeta^k_i - 
\sum_{(k,i)\in\bar{I}} \multip{^{Z}_{(k,i)}} Z(\zzz)^k_i - 
\sum_{(k,i)\in{I}} \multip{^{\theta}_{(k,i)}} \theta(\zzz)^k_i \\ -
\sum_{(k,i)\in{I}} \multip{^{G_1}_{(k,i)}} \ablue{^k_i} - 
\sum_{(k,i)\in{I}} \multip{^{H_1}_{(k,i)}} \mlower{^k_i} - 
\sum_{(k,i)\in{I}} \multip{^{G_2}_{(k,i)}} (C-\ablue{^k_i}) - 
\sum_{(k,i)\in{I}} \multip{^{H_2}_{(k,i)}} \mupper{^k_i} - 
\sum_{k=1,\dots,K} \multip{^{\varphi}_{k}} \varphi(\zzz){^k};
\end{multline}
where the decision variables and Lagrange multipliers are
\begin{align*}
\zzz&= \left[
    C,
    \gamma,
    \zeta,
    \ablue{},
    \mlower{},
    \mupper{},
    \uu{}
\right]\in \R^{3(K-1)n+n+K+2};\\
\multip{}&= \left[
    \multip{^{C}_{}},
    \multip{^{\gamma}_{}},
    \multip{^{\zeta}},
    \multip{^{Z}}, 
    \multip{^{\theta}},
    \multip{^{G_1}},
    \multip{^{H_1}},
    \multip{^{G_2}},
    \multip{^{H_2}},
    \multip{^{\varphi}}
\right]\in \R^{5(K-1)n+2n+K+2};
\end{align*}
with dimensions in terms of the total number of data points $n$ and the number of folds $K$.

The first order derivative with respect to all input variables is a vector to vector function  $\nabla\mathcal{L}:\R^{L}\rightarrow\R^{L}$ where $L=8(K-1)n+3n+2K+4$. We define it in terms of each partial derivative. 
\begin{equation*}
\label{eq:derivativesOfLagrangian}
\begin{split}
\nabla_{C} \mathcal{L}(\zzz,\multip{}) =& 
-\multip{^{C}_{}} 
- \sum_{k=1}^{K} \sum_{i=1}^{\bar{n}^{k}} \multip{^{Z}_{(k,i)}} \bar{y}_i^k \nabla_C b
- \sum_{k=1}^{K}\sum_{i=1}^{n^{k}} \multip{^{G_2}_{(k,i)}};
\\
\nabla_{\gamma} \mathcal{L}(\zzz,\multip{}) =&
-\multip{^{\gamma}_{}} 
+ \sum_{k=1}^{K} \sum_{i=1}^{\bar{n}^{k}} \multip{^{Z}_{(k,i)}} 
\left(
\sum_{j=1}^{n^k} \ablue{_j^k} \norm{\bar{x}_i^k-x_j^k}^2 \bar{Q}(\gamma)^k_{ij} - \bar{y}_i^{k}\nabla_{\gamma}b
\right)\\
&
\phantom{-\multip{^{\gamma}_{}}} +
\sum_{k=1}^{K} \sum_{i=1}^{n^{k}} \multip{^{\theta}_{(k,i)}} 
\left( 
\sum_{j=1}^{n^k} \ablue{^k_j} \norm{x_i^k-x_j^k}^2 Q(\gamma)^k_{ij}
\right);
\\
\nabla_{\zeta^k_i} \mathcal{L}(\zzz,\multip{}) =& 1-\multip{^{\zeta}_{(k,i)}}-\multip{^{Z}_{(k,i)}};
\\
\nabla_{\ablue{_j^k}} \mathcal{L}(\zzz,\multip{}) =&
- \sum_{i=1}^{\bar{n}^k} \multip{^{Z}_{(k,i)}} 
\left(
\bar{Q}(\gamma)^k_{ij}+\bar{y}_i^{k}\nabla_{\ablue{_j^k}} b 
\right)
- \sum_{i=1}^{n^k} \multip{^{\theta}_{(k,i)}} Q(\gamma)^k_{ij}
-\multip{^{G_1}_{(k,j)}}
+\multip{^{G_2}_{(k,j)}}
-\multip{^{\varphi}_{k}} y_j^{k};
\\
\nabla_{\mlower{_i^k}} \mathcal{L}(\zzz,\multip{}) =& +\multip{^{\theta}_{(k,i)}}  -\multip{^{H_1}_{(k,i)}};
\\
\nabla_{\mupper{_i^k}} \mathcal{L}(\zzz,\multip{}) =& -\multip{^{\theta}_{(k,i)}}  -\multip{^{H_2}_{(k,i)}};
\\
\nabla_{\uu{^k}} \mathcal{L}(\zzz,\multip{}) =& - \sum_{i=1}^{n^k} \multip{^{\theta}_{(k,i)}} y_i^k.
\end{split}
\end{equation*}

The Hessian of the Lagrangian is:
\begin{equation*}
\nabla_{\zzz\zzz}^2\mathcal{L}(\zzz, \multip{}) =
\begin{bmatrix}
\phantom{\ } \nabla^2_{C C} \mathcal{L}          \phantom{\ } &
\phantom{\ } \nabla^2_{C \gamma} \mathcal{L}     \phantom{\ } & 
\phantom{\ } 0                                   \phantom{\ } & 
\phantom{\ } \nabla^2_{C \ablue{}} \mathcal{L}   \phantom{\ } & 
\phantom{\ } 0                                   \phantom{\ } & 
\phantom{\ } 0                                   \phantom{\ } &
\phantom{\ } 0                                   \phantom{\ } \\
\nabla^2_{\gamma C} \mathcal{L}        & 
\nabla^2_{\gamma \gamma} \mathcal{L}   & 0 & 
\nabla^2_{\gamma  \ablue{}} \mathcal{L}& 0 & 0 & 0 \\
0 & 0 & 0 & 0 & 0 & 0 & 0 \\
\nabla^2_{\ablue{} C} \mathcal{L}        & 
\nabla^2_{\ablue{} \gamma} \mathcal{L}   & 0 & 
\nabla^2_{\ablue{} \ablue{}} \mathcal{L} & 0 & 0 & 0 \\
0 & 0 & 0 & 0 & 0 & 0 & 0 \\
0 & 0 & 0 & 0 & 0 & 0 & 0 \\
0 & 0 & 0 & 0 & 0 & 0 & 0 \\
\end{bmatrix}
\end{equation*}

The partial second-order derivatives can be written explicitly as:
\begin{align*}
&\nabla^2_{C C} \mathcal{L}(\zzz,\multip{}) =
 - \sum_{k=1}^{K} \sum_{i=1}^{\bar{n}^{k}} \multip{^{Z}_{(k,i)}} \bar{y}_i^k \nabla^2_{C C} b\left(\ablue{^k}, \gamma, C\right);
\\
&\nabla^2_{\gamma\gamma} \mathcal{L}(\zzz,\multip{}) =
 -\sum_{k=1}^{K} \sum_{i=1}^{\bar{n}^{k}} \multip{^{Z}_{(k,i)}} \left(\sum_{j=1}^{n^k} \ablue{_j^k} \norm{\bar{x}_i^k-x_j^k}^4 \bar{Q}(\gamma)^k_{ij}+\bar{y}_i^{k}\nabla^2_{\gamma\gamma}b\right)
\\
&\phantom{\nabla^2_{\gamma\gamma} \mathcal{L}(\zzz,\multip{}) =}
-\sum_{k=1}^{K} \sum_{i=1}^{n^{k}} \multip{^{\theta}_{(k,i)}} \sum_{j=1}^{n^k} \ablue{_j^k} \norm{x_i^k-x_j^k}^4 Q(\gamma)^k_{ij};
\\
&\nabla^2_{\ablue{_j^k} \ablue{_l^k}} \mathcal{L}(\zzz,\multip{}) =
 \sum_{i=1}^{\bar{n}^{k}} \multip{^{Z}_{(k,i)}} \bar{y}_i^{k}\nabla^2_{\ablue{_j^k} \ablue{_l^k}}b\left(\ablue{^k}, \gamma, C\right);
\\
&\nabla^2_{C \gamma} \mathcal{L}(\zzz,\multip{}) =
 - \sum_{k=1}^{K} \sum_{i=1}^{\bar{n}^{k}} \multip{^{Z}_{(k,i)}} \bar{y}_i^k \nabla^2_{C \gamma} b\left(\ablue{^k}, \gamma, C\right);
\\
&\nabla^2_{C \ablue{_l^k}} \mathcal{L} (\zzz,\multip{}) =
 - \sum_{i=1}^{\bar{n}^{k}} \multip{^{Z}_{(k,i)}} \bar{y}_i^k \nabla^2_{C\ablue{_l^k}} b\left(\ablue{^k}, \gamma, C\right);
\\
&\nabla^2_{ \gamma \ablue{_j^k}} \mathcal{L}(\zzz,\multip{}) = 
 \sum_{i=1}^{\bar{n}^{k}} \multip{^{Z}_{(k,i)}}  \left( \norm{\bar{x}_i^{k}-x_j^{k}}^2 \bar{Q}(\gamma)^k_{ij}-\bar{y}_i^{k}\nabla^2_{\ablue{_j^k} \gamma}b \right)
+\sum_{i=1}^{n^k} \multip{^{\theta}_{(k,i)}} \norm{x_i^k-x_j^k}^2 Q(\gamma)^k_{ij}.
\end{align*}

%% file: Supplementary/supplementary_40_svm.tex
\section{Preliminaries on the Support Vector Machine (SVM)}
\label{SM:SVM}
% ...
The concept of support vector machine (SVM) was introduced by Cortes and Vapnik \cite{Cortes1995,Vapnik1999} and has since been a prolific research subject; see, e.g., \cite{Cervantes2020}. Here, we provide an introduction to the topic as well as a setup of notation that is used throughout the final half of the paper.

%% file: Supplementary/supplementary_41_primal_svm.tex
\subsection{Primal SVM}
\label{SM:primal_svm}
% ...
\begin{figure}
    \centering
    \includegraphics[width=\textwidth]{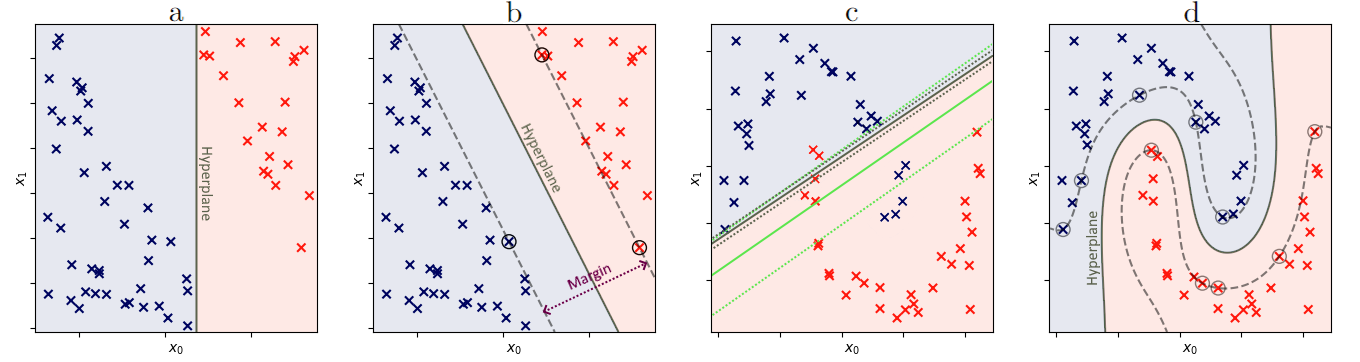}
    \caption{Each of the four plots shows blue and red data points in a two-dimensional space. Plots a and b (left) are the same linearly separable data set.  Plots c and d (right) are sampled from the Moons data set. The first plot shows a naive separating hyperplane. Plot b shows a more intelligent hyperplane with a maximum margin found using the (\ref{eq:hard_margin}). Plot c shows two possible hyperplanes in black and green found by the (\ref{eq:soft_margin}) with different choices of hyperparameter $C$. Plot d shows a nonlinear boundary with a maximum margin. We found this using the (\ref{eq:SVMdual}) attached with an RBF kernel.}
    \label{fig:SVM_margin}
\end{figure}
% ...
The SVM is defined  by a set of $n$ training examples, denoted by $(\xxx_i,y_i)$ for ${i=1,\dots,n}$, such that each example $i$ is a vector of features collected by $\xxx_i\in\R^d$ and a class label $y_i\in\left\{-1,+1\right\}$, where the aim is to find weights $\www\in \R^d$ and a bias $b\in\R$ describing a hyperplane $\{\xxx\in \R^d : \www^\top \xxx+b=0\}$ and a margin around it, $\{\xxx\in \R^d : -1 \leq \www^\top \xxx+b\leq 1\}$. The weights should be such that the positive label points lie above the margin ($\www^\top \xxx_i+b\geq 1 \quad \forall i:y_i=+1$) and the negative label points lie below the margin ($\www^\top \xxx_i+b\leq -1 \quad \forall i:y_i=-1$). This is equivalent to the constraint $1 - y_i(\www^\top \xxx_i+b)\leq 0$ for any $i$. Under these constraints, our goal is to maximise the width of the margin given by $\frac{2}{\|\www\|}$. Intuitively, this means new unseen data points have the most leeway to deviate beyond the range of the training data without crossing to the wrong side of the hyperplane; see Figure \link{fig:SVM_margin}{\ref*{fig:SVM_margin}a-b}. Once optimal weights $\www^*$ and bias $b^*$ are found, we can predict the label of new data using the classifier $\xxx\mapsto \sign\left(\www^{*\top}\xxx+b^*\right)$ where the sign function returns -1 for negative inputs and +1 otherwise. This problem can be formulated as 
\begin{equation}
\label{eq:hard_margin}
\tag{Hard-Margin-SVM}
\begin{aligned}
&\minimise_{\www\in\R^d,b\in\R} \quad & 
&\frac{1}{2}\www^{\top}\www&
\\
&\subjectto \quad & 
&1 - y_{i}\left(\www^\top \xxx_i+b\right) \leq 0&
&\text{ for } i = 1,\dots,n.&
\end{aligned}
\end{equation}

The weakness of the (\ref{eq:hard_margin}) becomes apparent when data is not linearly separable; the above problem is then infeasible. Thus, the (\ref{eq:hard_margin}) is rarely used in practice. A far more common approach is to use the hinge loss function 
\[
\hingeloss(\xxx_i,y_i)  \coloneqq  \max\left(0, 1 - y_{i}\left(\www^\top \xxx_i+b\right)\right),
\]
which measures how much example $(\xxx_i,y_i)$ violates the margin constraint for a given weight and bias. It takes the value zero when the constraint is satisfied. By replacing each constraint with a hinge loss penalty in the objective function, we can form a model that discourages but allows some points to violate the constraint in favour of ensuring a large margin.
A regularisation parameter $C$ can be introduced to instruct the balance between these two objectives, with a large $C$ ensuring that the margin constraint is less violated, whereas a smaller $C$ is more tolerant to such violations where it can find a larger margin. The new objective is now to minimise the function
\[
\frac{1}{2}\www^{\top}\www+C \sum_{i=1}^{n} \max\left(0, 1 - y_{i}\left(\www^\top \xxx_i+b\right)\right).
\]
The non-smooth max operator in this new objective function can be eliminated by introducing a variable $\xi_{i}$ for each example $(\xxx_i,y_i)$.  This leads to the soft margin formulation of the SVM problem
\begin{equation}
\label{eq:soft_margin}
\tag{Soft-Margin-SVM}
\begin{aligned}
&\minimise_{\www\in\R^d,\xibold\in\R^n,b\in\R}\quad &
&\frac{1}{2}\www^{\top}\www+C \sum_{i=1}^{n}{\xi_{i}}&
\\
&\subjectto \quad &
&\xi_{i} \geq 1-y_{i}\left(\www^\top \xxx_i+b\right)&
&\text{ for } i=1,\dots,n,&
\\
&&
& \xi_i \geq0&
&\text{ for } i=1,\dots,n.&
\end{aligned}
\end{equation}

\begin{figure}
    \centering
    \includegraphics[width=0.66\textwidth]{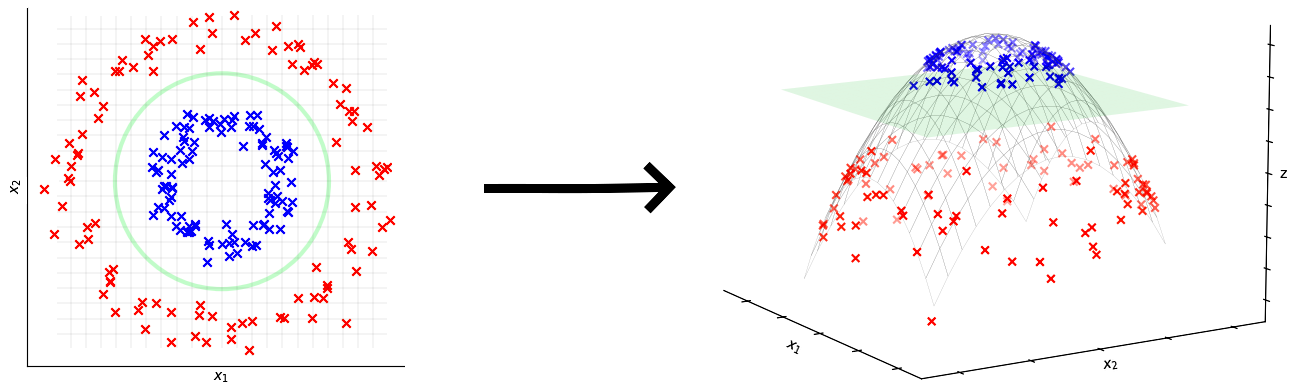}
    \caption{On the left, data in two dimensions that is not linearly separable is shown. On the right, the same data is lifted into a third dimension, where a separating hyperplane exists.}
    \label{fig:kernel_trick}
\end{figure}

Now, we wish to extend this model so that it can find nonlinear boundaries to separate data. To achieve this, we consider a mapping $\xxx_i \mapsto \phi(\xxx_i)$ that lifts the features of each data point $\xxx_i$, for $i=1, \ldots, n$, to a higher dimension. We can then find a hyperplane in this higher dimension, which corresponds to a nonlinear boundary in the original feature space; see Figure~\ref{fig:kernel_trick} for an illustration. However, it is often computationally hard to compute $\phi(\xxx_i)$. The solution is the so-called kernel trick, which allows us to find a separating hyperplane in the higher dimension without ever explicitly evaluating $\phi(\xxx_i)$. In our numerical experiments, we use the Radial Basis Function (RBF), which is defined by
\begin{equation}
    K_\gamma(\xxx_1, \xxx_2)  \coloneqq  \exp\left(-\gamma\norm{\xxx_1-\xxx_2}^2\right)
    \label{eq:rbf}
\end{equation}
(with hyperparameter $\gamma$) 
and maps features to an infinite-dimensional space, making the kernel trick necessary. Other  kernel functions commonly used in the literature \cite{Boser1992,Cervantes2020} include the linear, polynomial, and Sigmoid ones, respectively given as follows (for hyperparameters $\gamma, \delta$):
\[
K(\xxx_1, \xxx_2)  \coloneqq  \xxx_1^\top \xxx_2, \quad K_{\gamma,\delta}(\xxx_1, \xxx_2)  \coloneqq  (\gamma \xxx_1^\top \xxx_2 + \delta)^{2}, \; \mbox{ and }\;K_{\gamma,\delta}(\xxx_1, \xxx_2)  \coloneqq  \tanh\left(\gamma \xxx_1^\top \xxx_2 + \delta\right).
\]

To apply the kernel trick, we need to work with the dual of the soft margin SVM optimisation problem. This was first discovered by Guyon, Boser and Vapnik in 1992 \cite{Boser1992,Guyon1992} and is now standard. Nevertheless, for the paper to be self-contained, in the following Subsection~\ref{SM:dual_svm}, we provide the details of the conversion from the primal problem (\ref{eq:soft_margin}) to the dual (\ref{eq:SVMdual}).

%% file: Supplementary/supplementary_42_dual_svm.tex
\subsection{Dual SVM}
\label{SM:dual_svm}
% ...
The first step in dualising the SVM is to write the Lagrangian function $\mathcal{L}:\R^{d+3n+1}\rightarrow\R$ of the primal program (\ref{eq:soft_margin}). This is the sum of the original objective function plus each constraint scaled by corresponding Lagrange multipliers $\ablue{}\in\R^n$ and  $\etablue{}\in\R^n$. For the textbook theory on Lagrange functions, see Boyd and Vandenberghe \cite[Section 5.1.2]{boyd2004}.
% ...
\begin{equation}
\label{eq:LagrangeLL}
\mathcal{L}(\www,\xibold,b,\ablue{},\etablue{}) = 
\frac{1}{2} \www^{\top}  \www +
C \sum_{i=1}^{n}\xi_{i} +
\sum_{i=1}^{n} \ablue{_i} \Big( 1-\xi_{i}-y_{i}(\www^{\top}\phi(\xxx_{i})+b) \Big) -
\sum_{i=1}^{n} \etablue{_i} \xi_i
\end{equation}
% ...
We note three properties of the primal problem (\ref{eq:soft_margin}). First, the objective function is convex (in particular, positive quadratic). Secondly, the constraints are each linear in $\www,\xibold,b$, and so the feasible set is an (unbounded) convex polyhedron. Thirdly the feasible set has a non-empty interior since we can pick $\xi_i=\max\left(0, 1 - y_{i}\left(\www^{\top} \xxx_{i}+b\right)\right)+10^{-6}$ which strictly satisfies each inequality. Therefore the problem satisfies Slater's condition and strong duality. Any feasible point $\www,\xibold,b$ is a local minimiser of (\ref{eq:soft_margin}) (equivalent to a global minimiser by convexity) only if there exist Lagrange multipliers $\ablue{}, \etablue{}$ that satisfy the following first-order stationarity conditions
\begin{subequations}
\begin{align}
\label{eq:sc1}
&\nabla_{\www} \mathcal{L}(\www,\xibold,b,\ablue{},\etablue{})=  \www - \sum_{i=1}^{n}\ablue{_i} y_i \phi(\xxx_i) =  \veczero{},\\
\label{eq:sc2}
&\nabla_b \mathcal{L}(\www,\xibold,b,\ablue{},\etablue{})= -\sum_{i=1}^{n}\ablue{_i} y_i = 0, \\
\label{eq:sc3}
&\nabla_{\xibold} \mathcal{L}(\www,\xibold,b,\ablue{},\etablue{})= C\vecone{}-\ablue{} - \etablue{} = \veczero{}.
\end{align}
\end{subequations}

We now manipulate the Lagrangian (\ref{eq:LagrangeLL}) by substituting in the stationary conditions (\ref{eq:sc1}, \ref{eq:sc2}, \ref{eq:sc3}). This will give us a function of the Lagrangian in terms of $\ablue{}, \etablue{}$ assuming optimal choices for  $\www, \xibold, b$. The term $\sum_{i=1}^{n} (C - \ablue{_i} - \etablue{_i}) \xi_i$ vanishes due to (\ref{eq:sc3}). The term $\sum_{i=1}^{n} \ablue{_i} y_i b$ vanishes due to (\ref{eq:sc2}). And we eliminate $\www$ by substituting in (\ref{eq:sc1}). The conversion of the problem from primal space variables ($\www$, $\xibold$ and $b$) to dual space variables ($\ablue{}$ and $\etablue{}$) like this allows us to apply the kernel trick $\phi(\xxx_i)^{\top}\phi(\xxx_j)=K_\gamma\left(\xxx_i,\xxx_j\right)$. Below, we work through the manipulation line by line.
\begingroup
\allowdisplaybreaks
\begin{align*}
&\begin{aligned}
&&&\mathcal{L}^*(\ablue{}, \etablue{})&
\\
&=&
&\mathrlap{\frac{1}{2} \www^{\top}  \www +
C \sum_{i=1}^{n}\xi_{i} +
\sum_{i=1}^{n} \ablue{_i} \Big( 1-\xi_{i}-y_{i}(\www^{\top}\phi(\xxx_{i})+b) \Big) -
\sum_{i=1}^{n} \etablue{_i} \xi_i}&
\\
&=&
&\frac{1}{2} \www^{\top}  \www& + 
&\sum_{i=1}^{n} (C - \ablue{_i} - \etablue{_i}) \xi_i& + 
&\sum_{i=1}^{n} \ablue{_i}& -
&\sum_{i=1}^{n} \ablue{_i} y_i \www^{\top} \phi(\xxx_{i})& -
&\sum_{i=1}^{n} \ablue{_i} y_i b&
\\
&=&
&\frac{1}{2} \www^{\top}  \www& + 
&0& + 
&\sum_{i=1}^{n} \ablue{_i}& -
&\sum_{i=1}^{n} \ablue{_i} y_i \www^{\top} \phi(\xxx_{i})& -
&0&
\end{aligned}\\
&\begin{aligned}
&=&
&\frac{1}{2} \Big(\sum_{i=1}^{n}\ablue{_i} y_i \phi(\xxx_i)\Big)^{\top} \Big(\sum_{j=1}^{n}\ablue{_j} y_j \phi(\xxx_j)\Big)& + 
&\sum_{i=1}^{n} \ablue{_i}& - 
&\sum_{i=1}^{n} \ablue{_i} y_i \Big(\sum_{j=1}^{n}\ablue{_j} y_j \phi(\xxx_j)\Big)^{\top} \phi(\xxx_{i})&
\\
&=&
&\frac{1}{2} \sum_{i=1}^{n} \sum_{j=1}^{n} \ablue{_i} \ablue{_j} y_i y_j \phi(\xxx_{i})^{\top} \phi(\xxx_{j})& +
&\sum_{i=1}^{n} \ablue{_i}& -
&\sum_{i=1}^{n} \sum_{j=1}^{n} \ablue{_i} \ablue{_j} y_i y_j  \phi(\xxx_{i})^{\top} \phi(\xxx_{j})&
\\
&=&
&-\frac{1}{2} \sum_{i=1}^{n} \sum_{j=1}^{n} \ablue{_i} \ablue{_j} y_i y_j K_\gamma\left(\xxx_i,\xxx_j\right)& + 
&\sum_{i=1}^{n} \ablue{_i}&
\end{aligned}
\end{align*}
\endgroup

By considering the infimum of the $\{ \mathcal{L}(\www,\xibold,b,\ablue{},\etablue{})\}$  (\ref{eq:LagrangeLL}) we now write the dual of (\ref{eq:soft_margin}). For the textbook duality theory of convex programs, see Boyd and Vandenberghe \cite[Section 5.2]{boyd2004}.
\begin{align*}
\begin{split}
\maximise_{\ablue{}\in\R^n, \etablue{}\in\R^n}  \quad & \inf_{\www,\xibold\in\R^n,b\in\R}\Big\{ \mathcal{L}(\www,\xibold,b,\ablue{},\etablue{}) \Big\}
\\
\subjectto \quad & 0 \leq \ablue{_i} \quad \text{ for } i=1,\dots,n, \\
& 0 \leq \etablue{_i} \quad \text{ for } i=1,\dots,n. \\
\end{split}
\end{align*}

The above is impractical for us. We replace the infimum of the Lagrangian $(\inf\{ \mathcal{L}(\www,\xibold,b,\ablue{},\etablue{}) \})$ with the expression $\mathcal{L}^*(\ablue{}, \etablue{})$ derived from the stationary conditions above. Next, we multiply by $-1$ to convert from maximisation to a minimisation problem. Finally, the constraints $\ablue{}\geq0, \quad \etablue{} \geq0$ can be rearranged to $0 \leq \ablue{} \leq C$ by substituting out $\etablue{_i}=C-\ablue{_i}$ according to (\ref{eq:sc3}). Thus, the dual problem can be written as follows

\begin{align}
\label{eq:SVMdual}
\tag{Dual-SVM}
\begin{aligned}
&\minimise_{\ablue{} \in \R^n} \quad & 
&\frac{1}{2} \sum_{i=1}^{n} \sum_{j=1}^{n} \ablue{_i} \ablue{_j} y_i y_j K_\gamma({\xxx_i},{\xxx_j}) - \sum_{i=1}^{n} \ablue{_i}&
\\
&\subjectto \quad & 
&0 \leq \ablue{_i} \leq C \quad \text{ for } i=1,\dots,n,&
\\
&&
& \sum_{i=1}^n \ablue{_i} y_i =0,&
\end{aligned}
\end{align}
which will be the main focus of our attention for the remainder of this section.

The (\ref{eq:SVMdual}) is highly effective at learning to separate real-world data. Many machine-learning libraries have become proficient at solving this; one of the most widely used is the open-source library LIBSVM \cite{Chang2011}. At this point, it is worth discussing the link between the primal variable $(\www, b)$ and dual variable $\ablue{}$. For each index $i$, if $\ablue{_i}=0$, the solution represents a boundary that ensures data point $i$ lies on the correct side of the margin.  If $\ablue{_i}>0$, then the data point $i$ lies on the edge of the margin, within the margin or on the wrong side. Details of how to reconstruct the primal variables and further explanations are provided in Subsection~\ref{SM:reconstructing}.

%% file: Supplementary/supplementary_43_kkt_svm.tex
\subsection{KKT system of the dual SVM}
\label{SM:KKT_of_dual_SVM}
% ...
We conclude this subsection by stating the Karush-Kuhn-Tucker~(KKT) conditions for problem \eqref{eq:SVMdual}, which will be useful later. For the theory on KKT conditions, we refer again to Boyd and Vandenberghe \cite[Section 5.5.3]{boyd2004}. Since the program is convex and the constraints are linear, a point $\ablue{}$ is optimal if and only if there exist Lagrange multipliers $\mlower{}, \mupper{} \in \R^{n} , \uu{} \in \R $ such that 
\begin{equation}
\tag{KKT-SVM}
\label{eq:lowerLevelKKT}
\begin{aligned}
&\sum_{j=1}^{n}   \ablue{_j} Q(\gamma)_{ij} - \mlower{_i} + \mupper{_i} + \uu{} y_i = 1
&\text{ for } i=1,\dots,n,\\
& 0 \leq \ablue{_i} \bot \mlower{_i} \geq 0
&\text{ for } i=1,\dots,n,\\
& 0 \leq (C - \ablue{_i})  \bot \mupper{_i} \geq 0  \quad
&\text{ for } i=1,\dots,n,\\
&\sum_{i=1}^n \ablue{_i}y_i = 0,
\end{aligned}
\end{equation}
where $Q(\gamma)_{ij} \coloneqq y_i y_j K_\gamma({\xxx_i},{\xxx_j})$, $i=1, \ldots, n$, $j=1, \ldots, n$, denotes the so-called \textit{kernel} matrix which we write as a function of $\gamma$ as this will become a decision variable of our MPCC model.

%% file: Supplementary/supplementary_44_reconstruction.tex
\subsection{Reconstruction of primal variables from dual solutions}
\label{SM:reconstructing}
% ...
The weights $w$ and bias $b$ are primal variables. In order to use the kernel trick, we will be solving the dual optimisation problem (\ref{eq:SVMdual}) in terms of dual variables $\ablue{}$. This makes it less obvious what our final classifier should predict as we cannot explicitly calculate $x\mapsto \sign\left(w^{\top}\phi(x)+b\right)$. From the stationary conditions (\ref{eq:sc1}), we have
\begin{equation}
\label{eq:reconst_w}
 w = \sum_{i=1}^{n}\ablue{_i} y_i \phi(x_i).
\end{equation}
Using this formula (\ref{eq:reconst_w}), we can substitute out $w$ from the classifier. Our classifier should predict 
\[
x\mapsto \sign\left(\sum_{i=1}^n \ablue{_i}y_i \phi(x_i)^{\top} \phi( x)+b\right).
\]

We must work harder to find an expression for the bias $b$. Below, we write the complementarity slackness conditions for (\ref{eq:soft_margin})
\begin{align}
    \label{eq:cs1}
    \left(1-\xi_i - y_i (w^{\top}  \phi(x_i) +b)\right)\ablue{_i} &= 0,\\
    \label{eq:cs2}
    \xi_i \etablue{_i} &= 0.
\end{align}
Suppose $\ablue{_i}>0$. Then by complementarity slackness (\ref{eq:cs1}) we must have
\begin{equation}
\label{eq:csc1}
        1-\xi_i - y_i (w^{\top}  \phi(x_i) +b)=0 \quad \quad \text{ for } i\in\{1,...,n\}: \ablue{_i}>0.
\end{equation}
Suppose $\etablue{_i}>0$. Notice, this is equivalent to $\ablue{_i}<C$ by the stationary condition (\ref{eq:sc3}). Again, by complementarity slackness (\ref{eq:cs2}) we must have
\begin{equation}
\label{eq:csc2}
    \xi_i=0 \quad \quad \text{ for } i\in\{1,...,n\}: \ablue{_i}<C.
\end{equation}
Putting these two (\ref{eq:csc1}, \ref{eq:csc2}) together we get
\begin{equation*}
    1- y_i \left(w^{\top}  \phi(x_i) +b\right)=0 \quad \quad \text{ for } i\in\{1,...,n\}: 0<\ablue{_i}<C.
\end{equation*}
This gives us an important insight. If $\ablue{_i}$ takes a value between $0$ and $C$ then we can conclude $y_i(w^{\top}  \phi(x_i) + b)=1$. Geometrically, this means that example $i$ sits exactly on the edge of the margin. We call these ``free support vectors'', and they are our classifier's most crucial data points. Now recall that $y_i \in \{-1, +1\}$, so we can divide through by $y_i$, using $\frac{1}{y_i}= y_i$. This leads to
 \begin{equation*}
    b = y_i - w^{\top}  \phi(x_i) \quad \quad \text{ for } i\in\{1,...,n\}: 0<\ablue{_i}<C.
\end{equation*}

We can now substitute out $w$ using (\ref{eq:sc1}), giving:
 \begin{equation}
 \label{eq:biasSingle}
    b = y_i - \sum_{j=1}^{n}\ablue{_j} y_j \phi(x_j)^{\top}  \phi(x_i) \quad \quad \text{ for } i\in\{1,...,n\}: 0<\ablue{_i}<C.
\end{equation}
Theoretically, we can reconstruct the bias according to (\ref{eq:biasSingle}) from any support vector $i$ on the margin's edge, identified by $0<\ablue{_i}<C$. However, most libraries make the bias reconstruction calculation (\ref{eq:biasSingle}) for each such support vector and then take the average. With this in mind, let $\Omega(\ablue{_i}, C)$ be the indicator function that equals $1$ if $0<\ablue{_i}<C$ and $0$ otherwise. Whenever we write the bias as a function of dual variables $b(\ablue{},C)$ we mean the following reconstruction
\begin{equation}
\label{eq:bias}
    b= \frac{1}{\sum_{i=1}^{n}\Omega(\ablue{_i}, C)} \sum_{i=1}^{n}\Omega(\ablue{_i}, C)\left(y_i - \sum_{j=1}^{n} \ablue{_j} y_j K \left (x_i, x_j\right) \right).
\end{equation}

\begin{figure}
\begin{center}
\includegraphics[width=0.25\textwidth]{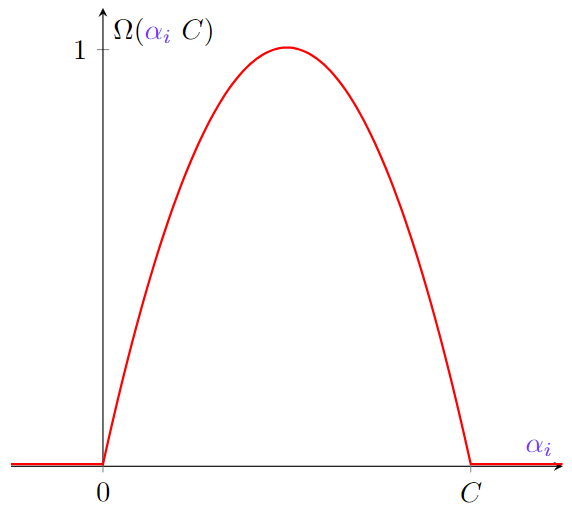}
\end{center}
    \caption{A graph of the function $\Omega(\ablue{_i}, C)$}
    \label{fig:omega}
\end{figure}

The algorithms we propose require a further adjustment. Firstly, because we want to avoid the combinatorial nature of a 0,1 indicator function. Secondly, it is hard to distinguish between $\ablue{_i}$ very close to zero and $\ablue{_j}$ that should be zero but carries some small numerical error. We, therefore, use an alternative indicator function as proposed by Anthony Dunn \cite{Coniglio2022}. 
\begin{equation}
\label{eq:omega}
\begin{aligned}
\Omega(\ablue{_i}, C) &= 
\begin{cases}
1-\left(\frac{2\ablue{_i}}{C}-1\right)^2 + \tau & \quad\text{if $0\leq \ablue{_i} \leq C$},\\
0 & \quad \text{otherwise.}\\
\end{cases}
\end{aligned}
\end{equation}

Similar to the original indicator, this takes the value zero for $\ablue{_i}\notin(0, C)$ and a positive value for $\ablue{_i}\in(0, C)$. Importantly, it assigns a much smaller value to $\ablue{_i}$ close to $0$ or $C$. With this new indicator, the reconstruction (\ref{eq:bias}) still averages across the examples where $0<\ablue{_i}<C$ and so is still theoretically correct but far more numerically stable.  

%% file: Supplementary/supplementary_50_futher_numerics.tex
\section{Further details on the numerical experiments}
\label{SM:futher_numerics}
% ...
To allow reproducibility of our experiments, this section of supplementary material provides full details on the starting point selection (Subsection~\ref{SM:starting_point}), datasets (Subsection~\ref{SM:datasets}) and derivative-free methods (Subsection~\ref{SM:derivative_free}) we employed.

Our experiments run on Southampton University's Iridis Compute Cluster \cite{Nicole1999}.
Each experiment is scheduled on a single Intel Xeon Gold 6138 CPU.
This CPU is comparable with most modern home PCs.

%% file: Supplementary/supplementary_51_starting_point.tex
\subsection{Starting point selection}
\label{SM:starting_point}
% ...
In this subsection, we discuss multi-starts and initialisation strategies for the algorithms in Section~\ref{section:solution_methods}, since each requires the user to provide a starting point $z^{0}=[C^0, \gamma^0, \zeta^0, \ablue{^0}, \mlower{^0}, \mupper{^0}, \uu{^0}]^\top$ in the context of our problem \eqref{eq:final-model}. The methods are deterministic but perform very differently depending on the choice of $z^0$. Recall also that the penalisation and relaxation methods converge to
C-stationary points, and to S-stationary points under an additional assumption at the biactive indices (Theorems~\ref{theorem:penalisation_convergence} and~\ref{theorem:relaxation_convergence}).  Such points may not be global or even local minima. Therefore, initialising the methods from a range of different starting points can increase the chance of finding one. However, it is well-known that under an MPCC-specific second-order sufficient condition \cite[Definition 2.8]{Ralph2004}, S-stationary points correspond to strict local minima \cite[Theorem 3.1]{Ralph2004}. We remind the reader here of the theoretical advantage of
converging to S-stationary points and their potential to be local minima.

\begin{table}[tbp]
\centering
\begin{tabular}{l|c c c}
    &\texttt{LL-Zero}  & \texttt{LL-Centre} & \texttt{LL-Feasible}\\
    \hline
    Relative computational cost & low & low & high \\
    Feasible $z^0$ for \MODEL{}? & no & no & yes \\
    Feasible $z^0$ for \penalisation{\penaltyparam{}^0}/\relaxation{\relaxparam{}^0}? & no & yes & yes 
\end{tabular} 
\caption{Summary of the lower-level initialisation strategies.}
\label{table:ll-init}
\end{table}

We propose a multi-start, where the upper-level variable $(C^0, \gamma^0)$ shall be initiated according to a grid spread uniformly across the log space such as $\{10^n:n=-5,\dots,4\}\times\{10^n:n=-3,\dots,6\}$, for example. This encourages the solution paths to explore as much of the hyperparameter space as possible. Note that this grid is far coarser than the derivative-free grid search we compare against later. Figure~\ref{fig:multistart} shows a two-by-four grid of starting points and the resulting path of iterates that the sequential relaxation algorithm produced for the Moons dataset.  Of these, three converge to the global minimum, while the remaining five converge to other stationary points.

Considering that the variable  $\zeta$ is only involved in two constraints, given an initialisation of $\ablue{^0}$, we can simply choose  
\[
\zeta_{\ i}^{0k}=\max\left(0, 1-\sum_{j=1}^{n^k} \ablue{_{\ j}^{0k}} \bar{Q}(\gamma)^k_{ij} - \bar{y}_i^k b^k\left(\ablue{^{0k}}, \gamma^0, C^0\right)\right) \;\; {\mbox{ for }\;\; k=1\dots K,\ i=1,\dots \bar{n}^k},
\] 
which ensures initial feasibility of the constraints $\zeta_{\ i}^{0k}\geq0$ and $Z_i^k(z^0)\geq0$; see \eqref{eq:Z}. 

It remains to be discussed how to initialise the lower-level variables $\ablue{^0}$, $\mlower{^0}$, $\mupper{^0}$ and $\uu{^0}$. We propose three strategies.
The first and simplest strategy shall be referred to as \texttt{LL-Zero}. This strategy initialises all of the lower-level variables to zero.  Notice this satisfies the complementarity constraints $ 0 \leq \ablue{_i^k} \bot \mlower{_i^k} \geq 0$ and $0 \leq  (C - \ablue{_i^k}) \bot \mupper{_i^k} \geq 0$ but violates the stationarity constraint $\theta^k_i(z)=0$  for each $(k,i)\in I$.  The strategy has no computational cost and can perform well but must be treated with caution as it is infeasible for all the programs \MODEL{}, \penalisation{\penaltyparam{}} and \relaxation{\relaxparam{}}.

We call the second strategy \texttt{LL-Centre}, which, for a given choice of $C^0$ and $\gamma^0$, selects values for the lower-level variables as follows:
\begin{align}
    &\ablue{^{0k}_{\ i}}=\frac{C^0}{2n^k}\sum_{j=1}^{n^k}(1{}-y^k_iy^k_j),&
    &\mlower{^{0k}_{\ i}}=\max\left(0, \sum_{j=1}^{n^k}  \ablue{^{0k}_{\ j}} Q(\gamma)^k_{ij}  -1 +y_i^k\right),&\\
    &\uu{^{0k}}=1,&
    &\mupper{^{0k}_{\ i}}=\max\left(0, -\sum_{j=1}^{n^k} \ablue{^{0k}_{\ j}} Q(\gamma)^k_{ij} +1 -y_i^k\right)&
\end{align}
for $k=1,\dots,K$, $i=1,\dots,n^k$. 
By construction,   $\theta^k_i(z^0)=0$ for $k=1,\dots,K$, $i=1,\dots,n^k$, and $\varphi^k(z^0)=0$ for $k=1,\dots,K$; see \eqref{eq:constraint_function_definitions}.
For a balanced dataset, where there is a similar number of positive and negative labels $y^k_i$ in each fold, the variables $\alpha^k_i$ are initialised close to $\frac{C}{2}$, the midpoint of its feasible interval $0\leq \alpha^k_i \leq C$, for $k=1,\dots,K$, $i=1,\dots,n^k$. Since one of $\mupper{}$ or $\mlower{}$ is also positive, the complementarity constraints are violated.  This means that the initialisation is not MPCC feasible for problem \MODEL{} but is feasible for the corresponding penalisation \penalisation{\penaltyparam{}^0} and relaxation \relaxation{\relaxparam{}^0} formulations for suitable choices of $\penaltyparam{}^0$ and $\relaxparam{}^0$.

Finally, we call the third strategy \texttt{LL-Feasible}, which corresponds, for a given choice of $(C^0, \gamma^0)$, to find  $\ablue{^0}$, $\mlower{^0}$, $\mupper{^0}$ and $\uu{^0}$ that satisfy the KKT conditions 
\ifsuplementaryincluded
\eqref{eq:lowerLevelKKT},
\else
\cite[p. 1891]{Fan2005}
\fi
characterising the lower-level problem. 
This is essentially the task of solving a mixed system of equations with complementarity terms, as studied by Michael Ferris and Steven Dirkse~\cite{Dirkse1995MCPLIB}.
In our experiments we will call the PATH solver \cite{misc_PATH_solver,dirkse1995,Dirkse1996} $K$ times, where $K$ is the number of folds, to solve each of the $K$ lower-level KKT systems. 
It is the most computationally expensive initialisation but ensures $z^0$ is feasible for all the programs \MODEL{}, \penalisation{\penaltyparam{}^0} and \relaxation{\relaxparam{}^0}. The three strategies can be summarised in Table~\ref{table:ll-init}.

\begin{table}[tbp]
    \centering
    \begin{tabular}{lcc}
        \toprule
        &  \shortstack[c]{Penalisation\\(Sequential)} &\shortstack[c]{Relaxation\\(Sequential)}  \\
        \hline
        Parameter $\penaltyparam{}$/$\relaxparam{}$ $Q_1$    & $10^{4}$      &  $10^{-6}$    \\
        Parameter $\penaltyparam{}$/$\relaxparam{}$ $Q_2$    & $10^{6}$      &  $10^{-6}$    \\
        Parameter $\penaltyparam{}$/$\relaxparam{}$ $Q_3$    & $10^{6}$      &  $10^{-6}$    \\
        MPCC Feasibility $Q_1$              & $3.94e-07$    &  $1.00e-06$   \\
        MPCC Feasibility $Q_2$              & $1.93e-08$    &  $9.98e-07$   \\
        MPCC Feasibility $Q_3$              & $4.88e-09$    &  $2.01e-07$   \\
        \bottomrule
    \end{tabular}
    \caption{This table presents the observed distributions of final parameter $\penaltyparam{}$ and $\relaxparam{}$ as well as the MPCC feasibility at termination. $Q_1, Q_2, Q_3$ denote the first quartile, median and third quartile.}
    \label{table:pi_and_feasibility}
\end{table}

\begin{figure}
    \centering
    \includegraphics[width=.5\linewidth]{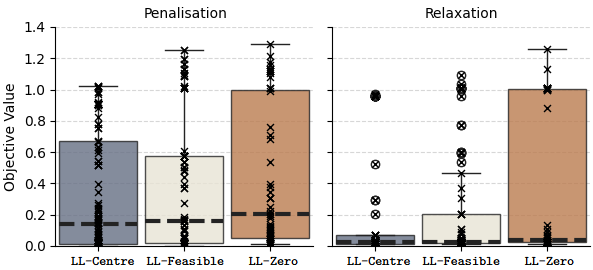}
    \caption{These box plots summarise the distribution of objective values achieved when using each of the three different lower-level initialisation strategies for the Circles dataset.}
    \label{figure:initialisation_strategies}
\end{figure}

When a feasible initialisation is made (i.e. \texttt{LL-Feasible}), at each index $i$, the complementarity choice (for example, choosing $G_i(z^0)=0$ or $H_i(z^0)=0$) is pre-decided.  We observe that, for many indices, this is unlikely to flip throughout the iterations of the solution methods.  However, when an infeasible initialisation (i.e. \texttt{LL-Centre}, \texttt{LL-Zero}) is chosen that violates the complementarity constraint (like having $G_i(z^0)=1$ and $H_i(z^0)=1$, for example), we observe that the same algorithms are more free to choose between $G_i(z^0)=0$ and $H_i(z^0)=0$ in the convergence process and arrive at better solutions. This phenomenon is demonstrated in the following experiment: 
For each dataset and solution method, we run 75 starting points, with 25 according to each of the strategies \texttt{LL-Zero}, \texttt{LL-Feasible} and \texttt{LL-Centre}.  We then assess the quality of the solution by the distribution of objective values found. These results for dataset \textit{Circles} are summarised as box plots in Figure~\ref{figure:initialisation_strategies}, where the effect of the lower-level initialisation is particularly pronounced for the relaxation method.  The algorithms are sensitive to the initial choice of $(C,\gamma)$. This is due to the highly non-convex nature of the hyperparameter space.  The choice of lower-level variable initialisation only has a very small effect on most datasets. However, overall, it is clear that \texttt{LL-Centre} encourages the algorithms to converge to higher quality solutions and have a low computational cost. Therefore, we will use the \texttt{LL-Centre} initialisation strategy for the remainder of our experiments.

\begin{figure}
\centering
\includegraphics[width=.6\linewidth]{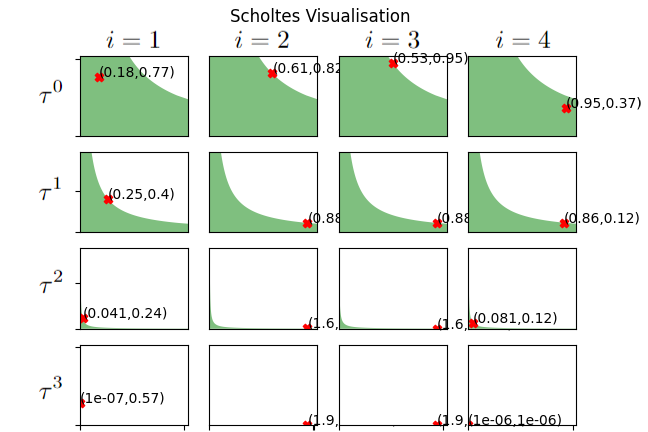}
\caption{The columns correspond to fold $k=1$ and indices $i=1, 2, 3, 4$ (Left to right). The rows correspond to sequentially smaller relaxation parameters $\relaxparam{}^t$ (top to bottom). The feasible regions for the relaxed program \relaxation{\relaxparam{}^t} are coloured in green. Each of the 16 small plots has a red point for the value of the complementarity variables $\ablue{_i}$ (horizontal axis) and {\color{dark_green}\underbar{$v$}${_i}$} (vertical axis) that were found to solve \relaxation{\relaxparam{}^t}.}
\label{fig:Scholtes}
\end{figure}

An interesting visualisation of the iteration of the sequential relaxation method solving problem \MODEL{}, parametrised by the Ionosphere dataset, is given in Figure~\ref{fig:Scholtes}.  This shows the values of the complementarity variables $\ablue{_i}$ and $\mlower{_i}$ for the first fold and first four indices $i=1,2,3,4$. Note that the variables are initially infeasible for the MPCC but feasible for the subproblem \eqref{eq:relaxation}. As we proceed down the rows, the variables $z$ converge to either the axis corresponding to $\ablue{_i}=0$ or the axis corresponding to $\mlower{_i}=0$, therefore generating an MPCC feasible final solution.

%% file: Supplementary/supplementary_52_datasets.tex
\subsection{Details of the classification datasets}
\label{SM:datasets}
We now present the data used for our experiments. We use eighteen datasets altogether.
Six of these (namely, \textit{Circles}, \textit{Hastie}, \textit{Linear}, \textit{Moons}, \textit{Swiss roll} and \textit{XOR}) are synthetically generated with code from scikit-learn \cite{Pedregosa2011}.  This is useful as globally and locally, optimal classifiers can be known from inspection, but could still be hard for an SVM to find.  This allows us to test the solution methods objectively.  The remaining twelve datasets come from real-world observations.  In these cases, the optimal values are unknown, and we must compare our solutions to those achieved by other authors in the literature. 
We perform standard preprocessing on these datasets before passing them to the solvers.  All rows with missing values are dropped.  For continuous features, we subtract the mean and divide by the standard deviation.  The resulting scaled data has a mean of zero and a standard deviation of one.  This helps avoid scaling issues that can cause our problem to become ill-conditioned.

Regarding categorical features, we either drop the column entirely or apply a one-hot encoding where $0-1$ indicator variables are added when they are relevant to the classification. 
Finally, we partition the multi-class labels into binary classes, which we replace with the labels -1 and +1.  A better classifier can sometimes be achieved through principal component analysis or other feature manipulation, however, we stick to these very simple preprocessing steps to allow a direct comparison to the original papers where the data was published.

Below, we list the datasets we use with a short description of any adjustments we make and a citation to the original publication. Each dataset is annotated with a tuple $(n,d)$ where $n$ is the number of observations and $d$ is the dimension of the feature vectors. 

\begin{enumerate}[leftmargin=*]
\setlength\itemsep{1em}

\item
\textbf{Banknote} $(1372, 4)$
The four features of this data set are digitally extracted from images of banknotes. Each banknote is then labelled as either genuine or forged. It is mostly balanced with 762 genuine and 610 forged examples. The banknote identification dataset was collected by Volker Lohweg \cite{banknote2013}.

\item 
\textbf{BUPA Liver Disorders} $(345, 6)$
Five features come from various blood tests relevant to liver disorders. A sixth feature comes from the patient's response to the question, ``How many half-pint equivalents of alcoholic beverages do you drink per day?''. If the patient is found to have any form of liver disorder, the patient is given a positive label. The data was published by Richard S. Forsyth and BUPA Medical Research Ltd. \cite{misc_BUPA_liver_disorders}

\item
\textbf{Circles} $(n, d)$
A synthetic dataset where half the examples are scattered uniformly over the surface of the $d$-dimensional unit ball and given a positive label. The other half are scattered on the surface of a smaller ball within the interior of the unit ball and given a negative label.  Gaussian noise with mean zero and standard deviation $\sigma^2$ is added to each point. For our experiments in Section~\ref{section:experiments}, we parametrise with $n=54, d=2, \sigma^2=0.2$.  See Figure~\ref{fig:kernel_trick} for an illustration.

\item
\textbf{Cleveland}
$(297, 14)$
The Cleveland Heart Disease data set was introduced by Detrano et al. \cite{Detrano1989}. Its features represent the physical attributes of a patient, such as age, sex, blood pressure, cholesterol, and heart rate, as recorded by their doctors. The labels represent a positive or negative diagnosis of heart disease. Originally, the dataset had 303 examples, 76 attributes and five classes. We adjust the dataset in the following manner. First, we select only  14 features that are commonly used by other authors. Then, we remove examples with missing values. Finally, we replace the different classes of heart disease with a positive label and the absence of heart disease with a negative label so that we can train a binary classifier.  

\item
\textbf{Glass}
$(138, 9)$
Each sample is described in terms of its oxide content (sodium, magnesium, aluminium, etc.) and refractive index. The samples are labelled by their glass type, which can be useful to determine for criminal investigations. This data was collected by B. German \cite{misc_glass_identification_42} and is openly available on the UCI Machine Learning Repository. Since it is a multi-class dataset, we make `building windows float processed' (70 examples) into the positive class and combine the vehicle, containers, tableware and lamp glass types (68 examples) into the negative class.

\item
\textbf{Hastie}
$(100, 10)$
This is the simulated data set from equation 10.2 of ``The Elements of Statistical Learning'' \cite{Hastie2009}. It has ten features that are independently and identically sampled according to a unit-independent Gaussian distribution. The labels are determined as $y_i=1$ if $\|x_i\|^2>9.34$ otherwise $y_i=-1$.

\item
\textbf{Ionosphere}
$(351, 34)$
This data was collected from radar antennas in the ionosphere in Goose Bay, Labrador \cite{misc_ionosphere_52} and used to classify structure in the ionosphere. It is somewhat unbalanced, with 225 positive labels and 126 negative labels.

\item
\textbf{Iris}
$(100, 4)$
The famous Iris flower data set was introduced by Anderson and Fisher \cite{Anderson1936,fisher1936}. Its four features represent the lengths and widths of the sepals and petals. The labels represent the species as either Setosa, Virginica or Versicolor. We discard all examples of the Versicolor class, so we end up with balanced binary labels of 50 Setosa and 50 Virginica.

\item
\textbf{Linear}
$(n, d)$
A simple synthetic dataset where the data points are scattered randomly and uniformly over the space $[0,1]^d$. A random hyperplane is chosen. The labels are then assigned as $+1$ if they lie above the hyperplane and $-1$ if they lie below. This is the only linearly separable dataset. In our experiments we choose $n=100, d=10$.

\item
\textbf{MNIST}
$(600, 784)$
Low resolution, 28x28 pixel, images of handwritten digits 0, 1, 2, ..., 9. The dataset is massive with 60,000 examples.  We subsample 300 examples of the digit 0 and 300 examples of the digit 1 in order to convert to a binary classification problem~\cite{LeCun1998}.

\item
\textbf{Mushrooms}
$(300, 117)$
This dataset details observations of gilled mushrooms in the Agaricus and Lepiota families \cite{misc_mushroom_73}.  All of the features are categorical, for example, the feature `cap shape' could take the value bell, conical, convex, flat, knobbed or sunken.  We add indicator variables for every category, which makes 117 features. Each example is labelled as either poisonous or edible.

\item
\textbf{Moons}
$(54, 2)$
A synthetic two-dimensional data set where points of the different classes are scattered in interweaving half-circles around each other, which can only be separated by a very nonlinear `S'-shaped boundary.

\item
\textbf{Palmer}
$(342, 5)$
Penguins had their features, such as the culmen length, culmen depth, flipper length, body mass and sex, recorded and were labelled by their species.  We perform a few modifications: drop the rows with missing values; replace the sex with a binary 1,0 indicator; and replace the species with a positive label if Gentoo and a negative label otherwise. The data was collected by Palmer Station Antarctica LTER and K. Gorman \cite{misc_palmer} and is openly available on GitHub.

\item
\textbf{Pima}
$(768, 8)$
The features of the Pima Indians Diabetes Database are measurements such as BMI, insulin level, and age from female patients of Pima Indian heritage. The labels represent the diagnosis of diabetes. It is quite unbalanced with 500 negative examples but only 268 positive. The dataset was provided by the National Institute of Diabetes and Digestive and Kidney Diseases and introduced to the literature by Smith et al. \cite{Smith1988}.  

\item
\textbf{Sonar}
$(208, 60)$
This data was collected to improve mine detection in the ocean.  The many features are those collected from the sonar signal. The labels are positive if the sonar signals bounced off a metal cylinder and negative if they bounced off rocks. Terry Sejnowski made the data available, and it was first published by Gorman and Sejnowski \cite{Gorman1988}.

\item
\textbf{Swiss Roll}
$(n, d)$
Swiss Roll is a common synthetic dataset that can be tricky for an SVM as the two classes are rolled inside one another. It is presented in the book Machine Learning: An Algorithmic Perspective \cite[Chapter 6, figure 6.13]{Marsland2011}. We use the Scikit-learn function make\_swiss\_roll to generate the data with Gaussian noise with standard deviation 0.1. In our experiments we choose $n=90, d=3$.



\item
\textbf{Wisconsin} $(569, 30)$
This first breast cancer dataset was presented in the work of Street et al. \cite{Street1993}. Its thirty features are detailed measurements of the breast and tumour, such as radius, concavity, and compactness. Each example is then labelled as malignant or benign.


\item
\textbf{XOR} $(n, 2)$
A two-dimensional synthetic data set where points are scattered uniformly at random across the plane $[-1,1]^2$. They are then labelled as positive if the product of their coordinates is positive and negative otherwise. This leads to the top right and bottom left quadrants being the opposite class to the top left and bottom right quadrants. It can be a particularly problematic dataset for classifiers. In our experiments we choose $n=100$.

\end{enumerate}
\vspace{1em}
% ...
These datasets are selected because they have quite different structures.  The optimal $(C,\gamma)$ hyperparameter choices are, for Ionosphere $(1.00e+01, 3.16e-02)$, Palmer ($1.00e+06, 1.00e-05$), Moons $(1.78e+01, 5.62e+00)$ and Wisconsin $(3.16e+03, 1.78e-04)$.  Furthermore, they show various different spreads of objective values and runtimes.

%% file: Supplementary/supplementary_53_derivative_free.tex
\subsection{Details of the derivative-free methods}
\label{SM:derivative_free}
% ...
Grid search \cite[p. 79]{Geron2022}\cite{Marsland2011} considers a sequence of plausible candidate choices for each hyperparameter and then forms a grid of every combination.
At each coordinate of the grid (representing a choice for each hyperparameter) the algorithm trains, evaluates and then discards $K$ fresh SVMs according to $K$-fold cross-validation.
This is repeated independently for every combination in the Cartesian product that constitutes the grid. 
Finally, the hyperparameter combination that results in the best mean validation accuracy across the $K$ classifiers is returned.
In our experiments we select a grid of size $10\times10$ where ten choices of $C$ are picked log-uniformly between $10^{-4}$ and $10^{6}$ and ten choices of $\gamma$ are picked log-uniformly between $10^{-5}$ and  $10^{4}$.

Random search \cite[p. 81]{Geron2022}\cite{Mantovani2015} proceeds similarly to grid search; however, rather than a systematic traversal of a grid, at each iteration the hyperparameters are drawn randomly and independently from predefined distributions.
This has the advantage of being able to search in-between the grid points of a grid search.
In our experiments we select 100 random points from a log-uniform distribution with the same upper and lower bounds as the grid search.

Bayesian optimisation \cite{garnett2023,Klein2017,Pelikan1999boa,Snoek2012} is a more sophisticated probabilistic method. 
It maintains a model approximation of the function $\tilde{f}$ that maps hyperparameter selections to their corresponding mean validation accuracy and an acquisition function representing the expected improvement in $\tilde{f}$.
At each iteration, the next choice of hyperparameters is determined by maximising the acquisition function; then the observed validation accuracy is used to update the approximation of $\tilde{f}$. 
The meta-hyperparameters of the surrogate model, including those controlling the balance between exploration and exploitation and the smoothness assumptions of $\tilde{f}$, must themselves be carefully chosen.
Bayesian optimisation has much more computationally expensive iterations, but is often able to find points near the global optimum.
We use the scikit-optimize implementation~\cite{Head2021scikit} with the following setup
\vskip.5em
\begin{verbatim}
skopt.gp_minimize(func,
  dimensions=[Real(1e-4, 1e+6, "log-uniform"), Real(1e-5, 1e+4, "log-uniform")],
  acq_func="gp_hedge", n_calls=100, random_state=42, xi=0.01, kappa=1.96,
)
\end{verbatim} 
\vskip.5em

\begin{figure}
    \centering
    \includegraphics[width=1\linewidth]{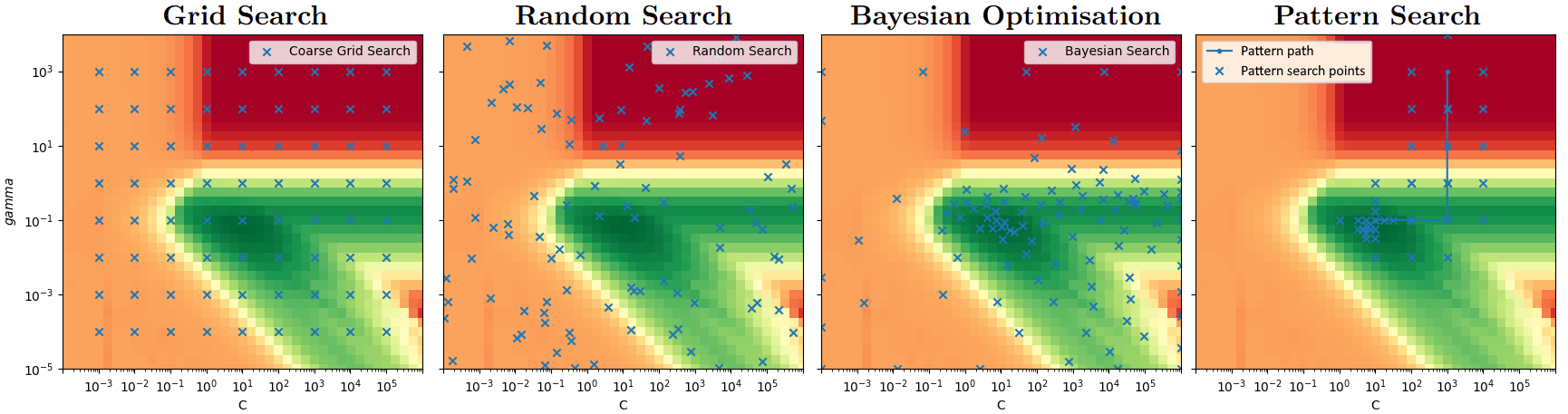}
    \caption{The heat map shows the hyperparameter space for the dataset Ionosphere.  Green regions correspond to $(C,\gamma)$ combinations, which allow the SVM to achieve high validation accuracy.  The heat map is repeated four times for the grid, random, Bayesian and pattern search methods.  On each, blue crosses are plotted on all the $(C,\gamma)$ combinations evaluated by that search method.}
    \label{fig:derivative_free}
\end{figure}

Pattern search \cite[Algorithm 9.2]{Nocedal1999} differs from the previous three DFO methods in that it is a local rather than a global search.  
The algorithm maintains a current best hyperparameter selection $(C,\gamma)$ and step size $h$.  
At each iteration, it evaluates new hyperparameters $\{(C\pm h,\gamma)\}\cup\{(C,\gamma\pm h)\}$ one step size away along each of the coordinate axes. 
If it finds a better objective value, it updates the current best selection and continues.  
Otherwise, it reduces its step size.  
Once the step size is sufficiently small, it terminates---at a local minimum if the function landscape is well-behaved.
This method can be made more sophisticated with pattern moves, momentum or a line search. 
Although it is less commonly used in machine learning, it provides an interesting comparison to the relaxation and penalisation algorithms, which are also local methods. 
In our experiments steps are made in log-space. 
The step size is initialised to $h=1$ and halved whenever a smaller objective cannot be found, down to $10^{-6}$.
All four methods are visualised in Figure~\ref{fig:derivative_free}.

Our MPCC optimisation method has two clear advantages over DFO. 
The first advantage comes from the idea that two choices of hyperparameters with similar values often result in similar SVM states.  
For example, an optimal SVM with $(C,\gamma)=(15, 0.3)$ is likely to have many of the same support vectors as the optimal SVM with $(C,\gamma)=(12, 0.2)$.  
This means only a few decision variables $\ablue{_i^k}$ for $(k,i)\in{I}$ must be changed to move between optima for the adjusted hyperparameters.  
When one of grid, random, Bayesian or pattern search moves between two nearby hyperparameter selections in this manner, it recalculates the optimal decision variables from scratch, whereas the MPCC optimisation method only adjusts the required variables.  
The second advantage comes from the fact that the MPCC optimisation approach makes use of the first and second-order derivatives of the objective function with respect to the hyperparameters to move in a descent direction through the hyperparameter space.


%% file: Supplementary/supplementary_54_results.tex
\begin{table}[ht]
\scriptsize
\centering
\begin{tabular}{ l l c c c c c c c c c}
    \toprule
&
& \vvhead{Penalisation}{(Exact)}
& \vvhead{Penalisation}{(Sequential)}
& \vvhead{Relaxation}{(Exact)}
& \vvhead{Relaxation}{(Sequential)}
& \vvhead{Semismooth}{Newton}
& \vhead{Grid}
& \vhead{Random}  
& \vhead{Bayesian}  
& \vhead{Pattern}  
\\
\hline
\multirow{9}{*}{\vhead{\textbf{Ionosphere}}}  
& Runtime $\ Q_1$ &\bf4    & 9      & 16     & 55     & 386    & 64     & 69     & 400    & 83    \\
& Runtime $\ Q_2$ &\bf11   & 12     & 29     & 90     & 589    & 98     & 106    & 880    & 132   \\
& Runtime $\ Q_3$ &\bf17   &\bf17   & 40     & 149    & 1369   & 135    & 136    & 2152   & 133   \\
& Objective $Q_1$ & 1.05   & 0.84   & 1.01   & 1.01   & 1.32   & 0.40   & 0.40   & 0.35   & 1.01  \\
& Objective $Q_2$ & 0.62   & 0.47   & 0.42   & 0.40   & 0.87   & 0.37   & 0.38   &\bf0.33 & 0.47  \\
& Objective $Q_3$ & 0.54   &\bf0.33 & 0.40   & 0.34   & 0.75   & 0.35   & 0.36   &\bf0.33 &\bf0.33\\
& Iterations      & 188    & 208    & 546    & 1569   & 891    & 49     & 50     & 50     & 8     \\
& Validation acc. & 84\%   &\bf89\% & 86\%   & 88\%   & 79\%   & 88\%   & 87\%  &\bf89\%  &\bf89\%\\
& Test acc.       & 83\%   &\bf89\% & 86\%   & 86\%   & 79\%   & 88\%   & 87\%  &\bf89\%  &\bf89\%\\
\hline
\multirow{9}{*}{\vhead{\textbf{Palmer}}}  
& Runtime $\ Q_1$ & \bf 1  & 106    & 111    & 200    & 579    & 71     & 73     & 389    & 105   \\
& Runtime $\ Q_2$ & \bf 1  & 133    & 168    & 616    & 798    & 108    & 111    & 968    & 160   \\
& Runtime $\ Q_3$ & \bf 27 & 201    & 569    & 808    & 1639   & 143    & 143    & 2200   & 206   \\
& Objective $Q_1$ & 0.89   & 0.69   & 0.91   & 0.42   & 1.26   & 0.20   & 0.04   & 0.04  & 0.90  \\
& Objective $Q_2$ & 0.80   & 0.36   & 0.19   & 0.19   & 0.63   & 0.09   & 0.04   &\bf0.03 & 0.90  \\
& Objective $Q_3$ & 0.44   &\bf0.03 & 0.19   & 0.19   & 0.46   & 0.04   & 0.04   &\bf0.03 & 0.04  \\
& Iterations      & 54     & 235    & 290    & 849    & 516    & 49     & 50     & 50     & 8     \\
& Validation acc. & 90\%   &\bf99\% & 95\%   & 95\%   & 90\%   &\bf99\% &\bf99\% &\bf99\% &\bf99\%\\
& Test acc.       & 90\%   &\bf99\% & 95\%   & 95\%   & 90\%   &\bf99\% &\bf99\% &\bf99\% &\bf99\%\\
\hline 
\multirow{9}{*}{\vhead{\textbf{Moons}}}  
& Runtime $\ Q_1$ &\bf1    & 2      & 2      & 6      & 3      & 49     & 49     & 447    & 72    \\
& Runtime $\ Q_2$ &\bf1    & 2      & 2      & 8      & 4      & 82     & 104    & 678    & 103   \\
& Runtime $\ Q_3$ &\bf1    & 3      & 3      & 16     & 112    & 104    & 111    & 979    & 217   \\
& Objective $Q_1$ & 0.95   & 0.92   & 0.49   & 0.49   & 0.76   & 0.35   & 0.35   & 0.32 & 0.81  \\
& Objective $Q_2$ & 0.60   & 0.38   & 0.48   & 0.32   & 0.54   & 0.34   & 0.34   &\bf0.31 & 0.47  \\
& Objective $Q_3$ & 0.56   &\bf0.31 & 0.34   &\bf0.31 & 0.48   & 0.32   & 0.33   &\bf0.31 & \bf0.31\\
& Iterations      & 65     & 335    & 348    & 1125   & 665    & 49     & 50     & 50     & 8     \\
& Validation acc. & 83\%   &\bf89\% & 88\%   &\bf89\% & 74\%   &\bf89\% &\bf89\% &\bf89\% &\bf89\% \\
& Test acc.       & 83\%   &   87\% & 87\%   &\bf89\% & 74\%   &\bf89\% &\bf89\% &   87\% &   87\% \\
\hline 
\multirow{9}{*}{\vhead{\textbf{Wisconsin}}}  
& Runtime $\ Q_1$ & 18     & 59     & 147    & 198    & 317    & 68     & 72     & 172    & 77    \\
& Runtime $\ Q_2$ & 27     & 70     & 247    & 467    & 547    & 82     & 85     & 210    & 97    \\
& Runtime $\ Q_3$ & 47     & 83     & 369    & 560    & 1431   & 94     & 95     & 591    & 104   \\
& Objective $Q_1$ & 0.76   & 0.09   & 0.40   & 0.17   & 0.96   & 0.10   & 0.09   &\bf0.08 & 1.00  \\
& Objective $Q_2$ & 0.52   &\bf0.08 & 0.18   & 0.13   & 0.51   & 0.09   &\bf0.08 &\bf0.08 & 0.12  \\
& Objective $Q_3$ & 0.27   &\bf0.08 & 0.18   & 0.12   & 0.39   &\bf0.08 &\bf0.08 &\bf0.08 &\bf0.08\\
& Iterations      & 85     & 230    & 852    & 1565   & 979    & 49     & 50     & 50     & 8     \\
& Validation acc. & 93\%   &\bf97\% & 94\%   & 96\%   & 90\%   &\bf97\% &\bf97\% &\bf97\% & 96\%  \\
& Test acc.       & 93\%   &\bf97\% & 94\%   & 96\%   & 90\%   &\bf97\% &\bf97\% &\bf97\% & 96\%  \\
\bottomrule
\end{tabular}
\caption{
A summary of the performance of the solution methods (columns) on four of the eighteen datasets (rows).  The first quartile ($Q_1$) is the runtime/objective value achieved by 75\% of initialisations; the second quartile ($Q_2$) is the median; and the third quartile ($Q_3$) is the runtime/objective value achieved by the top 25\% of initialisations.
}
\label{table:experiment_results}
\end{table}

%% file: references.bib
@article{Anderson1936,
  title     = {The species problem in {I}ris},
  author    = {Anderson, Edgar},
  journal   = {Annals of the Missouri Botanical Garden},
  volume    = {23},
  pages     = {457--509},
  year      = {1936},
  publisher = {St. Louis, Missouri Botanical Garden Press, 1914-},
  url       = {https://www.biodiversitylibrary.org/part/4079}
}

@article{Anitescu2000,
  author  = {Anitescu, Mihai},
  title   = {On Using the Elastic Mode in Nonlinear Programming Approaches to Mathematical Programs with Complementarity Constraints},
  journal = {{SIAM} Journal on Optimization},
  volume  = {15},
  number  = {4},
  pages   = {1203--1236},
  year    = {2005},
  note    = {Published version of Preprint ANL/MCS-P864-1200, Argonne National Laboratory, Argonne, IL, 2000},
  url     = {https://doi.org/10.1137/S1052623402401221}
}

@misc{banknote2013,
  author       = {Lohweg,Volker},
  title        = {Banknote Authentication},
  year         = {2013},
  howpublished = {{UCI} Machine Learning Repository},
  url          = {https://doi.org/10.24432/C55P57}
}

@article{Baumrucker2008,
  title     = {{MPEC} problem formulations and solution strategies with chemical engineering applications},
  author    = {Baumrucker, BT and Renfro, Jeffrey G and Biegler, Lorenz T},
  journal   = {Computers \& Chemical Engineering},
  volume    = {32},
  number    = {12},
  pages     = {2903--2913},
  year      = {2008},
  publisher = {Elsevier},
  url       = {https://doi.org/10.1016/j.compchemeng.2008.02.010}
}

@article{Bengio2000,
  title     = {Gradient-based optimization of hyperparameters},
  author    = {Bengio, Yoshua},
  journal   = {Neural Computation},
  volume    = {12},
  number    = {8},
  pages     = {1889--1900},
  year      = {2000},
  publisher = {{MIT} Press},
  url       = {https://doi.org/10.1162/089976600300015187}
}

@inproceedings{Bennett2006,
  title        = {Model Selection via Bilevel Optimization},
  author       = {Bennett, Kristin P and Hu, Jing and Ji, Xiaoyun and Kunapuli, Gautam and Pang, Jong-Shi},
  booktitle    = {The 2006 {IEEE} International Joint Conference on Neural Network Proceedings},
  year         = {2006},
  volume       = {},
  number       = {},
  pages        = {1922-1929},
  organization = {{IEEE}},
  url          = {https://doi.org/10.1109/IJCNN.2006.246935}
}

@inproceedings{Boser1992,
  title     = {A training algorithm for optimal margin classifiers},
  author    = {Boser, Bernhard E and Guyon, Isabelle M and Vapnik, Vladimir N},
  booktitle = {Proceedings of the fifth annual workshop on Computational learning theory},
  pages     = {144--152},
  year      = {1992},
  url       = {https://doi.org/10.1145/130385.130401}
}

@book{boyd2004,
  place     = {Cambridge},
  title     = {Convex Optimization},
  publisher = {Cambridge University Press},
  address   = {Cambridge},
  author    = {Boyd, Stephen and Vandenberghe, Lieven},
  year      = {2004},
  url       = {https://doi.org/10.1017/CBO9780511804441}
}

@incollection{Byrd2006,
  author    = {Byrd, Richard H. and Nocedal, Jorge and Waltz, Richard A.},
  editor    = {Di Pillo, G and Roma, M},
  title     = {{KNITRO}: {A}n Integrated Package for Nonlinear Optimization},
  booktitle = {Large-Scale Nonlinear Optimization},
  series    = {Nonconvex Optim. Appl.},
  year      = {2006},
  publisher = {Springer US},
  address   = {Boston, MA},
  volume    = {83},
  pages     = {35--59},
  isbn      = {978-0-387-30065-8},
  url       = {https://doi.org/10.1007/0-387-30065-1_4}
}

@article{Cervantes2020,
  title   = {A comprehensive survey on support vector machine classification: Applications, challenges and trends},
  journal = {Neurocomputing},
  volume  = {408},
  pages   = {189-215},
  year    = {2020},
  issn    = {0925-2312},
  author  = {Jair Cervantes and Farid Garcia-Lamont and Lisbeth Rodríguez-Mazahua and Asdrubal Lopez},
  url     = {https://doi.org/10.1016/j.neucom.2019.10.118}
}

@article{Chang2011,
  title     = {{LIBSVM}: a library for support vector machines},
  author    = {Chang, Chih-Chung and Lin, Chih-Jen},
  journal   = {ACM Transactions on Intelligent Systems and Technology (TIST)},
  volume    = {2},
  number    = {3},
  pages     = {1--27},
  year      = {2011},
  publisher = {Acm New York, NY, {USA}},
  url       = {https://www.csie.ntu.edu.tw/%7Ecjlin/papers/libsvm.pdf}
}

@article{Chapelle2001,
  author    = {Chapelle, Olivier and Vapnik, Vladimir and Bousquet, Olivier and Mukherjee, Sayan},
  title     = {Choosing Multiple Parameters for Support Vector Machines},
  year      = {2002},
  month     = {05},
  volume    = {46},
  pages     = {131--159},
  journal   = {Machine Learning},
  publisher = {Springer},
  url       = {https://doi.org/10.1023/A:1012450327387}
}

@article{Chen1995,
  author    = {Chen, Yang and Florian, Michael},
  title     = {The nonlinear bilevel programming problem: formulations, regularity and optimality conditions},
  journal   = {Optimization},
  volume    = {32},
  number    = {3},
  pages     = {193-209},
  year      = {1995},
  publisher = {Taylor \& Francis},
  url       = {https://doi.org/10.1080/02331939508844048}
}

@misc{Coniglio2022,
  author       = {Coniglio, Stefano and Dunn, Anthony and Li, Qingna and Zemkoho, Alain},
  title        = {Bilevel hyperparameter optimization for nonlinear support vector machines},
  year         = {2023},
  howpublished = {Optimization Online preprint},
  url          = {https://optimization-online.org/wp-content/uploads/2023/08/BHO-PaperV2.pdf},
  note         = {Accessed: 20/10/2026.}
}

@article{Cortes1995,
  title     = {Support-vector networks},
  author    = {Cortes, Corinna and Vapnik, Vladimir},
  journal   = {Machine Learning},
  volume    = {20},
  pages     = {273--297},
  year      = {1995},
  publisher = {Springer},
  url       = {https://doi.org/10.1007/BF00994018}
}

@article{Demiguel2005,
  title     = {A two-sided relaxation scheme for mathematical programs with equilibrium constraints},
  author    = {DeMiguel, Victor and Friedlander, Michael P. and Nogales, Francisco J. and Scholtes, Stefan},
  journal   = {{SIAM} Journal on Optimization},
  volume    = {16},
  number    = {2},
  pages     = {587--609},
  year      = {2005},
  publisher = {SIAM},
  url       = {https://doi.org/10.1137/04060754x}
}

@book{dempe2020bilevel,
  title     = {Bilevel Optimization: Advances and Next Challenges},
  editor    = {Dempe, Stephan and Zemkoho, Alain},
  series    = {Springer Optimization and Its Applications},
  volume    = {161},
  year      = {2020},
  publisher = {Springer},
  address   = {Cham},
  isbn      = {978-3-030-52119-6},
  url       = {https://link.springer.com/book/10.1007/978-3-030-52119-6}
}

@article{DempeDutta2012,
  title     = {Is bilevel programming a special case of a mathematical program with complementarity constraints?},
  author    = {Dempe, Stephan and Dutta, Joydeep},
  journal   = {Mathematical Programming},
  volume    = {131},
  pages     = {37--48},
  year      = {2012},
  publisher = {Springer},
  url       = {https://doi.org/10.1007/s10107-010-0342-1}
}

@article{DempeZemkoho2012,
  title    = {On the {K}arush–{K}uhn–{T}ucker reformulation of the bilevel optimization problem},
  journal  = {Nonlinear Analysis: Theory, Methods \& Applications},
  volume   = {75},
  number   = {3},
  pages    = {1202-1218},
  year     = {2012},
  issn     = {0362-546X},
  author   = {S. Dempe and A.B. Zemkoho},
  url      = {https://doi.org/10.1016/j.na.2011.05.097}
}

@article{Detrano1989,
  title   = {International application of a new probability algorithm for the diagnosis of coronary artery disease.},
  author  = {Robert C. Detrano and Andr{\'a}s J{\'a}nosi and Walter Steinbrunn and Matthias Emil Pfisterer and Johann-Jakob Schmid and Sarbjit Sandhu and Kern Guppy and Stella Lee and Victor Froelicher},
  journal = {American Journal of Cardiology},
  year    = {1989},
  volume  = {64},
  number  = {5},
  pages   = {304--310},
  url     = {https://doi.org/10.1016/0002-9149(89)90524-9}
}

@article{dirkse1995,
  author    = {Steven P. Dirkse and Michael C. Ferris},
  title     = {The {PATH} Solver: A Non-Monotone Stabilization Scheme for Mixed Complementarity Problems},
  journal   = {Optimization Methods and Software},
  publisher = {Informa {UK} Limited},
  year      = {1995},
  volume    = {5},
  number    = {2},
  pages     = {123--156},
  url       = {https://doi.org/10.1080/10556789508805606}
}

@article{Dirkse1995MCPLIB,
  title     = {{MCPLIB}: A collection of nonlinear mixed complementarity problems},
  author    = {Dirkse, Steven P and Ferris, Michael C},
  journal   = {Optimization Methods and Software},
  volume    = {5},
  number    = {4},
  pages     = {319--345},
  year      = {1995},
  publisher = {Taylor \& Francis},
  url       = {https://doi.org/10.1080/10556789508805619}
}

@article{Dirkse1996,
  title     = {A path search damped {N}ewton method for computing general equilibria},
  author    = {Dirkse, Steven P and Ferris, Michael C},
  journal   = {Annals of Operations Research},
  volume    = {68},
  number    = {2},
  pages     = {211--232},
  year      = {1996},
  publisher = {Springer},
  url       = {https://doi.org/10.1007/BF02209613}
}

@article{Dutta2013,
  title     = {Approximate {KKT} points and a proximity measure for termination},
  author    = {Dutta, Joydeep and Deb, Kalyanmoy and Tulshyan, Rupesh and Arora, Ramnik},
  journal   = {Journal of Global Optimization},
  volume    = {56},
  number    = {4},
  pages     = {1463--1499},
  year      = {2013},
  publisher = {Springer},
  url       = {https://doi.org/10.1007/s10898-012-9920-5}
}

@article{Fan2005,
  author  = {Rong-En Fan and Pai-Hsuen Chen and Chih-Jen Lin},
  title   = {Working Set Selection Using Second Order Information for Training Support Vector Machines},
  journal = {Journal of Machine Learning Research},
  year    = {2005},
  volume  = {6},
  number  = {63},
  pages   = {1889--1918},
  url     = {http://jmlr.org/papers/v6/fan05a.html}
}

@article{fisher1936,
  title     = {The use of multiple measurements in taxonomic problems},
  author    = {Fisher, Ronald A},
  journal   = {Annals of Eugenics},
  volume    = {7},
  number    = {2},
  pages     = {179--188},
  year      = {1936},
  publisher = {Wiley Online Library},
  url       = {https://doi.org/10.1111/j.1469-1809.1936.tb02137.x}
}

@phdthesis{Flegel2005,
  title  = {Constraint qualifications and stationarity concepts for mathematical programs with equilibrium constraints},
  author = {Flegel, Michael L},
  year   = {2005},
  school = {Universit{\"a}t W{\"u}rzburg},
  url    = {http://d-nb.info/975013661/34}
}

@techreport{fletcher2002numerical,
  title       = {Numerical experience with solving {MPEC}s as {NLP}s},
  author      = {Fletcher, Roger and Leyffer, Sven},
  type        = {Numerical Analysis Report},
  number      = {NA/210},
  institution = {Department of Mathematics, University of Dundee},
  address     = {Dundee, {UK}},
  year        = {2002},
  url     = {https://optimization-online.org/2002/08/522/}
}

@article{Fourer1990,
  title   = {{AMPL}: A mathematical programming language},
  author  = {Fourer, Robert and Gay, David M and Kernighan, Brian W},
  journal = {Management Science},
  volume  = {36},
  number  = {5},
  pages   = {519--554},
  year    = {1990},
  url     = {https://www.ampl.com/_archive/first-website/REFS/amplmod.pdf}
}

@inproceedings{Franceschi2018,
  title     = {Bilevel Programming for Hyperparameter Optimization and Meta-Learning},
  author    = {Franceschi, Luca and Frasconi, Paolo and Salzo, Saverio and Grazzi, Riccardo and Pontil, Massimiliano},
  booktitle = {Proceedings of the 35th International Conference on Machine Learning},
  pages     = {1568--1577},
  year      = {2018},
  editor    = {Dy, Jennifer and Krause, Andreas},
  volume    = {80},
  series    = {Proceedings of Machine Learning Research},
  month     = {10--15 Jul},
  publisher = {PMLR},
  address   = {Stockholm, Sweden},
  url       = {https://proceedings.mlr.press/v80/franceschi18a.html}
}

@article{Gabriel2010,
  title   = {Solving discretely-constrained {MPEC} problems with applications in electric power markets},
  author  = {Gabriel, Steven A and Leuthold, Florian U},
  journal = {Energy Economics},
  volume  = {32},
  number  = {1},
  pages   = {3-14},
  year    = {2010},
  issn    = {0140-9883},
  url     = {https://doi.org/10.1016/j.eneco.2009.03.008}
}

@inproceedings{Gao2022,
  title     = {Value Function based Difference-of-Convex Algorithm for Bilevel Hyperparameter Selection Problems},
  author    = {Gao, Lucy L and Ye, Jane and Yin, Haian and Zeng, Shangzhi and Zhang, Jin},
  booktitle = {Proceedings of the 39th International Conference on Machine Learning},
  pages     = {7164--7182},
  year      = {2022},
  editor    = {Chaudhuri, Kamalika and Jegelka, Stefanie and Song, Le and Szepesvari, Csaba and Niu, Gang and Sabato, Sivan},
  volume    = {162},
  series    = {Proceedings of Machine Learning Research},
  month     = {17--23 Jul},
  publisher = {PMLR},
  address   = {Baltimore, MD},
  url       = {https://proceedings.mlr.press/v162/gao22j.html}
}

@book{garnett2023,
  title     = {{B}ayesian optimization},
  author    = {Garnett, Roman},
  year      = {2023},
  publisher = {Cambridge University Press},
  address   = {Cambridge},
  url       = {https://bayesoptbook.com/}
}

@book{George2024,
  title     = {A Course in Linear Algebra},
  author    = {George, Raju K and Ajayakumar, Abhijith},
  year      = {2024},
  publisher = {Springer},
  address   = {Singapore},
  url       = {https://doi.org/10.1007/978-981-99-8680-4}
}

@book{Geron2022,
  title     = {Hands-on machine learning with Scikit-Learn, Keras, and TensorFlow},
  author    = {G{\'e}ron, Aur{\'e}lien},
  year      = {2022},
  publisher = {O'Reilly Media, Inc.},
  address   = {Sebastopol, CA},
  url       = {https://www.oreilly.com/library/view/hands-on-machine-learning/9781492032632/}
}

@article{gfrerer2017new,
  title     = {New constraint qualifications for mathematical programs with equilibrium constraints via variational analysis},
  author    = {Gfrerer, Helmut and Ye, Jane J},
  journal   = {{SIAM} Journal on Optimization},
  volume    = {27},
  number    = {2},
  pages     = {842--865},
  year      = {2017},
  publisher = {SIAM},
  url       = {https://doi.org/10.1137/16M1088752}
}

@article{Gorman1988,
  title     = {Analysis of hidden units in a layered network trained to classify sonar targets},
  author    = {Gorman, R Paul and Sejnowski, Terrence J},
  journal   = {Neural Networks},
  volume    = {1},
  number    = {1},
  pages     = {75-89},
  year      = {1988},
  publisher = {Elsevier},
  url       = {https://doi.org/10.1016/0893-6080(88)90023-8}
}

@inproceedings{Guyon1992,
  publisher = {Morgan-Kaufmann},
  address   = {San Mateo, CA},
  title     = {Automatic Capacity Tuning of Very Large {VC}-Dimension Classifiers},
  booktitle = {Advances in Neural Information Processing Systems},
  author    = {Guyon, Isabelle and Boser, B and Vapnik, Vladimir},
  journal   = {Advances in Neural Information Processing Systems},
  editor    = {S. Hanson and J. Cowan and C. Giles},
  pages     = {},
  volume    = {5},
  year      = {1992},
  url       = {https://proceedings.neurips.cc/paper_files/paper/1992/file/eaae339c4d89fc102edd9dbdb6a28915-Paper.pdf}
}

@article{Harder2021,
  author  = {Harder, Felix and Mehlitz, Patrick and Wachsmuth, Gerd},
  title   = {Reformulation of the M-Stationarity Conditions as a System of Discontinuous Equations and Its Solution by a Semismooth {N}ewton Method},
  journal = {{SIAM} Journal on Optimization},
  volume  = {31},
  number  = {2},
  pages   = {1459-1488},
  year    = {2021},
  url     = {https://doi.org/10.1137/20M1321413}
}

@article{Harris2020,
  title     = {Array programming with {NumPy}},
  author    = {Charles R. Harris and K. Jarrod Millman and St{\'{e}}fan J.
               van der Walt and Ralf Gommers and Pauli Virtanen and David
               Cournapeau and Eric Wieser and Julian Taylor and Sebastian
               Berg and Nathaniel J. Smith and Robert Kern and Matti Picus
               and Stephan Hoyer and Marten H. van Kerkwijk and Matthew
               Brett and Allan Haldane and Jaime Fern{\'{a}}ndez del
               R{\'{i}}o and Mark Wiebe and Pearu Peterson and Pierre
               G{\'{e}}rard-Marchant and Kevin Sheppard and Tyler Reddy and
               Warren Weckesser and Hameer Abbasi and Christoph Gohlke and
               Travis E. Oliphant},
  year      = {2020},
  month     = {sep},
  journal   = {Nature},
  volume    = {585},
  number    = {7825},
  pages     = {357--362},
  publisher = {Springer Science and Business Media {LLC}},
  url       = {https://doi.org/10.1038/s41586-020-2649-2}
}

@book{Hastie2009,
  title     = {The elements of statistical learning: data mining, inference, and prediction},
  author    = {Hastie, Trevor and Tibshirani, Robert and Friedman, Jerome H},
  edition   = {2nd},
  year      = {2009},
  publisher = {Springer},
  address   = {New York},
  url       = {https://doi.org/10.1007/978-0-387-84858-7}
}

@misc{Head2021scikit,
  author       = {Head, Tim and
                  Kumar, Manoj and
                  Nahrstaedt, Holger and
                  Louppe, Gilles and
                  Shcherbatyi, Iaroslav},
  title        = {scikit-optimize},
  year         = {2021},
  publisher    = {Zenodo},
  version      = {v0.9.0},
  url          = {https://doi.org/10.5281/zenodo.5565057}
}

@article{Hoheisel2013,
  title     = {Theoretical and numerical comparison of relaxation methods for mathematical programs with complementarity constraints},
  author    = {Hoheisel, Tim and Kanzow, Christian and Schwartz, Alexandra},
  journal   = {Mathematical Programming},
  volume    = {137},
  pages     = {257--288},
  year      = {2013},
  publisher = {Springer},
  url       = {https://doi.org/10.1007/s10107-011-0488-5}
}

@article{Hu2004,
  title   = {Convergence of a Penalty Method for Mathematical Programming with Complementarity Constraints},
  author  = {Hu, Xinmin M and Ralph, Daniel},
  journal = {Journal of Optimization Theory and Applications},
  volume  = {123},
  pages   = {365--390},
  year    = {2004},
  url     = {https://doi.org/10.1007/s10957-004-5154-0}
}

@article{Hunter2007,
  author    = {Hunter, J. D.},
  title     = {Matplotlib: A 2D graphics environment},
  journal   = {Computing in Science \& Engineering},
  year      = {2007},
  volume    = {9},
  number    = {3},
  pages     = {90--95},
  publisher = {{IEEE} COMPUTER SOC},
  url       = {https://doi.org/10.1109/MCSE.2007.55}
}

@article{Jiang2020,
  title     = {Hyper-parameter optimization for support vector machines using stochastic gradient descent and dual coordinate descent},
  author    = {Wei Jiang and Sauleh Siddiqui},
  journal   = {{EURO} Journal on Computational Optimization},
  volume    = {8},
  number    = {1},
  pages     = {85--101},
  year      = {2020},
  publisher = {Springer},
  url       = {https://doi.org/10.1007/s13675-019-00115-7}
}

@article{Kanzow2015,
  title     = {The price of inexactness: convergence properties of relaxation methods for mathematical programs with complementarity constraints revisited},
  author    = {Kanzow, Christian and Schwartz, Alexandra},
  journal   = {Mathematics of Operations Research},
  volume    = {40},
  number    = {2},
  pages     = {253--275},
  year      = {2015},
  publisher = {INFORMS},
  url       = {https://doi.org/10.1287/moor.2014.0667}
}

@article{Kenneth1954,
  issn      = {00129682, 14680262},
  author    = {Kenneth J. Arrow and Gerard Debreu},
  journal   = {Econometrica},
  number    = {3},
  pages     = {265--290},
  publisher = {[Wiley, Econometric Society]},
  title     = {Existence of an Equilibrium for a Competitive Economy},
  volume    = {22},
  year      = {1954},
  url       = {http://www.jstor.org/stable/1907353}
}

@article{KeerthiLin2003asymptotic,
  author  = {Keerthi, S. Sathiya and Lin, Chih-Jen},
  title   = {Asymptotic Behaviors of Support Vector Machines with {G}aussian Kernel},
  journal = {Neural Computation}, 
  volume = {15}, 
  number = {7},
  pages   = {1667--1689}, 
  year = {2003}, 
  month = {7},
  issn = {0899-7667},
  url = {https://doi.org/10.1162/089976603321891855}
}

@incollection{Kim2020,
  title     = {{MPEC} methods for bilevel optimization problems},
  author    = {Kim, Youngdae and Leyffer, Sven and Munson, Todd},
  editor    = {Dempe, Stephan
               and Zemkoho, Alain},
  booktitle = {Bilevel Optimization: Advances and Next Challenges},
  year      = {2020},
  publisher = {Springer International Publishing},
  address   = {Cham},
  pages     = {335--360},
  url       = {https://doi.org/10.1007/978-3-030-52119-6_12}
}

@inproceedings{Klein2017,
  title     = {Fast {B}ayesian optimization of machine learning hyperparameters on large datasets},
  author    = {Klein, Aaron and Falkner, Stefan and Bartels, Simon and Hennig, Philipp and Hutter, Frank},
  booktitle = {Proceedings of the 20th International Conference on Artificial Intelligence and Statistics},
  pages     = {528--536},
  year      = {2017},
  editor    = {Singh, Aarti and Zhu, Jerry},
  volume    = {54},
  series    = {Proceedings of Machine Learning Research},
  month     = {20--22 Apr},
  publisher = {PMLR},
  address   = {Fort Lauderdale, FL},
  url       = {https://proceedings.mlr.press/v54/klein17a.html}
}

@article{Kunapuli2008a,
  author    = {Kunapuli, G. and Bennett, K. P. and Hu, Jing and Pang, Jong-Shi},
  title     = {Classification model selection via bilevel programming},
  journal   = {Optimization Methods and Software},
  volume    = {23},
  number    = {4},
  pages     = {475-489},
  year      = {2008},
  publisher = {Taylor \& Francis},
  url       = {https://doi.org/10.1080/10556780802102586}
}

@inproceedings{Kunapuli2008b,
  author    = {Kunapuli, Gautam and Bennett, K and Hu, Jing and Pang, Jong-Shi},
  title     = {Bilevel model selection for support vector machines},
  booktitle = {{CRM} proceedings and lecture notes},
  volume    = {45},
  pages     = {129--158},
  year      = {2008},
  url       = {https://gkunapuli.github.io/files/08bilevelML.pdf}
}

@article{Lawphongpanich2004,
  title     = {An {MPEC} approach to second-best toll pricing},
  author    = {Lawphongpanich, Siriphong and Hearn, Donald W},
  journal   = {Mathematical Programming},
  volume    = {101},
  number    = {1},
  pages     = {33--55},
  year      = {2004},
  publisher = {Springer},
  url       = {https://doi.org/10.1007/s10107-004-0536-5}
}

@article{LeCun1998,
  title     = {Gradient-based learning applied to document recognition},
  author    = {LeCun, Yann and Bottou, L{\'e}on and Bengio, Yoshua and Haffner, Patrick},
  journal   = {Proceedings of the {IEEE}},
  volume    = {86},
  number    = {11},
  pages     = {2278--2324},
  year      = {1998},
  publisher = {Ieee},
  url       = {https://doi.org/10.1109/5.726791}
}

@incollection{Leyffer2006complementarity,
  title     = {Complementarity constraints as nonlinear equations: Theory and numerical experience},
  author    = {Leyffer, Sven},
  editor    = {Dempe, Stephan and Kalashnikov, Vyacheslav},
  booktitle = {Optimization with Multivalued Mappings: Theory, Applications, and Algorithms},
  year      = {2006},
  publisher = {Springer US},
  address   = {Boston, MA},
  pages     = {169--208},
  isbn      = {978-0-387-34221-4},
  url       = {https://doi.org/10.1007/0-387-34221-4_9}
}

@article{Leyffer2006interior,
  title     = {Interior methods for mathematical programs with complementarity constraints},
  author    = {Leyffer, Sven and L{\'o}pez-Calva, Gabriel and Nocedal, Jorge},
  journal   = {{SIAM} Journal on Optimization},
  volume    = {17},
  number    = {1},
  pages     = {52--77},
  year      = {2006},
  publisher = {SIAM},
  url       = {https://doi.org/10.1137/040621065}
}

@article{Li2022linear,
  title     = {Bilevel hyperparameter optimization for support vector classification: theoretical analysis and a solution method},
  author    = {Li, Qingna and Li, Zhen and Zemkoho, Alain},
  journal   = {Mathematical Methods of Operations Research},
  volume    = {96},
  number    = {3},
  pages     = {315--350},
  year      = {2022},
  publisher = {Springer},
  url       = {https://doi.org/10.1007/s00186-022-00798-6}
}

@inproceedings{Li2022unified,
  title        = {A Unified Framework and a Case Study for Hyperparameter Selection in Machine Learning via Bilevel Optimization},
  author       = {Li, Zhen and Qian, Yaru and Li, Qingna},
  booktitle    = {2022 5th International Conference on Data Science and Information Technology ({DSIT})},
  pages        = {1--8},
  year         = {2022},
  organization = {{IEEE}},
  url          = {http://doi.org/10.1109/DSIT55514.2022.9943929}
}

@article{Lin2005,
  title     = {A modified relaxation scheme for mathematical programs with complementarity constraints},
  author    = {Lin, Gui-Hua and Fukushima, Masao},
  journal   = {Annals of Operations Research},
  volume    = {133},
  number    = {1},
  pages     = {63--84},
  year      = {2005},
  publisher = {Springer},
  url       = {https://doi.org/10.1007/s10479-004-5024-z}
}

@book{Luo1996,
  title     = {Mathematical Programs with Equilibrium Constraints},
  publisher = {Cambridge University Press},
  author    = {Luo, Zhi-Quan and Pang, Jong-Shi and Ralph, Daniel},
  year      = {1996},
  address   = {Cambridge}
}

@misc{MacMPEC,
  author = {Sven Leyffer},
  title  = {{MacMPEC}},
  url    = {https://wiki.mcs.anl.gov/leyffer/index.php/MacMPEC}
}

@article{Mangasarian1969,
  title   = {Nonlinear Programming},
  author  = {Mangasarian, Olvi L.},
  journal = {McGraw-Hill, New York},
  year    = {1969},
  url     = {https://doi.org/10.1137/1.9781611971255}
}

@inproceedings{Mantovani2015,
  title        = {Effectiveness of random search in {SVM} hyper-parameter tuning},
  author       = {Mantovani, Rafael G and Rossi, Andr{\'e} LD and Vanschoren, Joaquin and Bischl, Bernd and De Carvalho, Andr{\'e} CPLF},
  booktitle    = {2015 International Joint Conference on Neural Networks (IJCNN)},
  pages        = {1--8},
  year         = {2015},
  organization = {Ieee},
  url          = {https://doi.org/10.1109/IJCNN.2015.7280664}
}

@book{Marsland2011,
  title     = {Machine learning: an algorithmic perspective},
  author    = {Marsland, Stephen},
  year      = {2011},
  publisher = {Chapman and Hall/CRC},
  address   = {Boca Raton, FL}
}

@misc{misc_BUPA_liver_disorders,
  author       = {Richard S. Forsyth},
  title        = {BUPA Liver Disorders},
  year         = {1990},
  howpublished = {{BUPA} Medical Research Ltd.},
  url          = {http://archive.ics.uci.edu/ml/machine-learning-databases/liver-disorders/bupa.data}
}

@misc{misc_glass_identification_42,
  author       = {German,B.},
  title        = {Glass Identification},
  year         = {1987},
  howpublished = {{UCI} Machine Learning Repository},
  url          = {https://doi.org/10.24432/C5WW2P}
}

@misc{misc_ionosphere_52,
  author       = {Sigillito, V. and Wing, S. and Hutton, L. and Baker, K.},
  title        = {Ionosphere},
  year         = {1989},
  howpublished = {{UCI} Machine Learning Repository},
  url          = {https://doi.org/10.24432/C5W01B}
}

@misc{misc_mushroom_73,
  author       = {G. H. Lincoff and Alfred A. Knopf},
  title        = {Mushroom},
  year         = {1981},
  howpublished = {The Audubon Society Field Guide to North American Mushrooms},
  url          = {https://doi.org/10.24432/C5959T}
}

@manual{misc_palmer,
  title  = {palmerpenguins: Palmer Archipelago ({A}ntarctica) penguin data},
  author = {Allison Marie Horst and Alison Presmanes Hill and Kristen B Gorman},
  year   = {2020},
  url    = {https://doi.org/10.5281/zenodo.3960218}
}

@misc{misc_PATH_solver,
  author = {Dirkse, Steven and Ferris, Michael C.  and Munson, Todd},
  title  = {The {PATH} Solver},
  url    = {https://pages.cs.wisc.edu/%7Eferris/path.html}
}

@inproceedings{Moore2009,
  title        = {Nonsmooth bilevel programming for hyperparameter selection},
  author       = {Moore, Gregory M and Bergeron, Charles and Bennett, Kristin P},
  booktitle    = {2009 {IEEE} International Conference on Data Mining Workshops},
  pages        = {374--381},
  year         = {2009},
  organization = {IEEE},
  url          = {https://doi.org/10.1109/ICDMW.2009.74}
}

@phdthesis{Moore2010,
  author   = {Moore, Gregory M.},
  year     = {2010},
  title    = {Bilevel programming algorithms for machine learning model selection},
  school   = {Rensselaer Polytechnic Institute},
  journal  = {ProQuest Dissertations and Theses},
  pages    = {107},
  isbn     = {978-1-124-54106-8},
  language = {English},
  url      = {https://www.proquest.com/dissertations-theses/bilevel-programming-algorithms-machine-learning/docview/860000550/se-2}
}

@article{Moore2011,
  title     = {Model selection for primal {SVM}},
  author    = {Moore, Gregory and Bergeron, Charles and Bennett, Kristin P},
  journal   = {Machine Learning},
  volume    = {85},
  pages     = {175--208},
  year      = {2011},
  publisher = {Springer},
  url       = {https://doi.org/10.1007/s10994-011-5246-7}
}

@article{Murray1971,
  title     = {Analytical expressions for the eigenvalues and eigenvectors of the {H}essian matrices of barrier and penalty functions},
  author    = {Murray, Walter},
  journal   = {Journal of Optimization Theory and Applications},
  volume    = {7},
  number    = {3},
  pages     = {189--196},
  year      = {1971},
  publisher = {Springer},
  url       = {https://doi.org/10.1007/BF00932477}
}

@incollection{Nicole1999,
  author    = {Nicole, Denis
               and Takeda, Kenji
               and Wolton, Ivan
               and Cox, Simon},
  title     = {Southampton High Performance Computing Centre},
  booktitle = {High-Performance Computing},
  year      = {1999},
  publisher = {Springer {US}},
  address   = {Boston, MA},
  pages     = {33--41},
  isbn      = {978-1-4615-4873-7},
  url       = {https://doi.org/10.1007/978-1-4615-4873-7_4}
}

@book{Nocedal1999,
  title     = {Numerical Optimization},
  author    = {Nocedal, Jorge and Wright, Stephen J},
  edition   = {2nd},
  year      = {2006},
  publisher = {Springer},
  address   = {New York},
  url       = {https://doi.org/10.1007/978-0-387-40065-5}
}

@article{Okuno2021,
  author  = {Takayuki Okuno and Akiko Takeda and Akihiro Kawana and Motokazu Watanabe},
  title   = {On lp-hyperparameter Learning via Bilevel Nonsmooth Optimization},
  journal = {Journal of Machine Learning Research},
  year    = {2021},
  volume  = {22},
  number  = {245},
  pages   = {1--47},
  url     = {http://jmlr.org/papers/v22/18-485.html}
}

@article{Pang1999,
  title     = {Complementarity constraint qualifications and simplified {B}-stationarity conditions for mathematical programs with equilibrium constraints},
  author    = {Pang, Jong-Shi and Fukushima, Masao},
  journal   = {Computational Optimization and Applications},
  volume    = {13},
  pages     = {111--136},
  year      = {1999},
  publisher = {Springer},
  url       = {https://doi.org/10.1023/A:1008656806889}
}

@article{Pedregosa2011,
  author  = {Fabian Pedregosa and Ga{{\"e}}l Varoquaux and Alexandre Gramfort and Vincent Michel and Bertrand Thirion and Olivier Grisel and Mathieu Blondel and Peter Prettenhofer and Ron Weiss and Vincent Dubourg and Jake Vanderplas and Alexandre Passos and David Cournapeau and Matthieu Brucher and Matthieu Perrot and {\'E}douard Duchesnay},
  title   = {Scikit-learn: Machine Learning in {P}ython},
  journal = {Journal of Machine Learning Research},
  year    = {2011},
  volume  = {12},
  number  = {85},
  pages   = {2825--2830},
  url     = {http://jmlr.org/papers/v12/pedregosa11a.html}
}

@inproceedings{Pelikan1999boa,
  title     = {{BOA}: The {B}ayesian optimization algorithm},
  author    = {Pelikan, Martin and Goldberg, David E. and Cant\'{u}-Paz, Erick},
  year      = {1999},
  isbn      = {1558606114},
  publisher = {Morgan Kaufmann Publishers Inc.},
  address   = {San Francisco, CA, USA},
  booktitle = {Proceedings of the 1st Annual Conference on Genetic and Evolutionary Computation - Volume 1},
  pages     = {525–532},
  numpages  = {8},
  location  = {Orlando, Florida},
  series    = {GECCO'99},
  url       = {https://dl.acm.org/doi/10.5555/2933923.2933973}
}

@book{Python2009,
  author    = {Van Rossum, Guido and Drake, Fred L.},
  title     = {Python 3 Reference Manual},
  year      = {2009},
  isbn      = {1441412697},
  publisher = {CreateSpace},
  address   = {Scotts Valley, CA},
  url       = {https://dl.acm.org/doi/book/10.5555/1593511}
}

@article{Qi2000,
  title     = {On the constant positive linear dependence condition and its application to {SQP} methods},
  author    = {Qi, Liqun and Wei, Zengxin},
  journal   = {{SIAM} Journal on Optimization},
  volume    = {10},
  number    = {4},
  pages     = {963--981},
  year      = {2000},
  publisher = {SIAM},
  url       = {https://doi.org/10.1137/S1052623497326629}
}

@misc{Qiu2015,
  title         = {Solving Mathematical Programs with Equilibrium Constraints as Nonlinear Programming: A New Framework},
  author        = {Songqiang Qiu and Zhongwen Chen},
  year          = {2015},
  eprint        = {1510.07145},
  archiveprefix = {arXiv},
  primaryclass  = {math.OC},
  url           = {https://arxiv.org/abs/1510.07145}
}

@article{Rahman2018,
  title     = {{isGPT}: An optimized model to identify sub-Golgi protein types using {SVM} and Random Forest based feature selection},
  author    = {Rahman, M Saifur and Rahman, Md Khaledur and Kaykobad, M and Rahman, M Sohel},
  journal   = {Artificial intelligence in medicine},
  volume    = {84},
  pages     = {90--100},
  year      = {2018},
  publisher = {Elsevier},
  url       = {https://doi.org/10.1016/j.artmed.2017.11.003}
}

@article{Ralph2004,
  title     = {Some properties of regularization and penalization schemes for {MPEC}s},
  author    = {Ralph, Daniel and Wright, Stephen J},
  journal   = {Optimization Methods and Software},
  volume    = {19},
  number    = {5},
  pages     = {527--556},
  year      = {2004},
  publisher = {Taylor \& Francis},
  url       = {https://doi.org/10.1080/10556780410001709439}
}

@article{Scheel2000,
  title     = {Mathematical programs with complementarity constraints: Stationarity, optimality, and sensitivity},
  author    = {Scheel, Holger and Scholtes, Stefan},
  journal   = {Mathematics of Operations Research},
  volume    = {25},
  number    = {1},
  pages     = {1--22},
  year      = {2000},
  publisher = {INFORMS},
  url       = {https://doi.org/10.1287/moor.25.1.1.15213}
}

@article{Scholtes2001,
  title     = {Convergence properties of a regularization scheme for mathematical programs with complementarity constraints},
  author    = {Scholtes, Stefan},
  journal   = {{SIAM} Journal on Optimization},
  volume    = {11},
  number    = {4},
  pages     = {918--936},
  year      = {2001},
  publisher = {SIAM},
  url       = {https://doi.org/10.1137/S1052623499361233}
}

@phdthesis{schwartz2011mathematical,
  title  = {Mathematical programs with complementarity constraints: Theory, methods and applications},
  author = {Schwartz, Alexandra},
  year   = {2011},
  school = {Universit{\"a}t W{\"u}rzburg},
  url    = {https://tu-dresden.de/mn/math/numerik/schwartz/ressourcen/dateien/Schwartz_Dissertation.pdf}
}

@inproceedings{Smith1988,
  title        = {Using the {ADAP} learning algorithm to forecast the onset of diabetes mellitus},
  author       = {Smith, Jack W and Everhart, James E and Dickson, WC and Knowler, William C and Johannes, Robert Scott},
  booktitle    = {Proceedings of the annual symposium on computer application in medical care},
  pages        = {261},
  year         = {1988},
  organization = {American Medical Informatics Association},
  url          = {https://pmc.ncbi.nlm.nih.gov/articles/PMC2245318/}
}

@inproceedings{Smith2001,
  author    = {Smith, N. and Gales, Mark},
  booktitle = {Advances in Neural Information Processing Systems},
  editor    = {T. Dietterich and S. Becker and Z. Ghahramani},
  pages     = {},
  publisher = {MIT Press},
  address   = {Cambridge, MA},
  title     = {Speech Recognition using {SVM}s},
  volume    = {14},
  year      = {2001},
}

@inproceedings{Snoek2012,
  title     = {Practical {B}ayesian Optimization of Machine Learning Algorithms},
  author    = {Snoek, Jasper and Larochelle, Hugo and Adams, Ryan P},
  booktitle = {Advances in Neural Information Processing Systems},
  editor    = {F. Pereira and C.J. Burges and L. Bottou and K.Q. Weinberger},
  pages     = {},
  publisher = {Curran Associates, Inc.},
  address   = {Red Hook, NY},
  volume    = {25},
  year      = {2012},
  url       = {https://proceedings.neurips.cc/paper_files/paper/2012/file/05311655a15b75fab86956663e1819cd-Paper.pdf}
}

@incollection{Solodov2010,
  title     = {Constraint Qualifications},
  author    = {Solodov, Mikhail V.},
  year      = {2011},
  booktitle = {Wiley Encyclopedia of Operations Research and Management Science},
  publisher = {John Wiley \& Sons},
  address   = {Hoboken, NJ},
  isbn      = {9780470400531},
  chapter   = {4},
  pages     = {341},
  url       = {https://doi.org/10.1002/9780470400531.eorms0978}
}

@book{Stackelberg,
  year      = 1934,
  publisher = {Springer},
  address   = {Vienna},
  title     = {Marktform und {G}leichgewicht},
  author    = {Stackelberg, H. von}
}

@article{Steffensen2010,
  title     = {A new relaxation scheme for mathematical programs with equilibrium constraints},
  author    = {Steffensen, Sonja and Ulbrich, Michael},
  journal   = {{SIAM} Journal on Optimization},
  volume    = {20},
  number    = {5},
  pages     = {2504--2539},
  year      = {2010},
  publisher = {SIAM},
  url       = {https://doi.org/10.1137/090748883}
}

@article{Stone1974,
  title     = {Cross-validatory choice and assessment of statistical predictions},
  author    = {Stone, Mervyn},
  journal   = {Journal of the royal statistical society: Series {B} (Methodological)},
  volume    = {36},
  number    = {2},
  pages     = {111--133},
  year      = {1974},
  publisher = {Wiley Online Library},
  url       = {https://doi.org/10.1111/j.2517-6161.1974.tb00994.x}
}

@inproceedings{Street1993,
  title        = {Nuclear feature extraction for breast tumor diagnosis},
  author       = {William Nick Street and William H. Wolberg and Olvi L. Mangasarian},
  volume       = {1905},
  booktitle    = {Biomedical Image Processing and Biomedical Visualization},
  editor       = {Raj S. Acharya and Dmitry B. Goldgof},
  organization = {International Society for Optics and Photonics},
  publisher    = {SPIE},
  address      = {Bellingham, WA},
  pages        = {861 -- 870},
  year         = {1993},
  doi          = {10.1117/12.148698},
  url          = {https://doi.org/10.1117/12.148698}
}

@article{Tefas2001,
  title     = {Using support vector machines to enhance the performance of elastic graph matching for frontal face authentication},
  author    = {Tefas, Anastasios and Kotropoulos, Constantine and Pitas, Ioannis},
  journal   = {{IEEE} Transactions on Pattern Analysis and Machine Intelligence},
  volume    = {23},
  number    = {7},
  pages     = {735--746},
  year      = {2001},
  publisher = {IEEE},
  url       = {https://doi.org/10.1109/34.935847}
}

@book{Vapnik1999,
  title     = {The nature of statistical learning theory},
  author    = {Vapnik, Vladimir},
  year      = {1999},
  publisher = {Springer science \& business media},
  address   = {New York},
  url       = {https://doi.org/10.1007/978-1-4757-3264-1}
}

@article{Wachter2006,
  title     = {On the implementation of an interior-point filter line-search algorithm for large-scale nonlinear programming},
  author    = {W{\"a}chter, Andreas and Biegler, Lorenz T},
  journal   = {Mathematical Programming},
  volume    = {106},
  pages     = {25--57},
  year      = {2006},
  publisher = {Springer},
  url       = {https://doi.org/10.1007/s10107-004-0559-y}
}

@article{Wainer2021,
  title     = {How to tune the {RBF} {SVM} hyperparameters? An empirical evaluation of 18 search algorithms},
  author    = {Wainer, Jacques and Fonseca, Pablo},
  journal   = {Artificial Intelligence Review},
  volume    = {54},
  number    = {6},
  pages     = {4771--4797},
  year      = {2021},
  publisher = {Springer},
  url       = {https://doi.org/10.1007/s10462-021-10011-5}
}

@article{Ward2025,
  title   = {Mathematical programs with complementarity constraints and application to hyperparameter tuning},
  author  = {Ward, Samuel and Zemkoho, Alain and Ahipasaoglu, Selin},
  journal = {arXiv preprint arXiv:2504.13006},
  year    = {2025},
  url={https://doi.org/10.48550/arXiv.2504.13006}
}

@article{Waskom2021,
  year      = {2021},
  publisher = {The Open Journal},
  volume    = {6},
  number    = {60},
  pages     = {3021},
  author    = {Michael L. Waskom},
  title     = {seaborn: statistical data visualization},
  journal   = {Journal of Open Source Software},
  url       = {https://doi.org/10.21105/joss.03021}
}

@book{Wendland2004,
  title     = {Scattered data approximation},
  author    = {Wendland, Holger},
  volume    = {17},
  year      = {2004},
  publisher = {Cambridge University Press},
  address   = {Cambridge},
  url       = {https://doi.org/10.1017/CBO9780511617539}
}

@article{Wolpert1996,
  author  = {Wolpert, David H.},
  title   = {The Lack of A Priori Distinctions Between Learning Algorithms},
  journal = {Neural Computation},
  volume  = {8},
  number  = {7},
  pages   = {1341-1390},
  year    = {1996},
  month   = {10},
  issn    = {0899-7667},
  url     = {https://doi.org/10.1162/neco.1996.8.7.1341}
}

@article{Wright1994,
  title     = {Some properties of the {H}essian of the logarithmic barrier function},
  author    = {Wright, Margaret H},
  journal   = {Mathematical Programming},
  volume    = {67},
  number    = {1},
  pages     = {265--295},
  year      = {1994},
  publisher = {Springer},
  url       = {https://doi.org/10.1007/BF01582224}
}

@article{Ye2005,
  title   = {Necessary and sufficient optimality conditions for mathematical programs with equilibrium constraints},
  journal = {Journal of Mathematical Analysis and Applications},
  volume  = {307},
  number  = {1},
  pages   = {350-369},
  year    = {2005},
  issn    = {0022-247X},
  author  = {Jane J. Ye},
  url     = {https://doi.org/10.1016/j.jmaa.2004.10.032}
}

@book{Zhang2006,
  title     = {The Schur complement and its applications},
  author    = {Zhang, Fuzhen},
  volume    = {4},
  year      = {2006},
  publisher = {Springer Science \& Business Media},
  address   = {New York},
  url       = {https://doi.org/10.1007/b105056}
}
